\documentclass[a4paper,12pt]{article}

\usepackage{amsmath}
\usepackage[backend=biber,style=numeric-comp,bibstyle=numeric,sorting=none,url=false,natbib=true,maxbibnames=99,firstinits=true]{biblatex}

\usepackage{setspace}

\usepackage{hyperref}
\usepackage[a4paper,left=2.5cm,right=2.5cm,top=3cm,bottom=2.5cm]{geometry}
\usepackage{authblk}
\usepackage{amsmath,bm}
\usepackage{graphicx}
\usepackage{mathrsfs}
\usepackage{amssymb}
\usepackage{physics}
\usepackage{tensor}
\usepackage{xcolor}
\usepackage{tikz}
\usepackage{pgfplots}
\usetikzlibrary{spy}

\pgfplotsset{grid style={line width=.1pt,dashed,lightgray},compat=newest}
\usepackage{caption}
\usepackage{subcaption}
\usepackage{float}

\usepackage[font=small,labelfont=bf]{caption}

\DeclareCaptionLabelFormat{subfigure}{#2)}
\usepackage{amsthm}
\usepackage{booktabs,dcolumn}
\usepackage{multirow}
\usepackage[english]{babel}
\usepackage{soul}
\usepackage{mathtools}

\usepackage[title]{appendix}

\usepackage{xstring}

\newtheoremstyle{rem}
{6pt}
{6pt}
{\small}
{}
{\bf}
{:}
{.5em}
{}

\theoremstyle{rem}
\newtheorem*{remark}{Remark} 

\title{A computational model for short-range van der Waals interactions between beams and shells}

\author[1,2]{A.~Borković}
\author[1]{M.~H.~Gfrerer}
\author[3,4,5]{R.~A.~Sauer}

\affil[1]{Institute of Applied Mechanics, Graz University of Technology, Technikerstraße 4/II, 8010 Graz, Austria, aleksandar.borkovic@aggf.unibl.org, aborkovic@tugraz.at}
\affil[2]{University of Banja Luka, Faculty of Architecture, Civil Engineering and Geodesy, Department of Mechanics and Theory of Structures, 78000 Banja Luka, Bosnia and Herzegovina}
\affil[3]{Institute for Structural Mechanics, Ruhr University Bochum, Universitätsstraße 150, 44801 Bochum }
\affil[4]{Department of Structural Mechanics, Gdansk University of Technology, ul. Narutowicza 11/12, 80-233 Gdansk, Poland} 
\affil[5]{Dept. of Mechanical Engineering, Indian Institute of Technology Guwahati, Assam 781039, India }
\date{}                     
\begin{document}
	
\newcommand{\ssub}[2]{{#1}_{#2}} 
\newcommand{\vsub}[2]{\textbf{#1}_{#2}} 
\newcommand{\ssup}[2]{{#1}^{#2}} 
\newcommand{\vsup}[2]{\textbf{#1}^{#2}} 
\newcommand{\ssupsub}[3]{{#1}^{#2}_{#3}} 
\newcommand{\vsupsub}[3]{\textbf{#1}^{#2}_{#3}} 

\newcommand{\veq}[1]{\bar{\boldsymbol{#1}}} 
\newcommand{\veqn}[1]{\bar{\textbf{#1}}} 
\newcommand{\seq}[1]{\bar{#1}} 
\newcommand{\ve}[1]{\boldsymbol{#1}} 
\newcommand{\veh}[1]{\hat{\boldsymbol{#1}}} 
\newcommand{\ven}[1]{\textbf{#1}} 
\newcommand{\vepre}[1]{\boldsymbol{#1}^\sharp} %
\newcommand{\sdef}[1]{#1^*} 
\newcommand{\vdef}[1]{{\boldsymbol{#1}}^*} 
\newcommand{\vdefeq}[1]{{\bar{\boldsymbol{#1}}}^*} 
\newcommand{\trans}[1]{\boldsymbol{#1}^\mathsf{T}} 
\newcommand{\transn}[1]{\textbf{#1}^\mathsf{T}} 
\newcommand{\transmd}[1]{\dot{\boldsymbol{#1}}^\mathsf{T}} 
\newcommand{\mdvdef}[1]{\dot{\boldsymbol{#1}}^*} 
\newcommand{\mdsdef}[1]{\dot{#1}^*} 
\newcommand{\mdv}[1]{\dot{\bm{#1}}} 
\newcommand{\mdvni}[1]{\dot{\textbf{#1}}} 
\newcommand{\mds}[1]{\dot{#1}} 

\newcommand{\loc}[1]{\hat{#1}} 
\newcommand{\iloc}[3]{\hat{#1}^{#2}_{#3}} 
\newcommand{\ilocmd}[3]{\dot{\hat{#1}}^{#2}_{#3}} 
\newcommand{\md}[1]{\dot{#1}} 

\newcommand{\ii}[3]{{#1}^{#2}_{#3}} 
\newcommand{\iv}[3]{\boldsymbol{#1}^{#2}_{#3}} 
\newcommand{\ivh}[3]{\hat{\boldsymbol{#1}}^{#2}_{#3}} 
\newcommand{\ivn}[3]{\boldsymbol{#1}^\text{#2}_\text{#3}} 
\newcommand{\ivnh}[3]{\hat{\boldsymbol{#1}}^\text{#2}_\text{#3}} 
\newcommand{\idef}[3]{{#1}^{* #2}_{#3}} 
\newcommand{\ivdef}[3]{\boldsymbol{#1}^{* #2}_{#3}} 
\newcommand{\ipre}[3]{{#1}^{\sharp #2}_{#3}} 
\newcommand{\ivpre}[3]{\textbf{#1}^{\sharp #2}_{#3}} 

\newcommand{\ieq}[3]{\bar{#1}^{#2}_{#3}} 
\newcommand{\ic}[3]{\tilde{#1}^{#2}_{#3}} 
\newcommand{\icdef}[3]{\tilde{#1}^{* #2}_{#3}} 
\newcommand{\icpre}[3]{\tilde{#1}^{\sharp #2}_{#3}} 
\newcommand{\iveq}[3]{\bar{\boldsymbol{#1}}^{#2}_{#3}} 
\newcommand{\iveqn}[3]{\bar{\boldsymbol{#1}}^\text{#2}_\text{#3}} 
\newcommand{\ieqdef}[3]{\bar{#1}^{* #2}_{#3}} 
\newcommand{\iveqdef}[3]{\bar{\textbf{#1}}^{* #2}_{#3}} 
\newcommand{\ieqmddef}[3]{\dot{\bar{#1}}^{* #2}_{#3}} 
\newcommand{\icmddef}[3]{\dot{\tilde{#1}}^{* #2}_{#3}} 
\newcommand{\iveqmddef}[3]{\dot{\bar{\textbf{#1}}}^{* #2}_{#3}} 
\newcommand{\iveqmdddef}[3]{\ddot{\bar{\textbf{#1}}}^{* #2}_{#3}} 

\newcommand{\ieqpre}[3]{\bar{#1}^{\sharp #2}_{#3}} 
\newcommand{\iveqpre}[3]{\bar{\textbf{#1}}^{\sharp #2}_{#3}} 
\newcommand{\ieqmdpre}[3]{\dot{\bar{#1}}^{\sharp #2}_{#3}} 
\newcommand{\icmdpre}[3]{\dot{\tilde{#1}}^{\sharp #2}_{#3}} 
\newcommand{\iveqmdpre}[3]{\dot{\bar{\textbf{#1}}}^{\sharp #2}_{#3}} 

\newcommand{\ieqmd}[3]{\dot{\bar{#1}}^{#2}_{#3}} 
\newcommand{\icmd}[3]{\dot{\tilde{#1}}^{#2}_{#3}} 
\newcommand{\iveqmd}[3]{\dot{\bar{\boldsymbol{#1}}}^{#2}_{#3}} 
\newcommand{\iveqmdd}[3]{\ddot{\bar{\boldsymbol{#1}}}^{#2}_{#3}} 

\newcommand{\imddef}[3]{\dot{#1}^{* #2}_{#3}} 
\newcommand{\ivmddef}[3]{\dot{\textbf{#1}}^{* #2}_{#3}} 
\newcommand{\ivmdddef}[3]{\ddot{\textbf{#1}}^{* #2}_{#3}} 

\newcommand{\imdpre}[3]{\dot{#1}^{\sharp #2}_{#3}} 
\newcommand{\ivmdpre}[3]{\dot{\textbf{#1}}^{\sharp #2}_{#3}} 

\newcommand{\imd}[3]{\dot{#1}^{#2}_{#3}} 
\newcommand{\imdd}[3]{\ddot{#1}^{#2}_{#3}} 

\newcommand{\ivmd}[3]{\dot{\textbf{#1}}^{#2}_{#3}} 

\newcommand{\iii}[5]{^{#2}_{#3}{#1}^{#4}_{#5}} 
\newcommand{\iiv}[5]{^{#2}_{#3}{\boldsymbol{#1}}^{#4}_{#5}} 
\newcommand{\iivn}[5]{^{#2}_{#3}{\tilde{\boldsymbol{#1}}}^{#4}_{#5}} 
\newcommand{\iiieq}[5]{^{#2}_{#3}{\bar{#1}}^{#4}_{#5}} 
\newcommand{\iiieqt}[5]{^{#2}_{#3}{\tilde{#1}}^{#4}_{#5}} 

\newcommand{\eqqref}[1]{Eq.~\eqref{#1}} 
\newcommand{\fref}[1]{Fig.~\ref{#1}} 
	
\maketitle
	
\section*{Abstract}

We consider potential-based interactions between beams (or fibers) and shells (or membranes) using a coarse-grained approach with focus on van der Waals attraction and steric repulsion. The involved 6D integral over volumes of a beam and a shell is split into a 5D analytical pre-integration over the beam's cross section and a surrogate plate tangential to the closest point on the shell, and the remaining 1D numerical integration along the beam's axis. This general inverse-power interaction potential is added to the potential energies of the Bernoulli-Euler beam and the Kirchhoff-Love shell. The total potential energy is spatially discretized using isogeometric finite elements, and the nonlinear weak form of quasi-static equilibrium is solved using the continuation method. We provide error estimates and convergence analysis, together with two intriguing numerical examples. The developed approach provides excellent balance between accuracy and efficiency for small separations.


\textbf{Keywords}: interaction potential; van der Waals attraction; coarse-grained approach; contact mechanics; fiber-membrane interaction; beam-shell contact

\section{Introduction}

Intermolecular forces are the underlying cause of many observable phenomena, such as friction, adhesion, surface tension, and mechanical contact. Among them, van der Waals (vdW) attraction is arguably the most universal because it exists between uncharged molecules due to transient fluctuations in charge distributions. This force remains significant across both short and long distances, causing not only the cohesion of condensed matter but also an adhesion between macroscopic bodies. Modeling vdW interactions is challenging, as it is influenced by multiple contributing factors, retardation effects, and exhibits non-additive behavior \cite{2005parsegian}. A widely used approximation for estimating vdW forces between bodies is the pairwise summation method, which treats the total interaction as the sum of all individual point-point interactions. The reliability of this approximation depends strongly on the molecular geometry and type \cite{2017venkataram, 2017ambrosetti}. In this study, we assume that the pairwise summation approach is valid and that the vdW potential can be described using an inverse sixth-power law w.r.t.~the point-pair distance.

Numerical modeling of intermolecular interactions is commonly performed using molecular dynamics or Monte Carlo simulations \cite{2011israelachvili}. 
An alternative to these approaches is a \emph{coarse-grained} model based on the homogenization and coarse-graining of the underlying molecular system. This method strikes a favorable balance between accuracy and computational efficiency \cite{1997argento, 2007sauer, 2008sauer, 2016fan} by combining the underlying physics of molecular interactions with the efficiency of continuum-based contact formulations \cite{2006wriggers}. In this framework, interactions are classified as either occurring within a single body (intrasolid) or between separate bodies (intersolid), allowing the interaction potential between two bodies to be expressed as a function of the gap vector. However, solving the boundary value problem for potential-based interactions at relevant time and length scales remains computationally intensive, largely due to the steep gradients of configuration-dependent interaction forces at small separations.

Molecular assemblies that resemble the shape of fibers and membranes are ubiquitous in nature. Biological examples of fiber-like macromolecules include proteins such as filamentous actin \cite{2012murrell} and collagen \cite{2021slepukhin}, nucleic acids like DNA and RNA \cite{2018franquelim}, as well as polysaccharides such as cellulose \cite{2018nishiyama} and fungal hyphae \cite{2017islama}. At comparable scales, biological membranes predominantly take the form of lipid mono- or bilayers that define cellular boundaries and internal organelles \cite{1995lipowsky}. Nature’s design has inspired the development of advanced materials and technologies using glass fibers \cite{2008alavinasab}, silicon nanotubes \cite{2012yoo}, carbon nanotubes \cite{2021čanadija}, mycelium \cite{2017islama}, and synthetic cell membranes \cite{2017lu}.

The interaction of biological membranes with fiber-like macromolecules plays a key role in processes such as adhesion of elastin fibers and smooth muscle cells \cite{1994perdomo}, interactions of membrane and actin cytoskeleton \cite{2018gov}, membrane shaping by BAR (Bin-Amphiphysin-Rvs) domain proteins \cite{2010saarikangas}, interaction of bacteria and virus with a cell \cite{2018charles-orszag}, etc. Protein-induced changes of the cell membrane are often approximated by prescribing a spontaneous curvature and assuming symmetry. However, the process is highly nonlinear and often leads to instabilities \cite{2015walani}. The state-of-the-art computational model for protein-induced budding of the cell membrane is presented in \cite{2017sauer} where asymmetric deformed shapes are obtained by prescribing the curvature. The authors in \cite{2018gov} emphasize that the deformability of a protein cannot be neglected and that more accurate nonlinear computational models are required. The interaction of nanorods and lipid membranes is modeled in \cite{2020maa} by the coarse-grained approach using the straight rod and flat shell finite elements. 
An alternative approach, based on the statistical thermodynamics of interactions of colloidal particles, is presented in \cite{2021tozzia}. Although the orientational order and the free energy of assemblies of fibers on membranes are successfully predicted, the approach is limited to 2D and axisymmetric analyses.

Fiber-membrane interaction also arises in the production of nanofiber membranes, which are widely used as biomimetic and mechanically stable surface coatings \cite{2016wang}. 
Nanofiber membranes are formed by an electrospinning process, which draws polymer solutions into ultrafine fibers and deposits them onto a collector. The adhesion between electrospun nanofibers and a collector is, in general, governed by the vdW interactions \cite{2020zhangb}. A stable integration of electrospun fibers onto the surfaces of elastomeric materials is considered in \cite{2020brunelli}. However, modeling fiber-substrate adhesion is often limited to ideal cases where the role of substrate geometry and mechanical properties is neglected. 
A numerical model that accounts for substrate roughness, patterning, curvature, and deformability using simple truss elements is presented in \cite{2018brely}. 
Furthermore, the adhesion of micromachined surfaces used in microelectro-mechanical systems (MEMS) at low roughness levels is primarily governed by vdW forces acting across extensive non-touching areas \cite{2005delrio}. The adhesion between carbon nanotubes (CNT) and elastic substrates with finite thickness is investigated in \cite{2018yuan}. CNTs are modeled as isotropic cylindrical shells, and substrates are treated as linear pure-bending plates. The obtained results are in partial agreement with those of molecular dynamics simulations. Finally, it should be noted that vdW interactions also play an important role in the electroless plating process that allows the combination of multiple heterogeneous structures and provides a many possibilities for 4D printed structures \cite{2025song}.

This brief literature overview suggests that vdW interactions between fibers and membranes are a critical topic across various fields of engineering and mechanics, which is yet to be fully understood. Aiming to advance the state-of-the-art in this area, the present study proposes an accurate and efficient computational model for vdW fiber-membrane interactions, based on the well-established finite element (FE) methodology. The literature on classical mechanical contact between beams and shells is scarce \cite{2010wriggers,2020gaynetob}, and the focus is mostly on beam-shell coupling \cite{2024sky}. Recently, the coarse-grained model has been incorporated into structural beam and shell theories \cite{2014sauerb, 2020grill, 2024mokhalingam}, offering a favorable compromise between accuracy and computational efficiency. A section-section approach, introduced in \cite{2020grill, 2023meier}, has been used to model short-range interactions between deformable planar beams \cite{2021grilla, 2024borkovićb, 2024borkovićg, 2025borković} and has been further extended to the long-range regime \cite{2026borkovića}. However, identifying an appropriate section-section law for short-range interactions between spatial beams remains a challenge. Alternative strategies, such as the section-beam model, have therefore been explored \cite{2023grill, 2024grill}. Motivated by these studies, we propose a computational formulation suited for short-range vdW fiber-membrane (or beam-shell) interactions. The approach is based on the analytical pre-integration of the interaction potential between the beam cross section and a surrogate plate representing the shell. For a decisive range of small separations, this approach promises good balance between accuracy and efficiency by: (i) capturing the main property of interaction force -- scaling w.r.t.~the gap, and (ii) reducing the 6D to 1D numerical integration. The resulting interaction potential is added to the total potential energy of a system consisting of a Bernoulli-Euler beam and a Kirchhoff-Love shell. The total potential is spatially discretized and solved using the isogeometric FE approach \cite{2005hughes}.

The modeling of intermolecular interactions between fibers and membranes is inherently complex and therefore simplified here by the following assumptions:
\begin{itemize}
	\item The total beam-shell interaction equals the pairwise summation (integration) of point-pair interactions.
	\item The point-point interaction potential is modeled as an inverse-power law of the point-pair distance.
	\item Only one beam and one shell are considered, and many-body effects are neglected.
	\item Any influences of a surrounding medium and retardation effects are neglected.
	\item There is no redistribution of particles or charges inside the bodies; that is, we are dealing with dielectric or nonconducting materials. 
	\item The density distributions of particles and physical constants over the interacting bodies are homogeneous at the initial configuration.
	\item Fibers have circular cross sections.
\end{itemize}

The main contributions of the present paper are: (i) development of a new approach for beam-shell vdW interactions using section-surrogate plate concept, (ii) derivation of the computational formulation that enables accurate and efficient modeling of highly nonlinear quasi-static interactions between deformable 3D fibers and membranes, and (iii) thorough numerical analysis of two fundamental and non-trivial examples concerning peeling of an adhering fiber from a membrane and bending of a membrane by an adhering fiber. To the best of our knowledge, no experimental or numerical examples suitable for direct comparison are available in the existing literature. These two examples are carefully selected to represent fundamental cases of fiber-membrane interaction and to serve as benchmarks for future studies. Furthermore, we provide all details on the variations and linearizations required for reproducibility, and we develop a novel disk-sphere vdW law that enables the computation of reference solutions.

The remaining paper is organized as follows: The problem of potential-based interactions between fibers and membranes is discussed in the next section. Computational models of beams and shells used in this research are scrutinized in Section \ref{sec:met}, and their potential energies are derived and varied. The beam-shell interaction potential based on the disk-infinite plate interaction is presented in Section \ref{secbsip}. Additionally, the error of the proposed approximation is thoroughly analyzed. The beam-shell potential is varied in Section \ref{sec:varsbs}, and the three formulations are suggested. The two proposed numerical experiments are given in Section \ref{secnum}, followed by conclusions in Section \ref{seccon}.

\section{Potential-based beam-shell interactions}


The concept of point-pair interaction potentials and their integration over the interacting beam and shell using the coarse-grained method are revisited in this section. The basic idea of the section-surrogate plate approach is introduced, and the general form of the equation of motion is presented.

\subsection{Coarse-graining of the beam-shell interaction potential}

Let us consider particles $i$ and $j$ that interact via a potential field modeled as an inverse power law w.r.t.~their distance $r_{ij}$. A point-point interaction potential of $m^\text{th}$ order, ${\Pi}^m_{\operatorname{P-P}}$, is the energy required to separate these particles from $r_{ij}$ to infinite separation, i.e.,
	\begin{equation}
		\label{eq: ip01}
		\begin{aligned}
			{\Pi}_{\operatorname{P-P}}^m (r_{ij}) = k_m \, r_{ij}^{-m},
		\end{aligned}
	\end{equation}
where $k_m$ is a physical constant. The interaction force is obtained as the gradient of this potential w.r.t.~the distance; it acts on both particles with the same intensity but in the opposite direction. By the pairwise summation concept \cite{2011israelachvili}, the volume interaction potential between a beam (B) and a shell (S) becomes
\begin{equation}
	\label{eq: ip01sum}
	\begin{aligned}
		\Pi_{\operatorname{{B-S}_{PW}}}^{m} &=  \sum_{i \in \mathrm{B}}^{} \sum_{j \in \mathrm{S}}^{} \Pi_{\operatorname{{P-P}}}^m (r_{ij}).
	\end{aligned}
\end{equation}
Let us approximate the pairwise summation as a volume integral over both bodies by using the coarse-graining procedure, introduced in \cite{2007sauer},
\begin{equation}
	\label{eq: ip01x}
	\begin{aligned}
		\Pi_{\operatorname{{B-S}_{PW}}}^{m} \approx \Pi_{\operatorname{B-S}}^{m} =  \int_{\ii{v}{}{\mathrm{B}}} \int_{\ii{v}{}{\mathrm{S}}} \beta_\mathrm{B} \, \beta_\mathrm{S} \, {\Pi}_{\operatorname{P-P}}^m \dd{v_\mathrm{S}} \dd{v_\mathrm{B}},
	\end{aligned}
\end{equation}
where $v_k$ $(k=\mathrm{B},\mathrm{S})$ are volumes, while $\beta_k$ are particle volume densities; both at the current configuration. Since we assume that the interacting property of an elementary volume is conserved during deformation, we have $\beta_{k \mathrm{0}} \dd{V}_k = \beta_{k} \dd{v_k}$.
This fact allows us to calculate an interaction potential at the current configuration by integrating over the reference volume $V_k$, using the reference particle densities $\beta_{k 0}$, i.e.,
\begin{equation}
	\label{eq: ip01xhh}
	\begin{aligned}
		\Pi_{\operatorname{B-S}}^{m} &= \int_{\ii{v}{}{\mathrm{B}}} \int_{\ii{v}{}{\mathrm{S}}} \beta_\mathrm{B} \,  \beta_\mathrm{S} \, k_m \, \norm{\tilde{\ve{x}}_\mathrm{B} - \tilde{\ve{x}}_\mathrm{S}}^{-m} \dd{v_\mathrm{S}} \dd{v_\mathrm{B}} \\
		&= \int_{V_\mathrm{B}} \int_{V_\mathrm{S}} \beta_\mathrm{B 0}  \, \beta_\mathrm{S 0} \, k_m \, \norm{\ve{x}_\mathrm{B} - \ve{x}_\mathrm{S}}^{-m} \dd{V_\mathrm{S}} \dd{V_\mathrm{B}},
	\end{aligned}
\end{equation}
where $\tilde{\ve{x}}_k$ and $\ve{x}_k$ are the current positions of elementary volumes w.r.t.~spatial and material coordinates, respectively. 
Calculating this integral for practical time and space resolutions is computationally expensive, which motivates us to employ a more efficient model.

For interactions between a beam and a shell, one approach would be to consider a pre-integration over a beam's cross-section and a shell's normal fiber. This would leave us with 3D numerical integration, which is still quite a demanding computational task. In the next subsection, we propose a decomposition of integral \eqref{eq: ip01xhh} into 5D analytical pre-integration and 1D numerical integration.

\subsection{Beam-shell interaction via disk-surrogate plate interaction}

\begin{figure}[h!]
	\centering
	\includegraphics[width=1\textwidth]{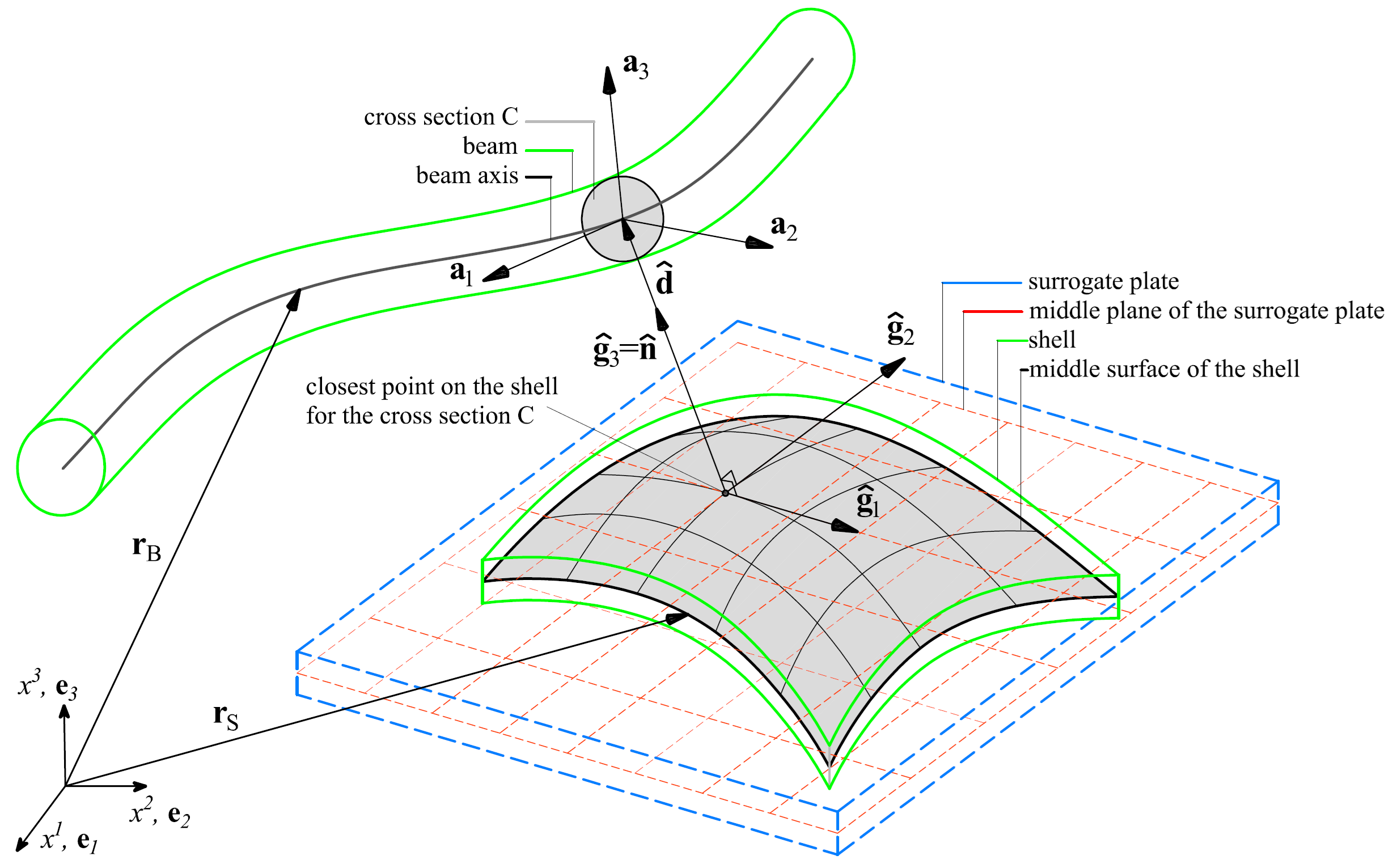} 
	\caption{Beam-shell interaction. The interaction between a beam's disk-shaped cross section C and a shell is approximated as an interaction between disk C and an infinite plate, tangential to the shell at the closest point. The hat accent designates a quantity evaluated at the closest point.}
	\label{fig:FigBS}
\end{figure}

Following the section-surrogate cylinder idea from \cite{2023grill}, we pursue an efficient computational approach for beam-shell interaction. Let us consider a beam and a shell, as illustrated in Fig.~\ref{fig:FigBS}.
By integrating the potential over the area of the beam's cross section, $A_\mathrm{B}$, the interaction reduces to a 4D integral: one along the length of the beam axis, $L_\mathrm{B}$, and three over the shell volume. Recall that the interaction potential is the inverse-power law w.r.t.~gap \eqref{eq: ip01}. If the gap between the beam and the shell is small, the contribution to the interaction potential comes only from very small region. This observation enables the approximation of the disk-shell interaction by a disk-surrogate plate interaction, where the surrogate plate is an infinite-area plate tangent to the shell at the point closest to the disk and has the same thickness as the shell. The main advantage of this approach is that we can analytically integrate over the volume of the surrogate plate, $V_\mathrm{SP}$, i.e.,
%
\begin{equation}
	\label{eq: ip01xhah}
	\begin{aligned}
		\Pi_{\operatorname{B-S}}^{m} &= \int_{L_\mathrm{B}} \int_{A_\mathrm{B}} \int_{V_\mathrm{S}}  \beta_\mathrm{B0}  \, \beta_\mathrm{S0} \, k_m \, \norm{\ve{x}_\mathrm{B} - \ve{x}_\mathrm{S}}^{-m} \dd{V_\mathrm{S}} \dd{A_\mathrm{B}}\dd{s_\mathrm{B}} \approx \int_{L_\mathrm{B}} \beta_\mathrm{B0}  \, \beta_\mathrm{S0} \, k_m  \Pi_{\operatorname{D-SP}}^m\dd{s_\mathrm{B}}, \\
		\Pi_{\operatorname{D-SP}}^m&:=\int_{A_\mathrm{B}} \int_{V_\mathrm{SP}}  \, \norm{\ve{x}_\mathrm{B} - \ve{x}_\mathrm{S}}^{-m} \dd{V_\mathrm{SP}} \dd{A_\mathrm{B}},
	\end{aligned}
\end{equation}
where $\Pi_{\operatorname{D-SP}}^m$ is the disk-surrogate plate interaction law. In this way, the beam-shell interaction potential is approximated with a 1D integral over the beam length. We will refer to this approach as the surrogate beam-shell model (SBS). An exact analytical integration of the interaction potential over a disk and an infinite plate has been obtained in \cite{2025borković} for an arbitrary exponent $m$.

Let us emphasize that the closest point-pairs dominate the interaction when bodies are in close proximity, and the exponent of the potential is $m>3$. This short-range effect leads to steep gradients of the interacting force, reflecting a complex equilibrium state between repulsive and attractive forces. Addressing the high sensitivity of interaction forces at the interface is crucial for accurate simulations but also computationally demanding. On the other hand, description of long- and moderate-range effects is less involved since the distributions of interaction forces are smoothed. Although these effects can be important for some problems \cite{2026borkovića}, they are not considered in this paper.

\subsection{Lennard-Jones potential and equation of motion}

To model adhesion due to vdW forces realistically, repulsive effects must also be included. The repulsion develops between bodies in close vicinity due to overlapping electron clouds. This repulsion effect is observed as \emph{contact} from a macroscopic point of view. vdW and steric interactions exist for practically all bodies, making them one of the most common forces in nature. By modeling the repulsive steric potential with an inverse-power law with $m=12$, and adding it to the vdW potential, we obtain the well-known Lennard-Jones (LJ) potential between two particles with distance $r$
	\begin{equation}
		\begin{aligned}
			\Pi_{\operatorname{P-P}}^{\operatorname{LJ}} (r)=4\epsilon \left[\left(\frac{\sigma}{r}\right)^{12}-\left(\frac{\sigma}{r}\right)^6\right] = k_6 \, r^{-6} + k_{12} \, r^{-12},
		\end{aligned}
		\label{eqLJex}
	\end{equation}
where $\sigma$ is the distance at which the potential is zero, while $\epsilon$ is the minimum value of the potential. 
In our numerical analysis, we exclusively consider LJ beam-shell interactions that stem from the integration of the point-point LJ potential \eqref{eqLJex}.

	
The computational modeling of potential-based interactions between bodies requires solving an appropriate (initial) boundary value problem. The strong form of this problem consists of the equation of motion (balance of linear momentum), boundary and initial conditions, and a constitutive relation. To derive the weak form, we define the total potential energy of a system involving elastic bodies, $\Pi_{\mathrm{tot}}$. If we neglect inertial effects and consider interaction between a beam and a shell, this energy consists of the strain energy, $\Pi_{\mathrm{str}}=\Pi_{\mathrm{B,str}}+\Pi_{\mathrm{S,str}}$, the work of external forces, $\Pi_{\mathrm{ext}}=\Pi_{\mathrm{B,ext}}+\Pi_{\mathrm{S,ext}}$, and the interaction potential, $\Pi_{\operatorname{B-S}}^m$. By the principle of stationary total potential energy, the weak form of the boundary value problem follows
\begin{equation}
	\label{eq: poten}
	\begin{aligned}
		\delta \Pi_{\mathrm{tot}} (\ve{x}_\mathrm{B},\ve{x}_\mathrm{S})= \delta \Pi_{\mathrm{str}} -\delta \Pi_{\mathrm{ext}}+\delta \Pi_{\operatorname{B-S}}^m=0,
	\end{aligned}
\end{equation}
and it must hold for all admissible variations of the configuration of beam and shell. Potential energies of beams and a shells are derived and varied in the next section, while the interaction potential energy is derived in Section \ref{secbsip} and varied in Section \ref{sec:varsbs}. 

The obtained equilibrium equation \eqref{eq: poten} is spatially discretized using the isogeometric approach \cite{2018borkovic}, details of which are omitted for brevity. The resulting system of algebraic equations is highly nonlinear, and a continuation approach is employed in which the forcing term is incremented, allowing the equilibrium configurations to be determined at a finite number of quasi-time steps. The Newton-Raphson method is employed to solve the system of equations within each increment, with details of the required linearizations provided in Appendices \ref{appendixa}, \ref{appendixb}, \ref{appendixd}, and \ref{appendixe}.

	\section{Potential energy of beams and shells}
	\label{sec:met}
For the sake of completeness, the beam and shell models that are employed in this work are briefly discussed. Since our goal is the efficient modeling of interactions between thin beams and shells, the Bernoulli-Euler (BE) beam and Kirchhoff-Love (KL) shell models are employed. Both of these theories reduce general deformable 3D continua via plausible assumptions. By assuming that one dimension of the 3D body is much smaller than the other two, we obtain the mechanical model of a shell. If we further assume that the fibers perpendicular to the midsurface of the shell are rigid and perpendicular to the deformed midsurface, we obtain the KL shell theory. This mechanical model allows the 3D body to be described solely by the position of the shell’s midsurface. 

On the other hand, by assuming that one dimension of the 3D body is much smaller than the other two, we obtain a mechanical model of a beam. By assuming that the cross sections are rigid and perpendicular to the deformed beam axis, the BE beam model is found. Due to additional dimensional reduction, in comparison with the KL shell, the BE beam model requires one rotational DOF (twist angle) besides the position of the axis. 

To describe the deformation of beams and shells, we introduce a reference configuration that coincides with the initial, undeformed state. The kinematic quantities for both the reference and current configurations can be obtained in an analogous manner. We will distinguish the former by using capital letters. 
In the following, Greek indices take the values $\alpha=1,2$ for the shell, and $\alpha=2,3$ for the beam model, while Latin indices take the values $i=1,2,3$.


\subsection{Kirchhoff-Love shell}

Metrics of the shell's midsurface, equidistant surface, and the variation of the KL shell's strain energy are given in this subsection, c.f.~\cite{2021radenković,2017sauerf,2023gfrerer}.

\subsubsection{Metric of a midsurface}

\begin{figure}[h]
	\centering
	\includegraphics[width=1\textwidth]{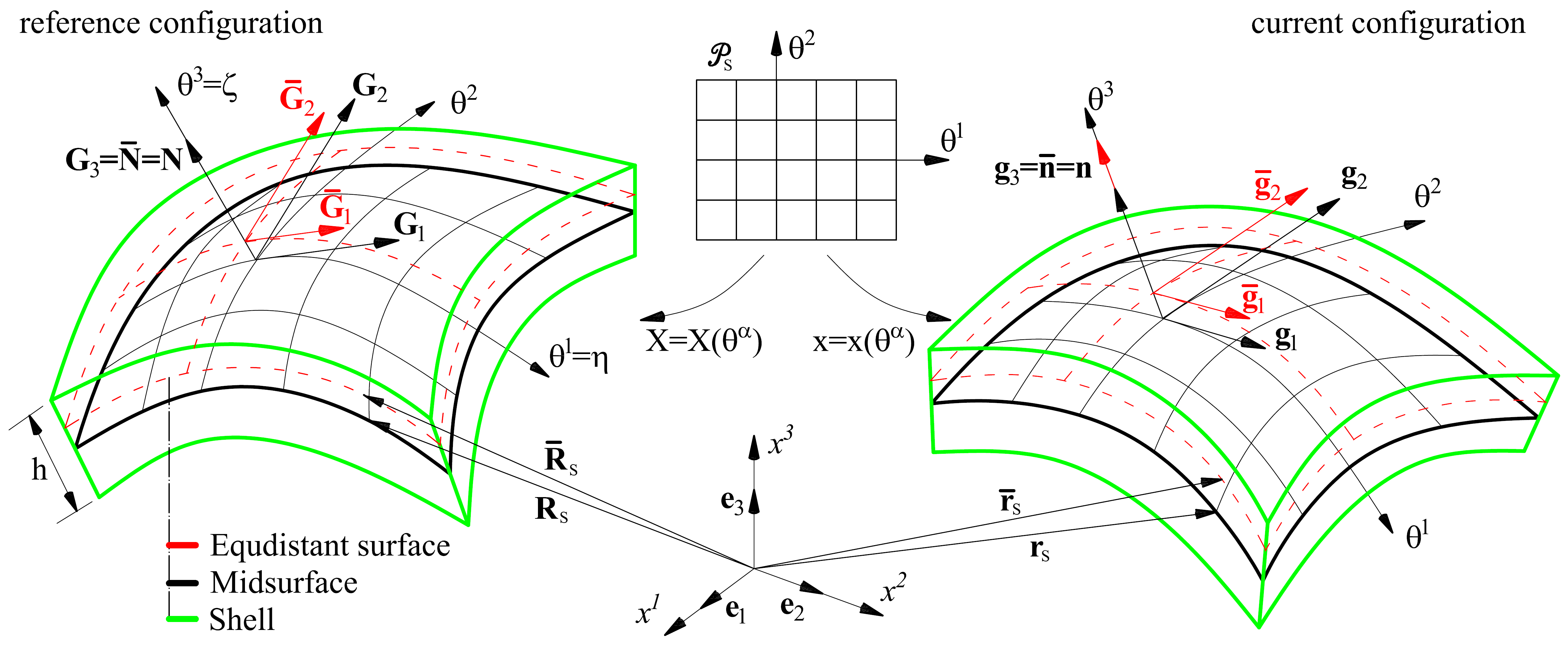} 
	\caption{Schematic representation of the shell model in reference and current configurations.}
	\label{fig:shell}
\end{figure}

The position of the shell's midsurface at the current configuration is defined by the mapping
\begin{equation}
	\begin{aligned}
		\ivn{r}{}{S} &= \ivn{r}{}{S} (\theta^\alpha) , \quad \alpha=1,2,
	\end{aligned}
\end{equation}
where $\ve{r}_\text{S} (\theta^\alpha) \in \mathbb{R}^3$ is the position vector while $\theta^\alpha$ are curvilinear coordinates
, Fig.~\ref{fig:shell}.
The current position of the midsurface is related to the reference configuration via displacement vector $\ivn{u}{}{S}$,
\begin{equation}
	\ivn{r}{}{S} (\theta^\alpha) = \ivn{R}{}{S} (\theta^\alpha) + \ivn{u}{}{S} (\theta^\alpha).
\end{equation}
Then we define covariant base vectors and covariant components of the surface metric
\begin{equation}
	\begin{aligned}
		\iv{g}{}{\alpha} &= \frac{\partial \ivn{r}{}{S}}{\partial \theta^\alpha}, \quad g_{\alpha \beta} = \iv{g}{}{\alpha} \cdot \iv{g}{}{\beta},
	\end{aligned}
\end{equation}
contravariant basis vectors and contravariant components of the metric
\begin{equation}
	\begin{aligned}
		\iv{g}{}{\alpha} \cdot \iv{g}{\beta}{} &= \delta_\alpha^\beta, \quad g^{\alpha \beta} = \iv{g}{\alpha}{} \cdot \iv{g}{\beta}{}, \quad [g^{\alpha \beta}] = [g_{\alpha \beta}]^{-1}	,
	\end{aligned}
\end{equation}
normal vector
\begin{equation}
	\begin{aligned}
		\iv{g}{}{3} & = \iv{n}{}{} = \frac{\iv{g}{}{1} \times \iv{g}{}{2}}{\norm{\iv{g}{}{1} \times \iv{g}{}{2}}},
	\end{aligned}
\end{equation}
and the Jacobian of the coordinate transformation between the area element in the parametric and physical domains
\begin{equation}
	\begin{aligned}
		\dd{a}_\mathrm{S} &= j_\text{S} \dd{\theta^1} \dd{\theta^2}, \quad j_\text{S} = \sqrt{\det[g_{\alpha\beta}]}	.
	\end{aligned}
\end{equation}
It is convenient to decompose the full identity tensor in $\mathbb{R}^3$ as
\begin{equation}
\mathbf{I} = \iv{i}{}{} + \iv{n}{}{} \otimes \iv{n}{}{}, \quad \iv{i}{}{} = \iv{g}{}{\alpha} \otimes \iv{g}{\alpha}{} = \iv{g}{\alpha}{} \otimes \iv{g}{}{\alpha},
\end{equation}
where $\iv{i}{}{}$ is the surface identity tensor. Furthermore, components of the curvature tensor follow from the Gauss-Weingarten equation
\begin{equation}
	\label{eq:shellcurcurr}
	\begin{aligned}
		\iv{b}{}{} & = b_{\alpha \beta} \,\iv{g}{\alpha}{} \otimes \iv{g}{\beta}{}, \quad b_{\alpha \beta} &:= \iv{n}{}{} \cdot \iv{g}{}{\alpha , \beta} = - \iv{n}{}{,\beta} \cdot \iv{g}{}{\alpha}.
	\end{aligned}
\end{equation}
As aforementioned, we can define analogous quantities in the reference configuration: $\ivn{R}{}{S}$, $\iv{G}{}{\alpha}$, $\iv{G}{\alpha}{}$, $J_\mathrm{S}$, $G_{\alpha \beta}$, $G^{\alpha\beta}$, $\iv{N}{}{}$, $\iv{B}{}{}$.

By having both fundamental metric forms defined in the reference, $G_{\alpha\beta}$ and $B_{\alpha\beta}$, and the current, $g_{\alpha\beta}$ and $b_{\alpha\beta}$, configurations, the deformation can be described by two quantities: the surface Green-Lagrange tensor $\iv{E}{}{}$ and the change of curvature tensor $\iv{K}{}{}$. By defining the surface deformation gradient as
\begin{equation}
	\iv{F}{}{} := \iv{g}{}{\alpha} \otimes \iv{G}{\alpha}{},
\end{equation}
the surface Green-Lagrange strain tensor follows
\begin{equation}
	\ivn{E}{}{} = E_{\alpha\beta} \,\iv{G}{\alpha}{} \otimes \iv{G}{\beta}{} := \frac{1}{2} \big(\ivn{F}{T}{}\iv{F}{}{} - \iv{i}{}{}\big) = \frac{1}{2} \big(g_{\alpha\beta} - G_{\alpha\beta}\big) \iv{G}{\alpha}{} \otimes \iv{G}{\beta}{}.
\end{equation}
The change of curvature tensor is defined as
\begin{equation}
	\iv{K}{}{} = K_{\alpha\beta} \,\iv{G}{\alpha}{} \otimes \iv{G}{\beta}{} := \pmb{\kappa} - \ivn{B}{}{} = \big(b_{\alpha\beta} - B_{\alpha\beta}\big) \iv{G}{\alpha}{} \otimes \iv{G}{\beta}{},
\end{equation}
where $\pmb{\kappa}$ is the pull-back of the curvature tensor of the current configuration to the reference basis, i.e.,
\begin{equation}
	\label{eq:shellcurref}
	\pmb{\kappa} := \ivn{F}{T}{} \iv{b}{}{} \iv{F}{}{} = b_{\alpha\beta} \, \iv{G}{\alpha}{} \otimes \iv{G}{\beta}{}.
\end{equation}

\subsubsection{Metric of an equidistant surface}

Since normal fibers remain normal and rigid during the motion, the geometry of an equidistant surface can be represented via the geometry of the midsurface, i.e.,
\begin{equation}
	\begin{aligned}
		\label{eqs11}
		\iveqn{r}{}{S} (\theta^\alpha) &= \ivn{r}{}{S} (\theta^\alpha) + \zeta\ve{n} (\theta^\alpha),
	\end{aligned}
\end{equation}
where $\zeta=\theta^3 \in [-h/2,h/2]$ is the shell's thickness coordinate and the overbar denotes a quantity at an equidistant surface. The normal vector of an equidistant surface is the same as at the midsurface, $\iveq{n}{}{} (\theta^\alpha)= \iv{n}{}{} (\theta^\alpha)$, while the tangent vectors and the metric tensor are
\begin{equation}
		\label{eqs2}
	\begin{aligned}
		\iveq{g}{}{\alpha} &:= \iveq{r}{}{\mathrm{S},\alpha} = \iv{g}{}{\alpha} - \zeta \,b_\alpha^\beta\, \iv{g}{}{\beta}, \quad 	\ieq{g}{}{\alpha\beta} := \iveq{g}{}{\alpha} \cdot \iveq{g}{}{\beta} \approx g_{\alpha\beta} -2\zeta \,b_{\alpha\beta}, 
	\end{aligned}
\end{equation}
where we approximate by neglecting the second-order term w.r.t.~the $\zeta$-coordinate.
Furthermore, let us define contravariant basis vectors
\begin{equation}
		\label{eqs3}
	\begin{aligned}
		\iveq{g}{}{\alpha} \cdot \iveq{g}{\beta}{} &= \delta_\alpha^\beta, \quad \ieq{g}{\alpha \beta}{} = \iveq{g}{\alpha}{} \cdot \iveq{g}{\beta}{}, \quad [\bar g^{\alpha \beta}] = [\bar g_{\alpha \beta}]^{-1}	,
	\end{aligned}
\end{equation}
and the Jacobian of the coordinate transformation between the area element in the parametric and physical domains
\begin{equation}
		\label{eqs4}
	\begin{aligned}
		\dd{\bar a}_\mathrm{S} &= \bar j_\text{S} \dd{\theta^1} \dd{\theta^2}, \quad \bar j_\text{S} = \sqrt{\det[\bar g_{\alpha\beta}]}	.
	\end{aligned}
\end{equation}
In the reference configuration, quantities $\iveqn{R}{}{S}$, $\iveq{G}{}{\alpha}$, $\iveq{G}{\alpha}{}$, $\ieq{G}{}{\alpha\beta}$, $\ieq{G}{\alpha\beta}{}$, and $\bar J_\mathrm{S}$ can be defined analogously to Eqs.~\ref{eqs11}, \ref{eqs2}, \ref{eqs3}, and \ref{eqs4}. The Green-Lagrange strain at an equidistant surface is now
\begin{equation}
	\begin{aligned}
		\iveq{E}{}{} := \ieq{E}{}{\alpha\beta}\big( \iveq{G}{\alpha}{}\otimes \iveq{G}{\beta}{}\big) = \frac{1}{2} \big(\ieq{g}{}{\alpha\beta} - \ieq{G}{}{\alpha\beta}\big) \big( \iveq{G}{\alpha}{}\otimes \iveq{G}{\beta}{}\big)= \big(\frac{1}{2} E_{\alpha\beta} - \zeta K_{\alpha\beta}\big) \big( \iveq{G}{\alpha}{}\otimes \iveq{G}{\beta}{}\big).
	\end{aligned}
\end{equation}

\subsubsection{Variation of the strain energy}

We assume the existence of the strain energy density for the Saint Venant--Kirchhoff model
\begin{equation}
	\tilde \Pi_\mathrm{S,str} := \frac{1}{2} \iveq{E}{}{}  \colon \mathbb{\bar{C}} \colon \iveq{E}{}{} = \frac{1}{2} \iveq{E}{}{}  \colon \iveq{S}{}{}, \quad \iveq{S}{}{} = \lambda \tr (\iveq{E}{}{}) \iv{I}{}{} + 2 \mu \iveq{E}{}{},
\end{equation}
where $\mathbb{\bar{C}}$ is the constitutive tensor, $\iveq{S}{}{}$ is the second Piola-Kirchhoff stress at an equidistant surface, while $\lambda$ and $\mu$ are Lam\'e constants. Furthermore, we assume that both the plane stress and the plane strain conditions are satisfied, which allows us to express components of the second Piola-Kirchhoff stress tensor as
\begin{equation}
	\ieq{S}{\alpha\beta}{} = 2\mu \big(\bar G^{\alpha\nu} \bar G^{\beta\gamma} + \frac{\nu}{1-\nu} \bar G^{\alpha\beta} \bar G^{\nu\gamma}\big) \ieq{E}{}{\nu\gamma} = \bar C^{\alpha\beta\nu\gamma} \ieq{E}{}{\nu\gamma}.
\end{equation}
The variation of the strain energy
\begin{equation}
	\delta 	\ii{\Pi}{}{\mathrm{S,str}} := \int_{V_\mathrm{S}}^{} \delta 	\tilde \Pi_\mathrm{S,str} \dd{\bar V}_\mathrm{S} = \int_{A_\mathrm{S}}^{} \int_{-h/2}^{h/2} \ieq{S}{\alpha\beta}{} \, \delta \bar E_{\alpha\beta} \dd{\zeta}  \bar J_\text{S} \dd{\theta^1} \dd{\theta^2}
\end{equation}
can be analytically integrated along the shell's thickness by neglecting higher order terms w.r.t.~$\zeta$, i.e.,
\begin{equation}
	\label{shellvp}
	\delta 	\Pi_\mathrm{S,str} =  \int_{A_\mathrm{S}}^{} \big(\ii{N}{\alpha\beta}{} \delta E_{\alpha\beta} + \ii{M}{\alpha\beta}{} \delta K_{\alpha\beta}\big) J_\text{S} \dd{\theta^1} \dd{\theta^2}, 
\end{equation}
where
\begin{equation}
	\label{shellsc}
	\begin{aligned}
		\ii{N}{\alpha\beta}{} &= h C^{\alpha\beta\nu\gamma} \ii{E}{}{\nu\gamma} , \quad 
		\ii{M}{\alpha\beta}{} = \frac{h^3}{12} C^{\alpha\beta\nu\gamma} \ii{K}{}{\nu\gamma},
	\end{aligned}
\end{equation}
are section forces and section moments. 

For the variation of the shell's strain energy in \eqref{shellvp}, we require the variation of the membrane strains 
%
\begin{equation}
	\delta E_{\alpha\beta} = \frac{1}{2} (\iv{g}{}{\alpha} \cdot \delta \iv{u}{}{\text{S},\beta} + \iv{g}{}{\beta}\cdot \delta \iv{u}{}{\text{S},\alpha}),
\end{equation}
and the variation of the change of curvatures
\begin{equation}
	\delta K_{\alpha\beta} = \delta b_{\alpha\beta} = \delta (\iv{g}{}{\alpha,\beta}\cdot\ve{n}) = \iv{n}{}{} \cdot \left(\delta \iv{u}{}{\text{S},\alpha\beta} - \ii{\Gamma}{\gamma}{\alpha\beta} \cdot \delta \iv{u}{}{\text{S},\gamma}\right),
\end{equation}
where we have used $\delta \ve{n}=-(\ve{n}\cdot\iv{g}{}{\alpha}) \iv{g}{\alpha}{}$, while $\ii{\Gamma}{\gamma}{\alpha\beta} = \iv{g}{}{\alpha,\beta} \cdot \iv{g}{\gamma}{}$ are Christoffel symbols of the second kind. Linearization of the variation of the shell's strain energy is given in Appendix \ref{appendixa}.

\subsection{Bernoulli-Euler beam}
	
\label{besec}

The metric of the beam axis and an equidistant line, parametrization of rotation, variation of the beam's strain energy, and variation of the work of the point moment are given in this subsection, c.f.~\cite{2023borković,2022borkovićc,2015meier,2020vo}.

\subsubsection{Metric of the beam axis}

\begin{figure}[h]
	\centering
	\includegraphics[width=1\textwidth]{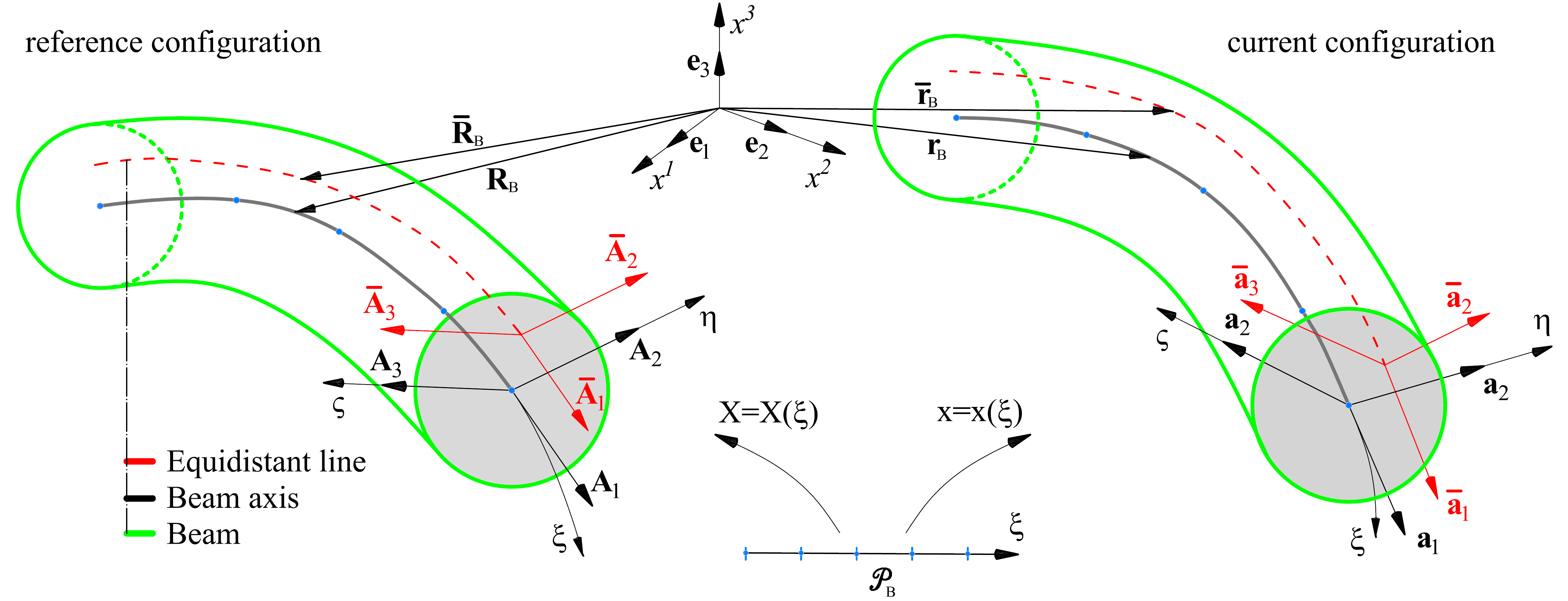} 
	\caption{Schematic representation of the beam model in reference and current configurations.}
	\label{fig:FigB}
\end{figure}

In the reference configuration, the beam axis is a continuous space curve $\ve{R}_\text{B} (\xi) \in \mathbb{R}^3$ with the tangential basis vector $\ve{A}_1 (\xi)=\partial \ve{R}_\mathrm{B} / \partial \xi$, while $\xi$ is a parametric convective coordinate, Fig.~\ref{fig:FigB}.
\noindent At each point along the beam axis, we define an orthonormal basis $\{\ve{T},\ve{A}_2,\ve{A}_3\}$ where $\ve{T}=\ve{A}_1 / \norm{\ve{A}_1}$ while $\ve{A}_\alpha$ are aligned with the principal axes of the beam cross sections. Vectors $\ve{A}_\alpha$, often referred to as \emph{directors}, can be defined in several ways, either by the Frenet-Serret frame, Bishop frame, or some other algorithm.

The position vector of the current configuration is $\ve{r}_\mathrm{B} = \ve{R}_\mathrm{B}+\ve{u}_\mathrm{B}$, where $\ve{u}_\mathrm{B}$ is the displacement vector. The tangent in the current configuration is
\begin{equation}
	\label{eq:tangb}
	\begin{aligned}
		\ve{t} =\frac{\partial \ivn{r}{}{B}}{\partial s} =\frac{\partial \ivn{r}{}{B}}{\partial \xi} \frac{\partial \ii{\xi}{}{}}{\partial s} = \iv{r}{}{\mathrm{B},1} \xi_{,s}=\iv{a}{}{1}/j_\mathrm{B}, \quad j_\mathrm{B}= \sqrt{a}=\norm{\iv{a}{}{1}} = \sqrt{\iv{a}{}{1}\cdot\iv{a}{}{1}},
	\end{aligned}	
\end{equation}
where $j_\mathrm{B}$ is the Jacobian of the coordinate transformation between physical (arc length) and parametric convective coordinate, $\dd s=j_\mathrm{B} \dd \xi$.

Directors in the current configuration are found by rotating the directors from the reference configuration as $\ve{a}_\alpha = \pmb{\Lambda} \ve{A}_\alpha$, where $\pmb{\Lambda}$ is a rotation tensor, i.e.,
\begin{equation}
	\pmb{\Lambda}:\mathbb{R}^3\rightarrow \mathbb{R}^3, \; \pmb{\Lambda}\pmb{\Lambda}^\mathrm{T}=\ve{I}, \;\det \pmb{\Lambda} =1.
\end{equation} 
There are many ways to parametrize the rotation tensor, one of them being via the rotation vector $\pmb{\theta}$, $\pmb{\Lambda} = \mathrm{exp} (\pmb{\theta}_\times)$,
%
where $(\cdot)_\times$ represents the skew-symmetric tensor of the vector $(\cdot)$.
Since we assume rigid cross sections that remain orthogonal to the beam axis, the current triad $\left\{\ve{t},\ve{a}_2,\ve{a}_3\right\}$ is orthonormal and uniquely defined by the position of the beam axis and a single scalar rotation parameter.

Let us consider the so-called smallest rotation (SR) algorithm,  which consists of two steps. First, we rotate the $\ve{A}_\alpha$ vectors into the plane with normal $\ve{t}$ by the smallest angle, and obtain vectors $\ve{a}_\alpha^\mathrm{SR}$, i.e.,
\begin{equation}
	\label{eq:sr}
	\ve{a}_\alpha^\mathrm{SR} = \ve{A}_\alpha - \frac{\ve{A}_\alpha \cdot \ve{t}}{t_T} \left(\ve{t} + \ve{T}\right), \quad t_T:= 1+\ve{T}\cdot\ve{t}, \quad \alpha = 2,3.
\end{equation}
Second, vectors $\ve{a}_\alpha^\mathrm{SR}$ are rotated by an angle $\varphi$ in the current cross-sectional plane
\begin{equation}
	\label{eq:sr2}
	\begin{bmatrix}
		\ve{a}_2 \\
		\ve{a}_3
	\end{bmatrix}  =
		\begin{bmatrix}
			\cos \varphi & \sin \varphi\\
			-\sin \varphi & \cos\varphi 
		\end{bmatrix}
	\begin{bmatrix}
			\ve{a}_2^\mathrm{SR} \\
			\ve{a}_3^\mathrm{SR}
	\end{bmatrix}.
\end{equation}
The SR algorithm is accurate and straightforward but has two limitations: (i) singularity for $\ve{T} \cdot \ve{t}$, and (ii) loss of objectivity for sparse FE meshes. The singularity issue is easily addressed with an incremental update of rotations, implemented here. Moreover, resolving the interaction force field in envisioned applications requires dense FE meshes, which alleviates the objectivity issue.

Covariant and contravariant components of the metric of the beam axis are 
\begin{equation}
	\label{eq:metriccur}
		\iv{a}{}{i} \cdot \iv{a}{j}{} = \delta_i^j, \quad a_{ij} = \iv{a}{}{i} \cdot \iv{a}{}{j}, \quad a^{ij} = \iv{a}{i}{} \cdot \iv{a}{j}{}, \quad [a^{ij}] = [a_{ij}]^{-1}, \quad \det[a_{ij}]=a=\ve{a}_1
\cdot \ve{a}_1,
\end{equation}
while the derivatives of the basis vectors are
\begin{equation}
	\label{eq: def: derivatives of base vectors}
	\ve{a}_{i,1} =   \Gamma^{k}_{i1} \ve{a}_k, \quad
	\begin{bmatrix}
		\iv{a}{}{1,1}\\
		\iv{a}{}{2,1}\\
		\iv{a}{}{3,1}
	\end{bmatrix}
	=
	\begin{bmatrix}
		\Gamma & a\ii{k}{}{3} & -a\ii{k}{}{2}\\
		-\ii{k}{}{3} & 0 & \ii{k}{}{1} \\
		\ii{k}{}{2} & -\ii{k}{}{1} & 0
	\end{bmatrix}
	\begin{bmatrix}
		\iv{a}{}{1}\\
		\iv{a}{}{2}\\
		\iv{a}{}{3}
	\end{bmatrix},
\end{equation}
where 
$\ii{k}{}{i}$ are the components of the curvature vector with respect to the reciprocal basis $\iv{a}{i}{}$. $\ii{k}{}{1}$ represents the torsion of the cross section, while $\ii{k}{}{\alpha}$ are components of the curvature of the beam axis,
\begin{equation}
	\label{eq: 25}
	\ii{k}{}{1} = \iv{g}{}{2,1}\cdot \iv{g}{}{3} = \iv{g}{\text{SR}}{2,1}\cdot \iv{g}{\text{SR}}{3} + \varphi_{,1}, \quad \ii{k}{}{2} = -\iv{a}{}{1,1} \cdot \iv{a}{}{3}/a, \quad \ii{k}{}{3} = \iv{g}{}{1,1} \cdot \iv{g}{}{2}/a.
\end{equation}
The components of the curvature vector are identical w.r.t.~both the reference and current bases, analogous to the shell model, see Eqs.~\eqref{eq:shellcurcurr} and \eqref{eq:shellcurref}.

\subsubsection{Metric of an equidistant line}

The position vector and basis vectors at an equidistant line are
\begin{equation}
	\label{eq:equi}
	\begin{aligned}
		\veq{r}_\mathrm{B} &=\ve{r} + \eta \ve{a}_2 + \varsigma \ve{a}_3,\\
		\veq{a}_1 &=\veq{r}_{\mathrm{B},1} = \ve{a}_1 + \eta \ve{a}_{2,1} + \varsigma \ve{a}_{3,1}= a_0 \ve{a}_1 - \varsigma k_1 \ve{a}_2 + \eta k_1 \ve{a}_3, \quad \veq{a}_\alpha = \ve{a}_\alpha,
	\end{aligned}	
\end{equation}
where $a_0=1-\eta k_3 + \varsigma k_2$ is the \emph{shifter}. The metric of an equidistant line is
\begin{equation}
	\begin{aligned}
	\label{eq:mettens}
\ieq{a}{}{ij}
	&:= \iveq{a}{}{i} \cdot \iveq{a}{}{j},
	\quad \det [\ieq{a}{}{ij}] = \ieq{a}{}{} = a_0^2 \; a , \quad \ieq{a}{}{11} = \iveq{a}{}{1} \cdot \iveq{a}{}{1} \approx (1-2\eta k_3 + 2 \varsigma k_2) a,\\
	\end{aligned}
\end{equation}
where we approximate by neglecting higher order terms w.r.t.~cross-sectional coordinates. The contravariant basis vectors are
\begin{equation}
	\label{eqsbrec}
	\begin{aligned}
		\iveq{a}{}{i} \cdot \iveq{a}{j}{} &= \delta_i^j, \quad \ieq{a}{ij}{} = \iveq{a}{i}{} \cdot \iveq{a}{j}{}, \quad [\bar a^{ij}] = [\bar a_{ij}]^{-1}	,
	\end{aligned}
\end{equation}
and the Jacobian of the coordinate transformation is
\begin{equation}
	\label{eqs4}
	\begin{aligned}
		\dd{\bar s} &= \bar j_\text{B} \dd{\xi} = \sqrt{\bar{a}} \dd{\xi}.
	\end{aligned}
\end{equation}
Analogous to Eqs.~\eqref{eq:metriccur}, \eqref{eq: 25}, and \eqref{eq:mettens}, we can calculate all required quantities in the reference configuration: $K_1$, $K_\alpha$, $A_{ij}$, $\ieq{A}{}{ij}$, $J_\mathrm{B}$, and $A$.

Due to BE assumptions, non-zero components of the Green-Lagrange strain tensor are $\epsilon_{11}$ and $\epsilon_{1\alpha}$, i.e.,
\begin{equation}
	\label{eq:beamstrain}
	\begin{aligned}
		\bar{\epsilon}_{11} &= \frac{1}{2}(\ieq{a}{}{11} - \ieq{A}{}{11}) = \epsilon_{11}  -\eta \kappa_3 + \varsigma \kappa_2,\\
		\bar{\epsilon}_{12} &= \frac{1}{2}(\ieq{a}{}{12} - \ieq{A}{}{12}) = - \varsigma \kappa_1, \quad \bar{\epsilon}_{13} = \frac{1}{2}(\ieq{a}{}{13} - \ieq{A}{}{13}) = - \eta \kappa_1, \\
		\epsilon_{11}&:= a_{11} - A_{11}, \quad \kappa_1 := k_1-K_1, \quad \kappa_\alpha := a (k_\alpha-K_\alpha),
	\end{aligned}
\end{equation}
where $\epsilon_{11}$ is the axial strain of the beam axis, $\kappa_1$ is the torsional curvature of the cross section, and $\kappa_\alpha$ are the bending curvatures of the beam axis. These quantities represent \emph{reference strains} of the BE beam.

%

\subsubsection{Variation of the beam potential energy}

To define components of the second Piola-Kirchhoff stress $\bar{S}^{ij}$, we employ the Saint Venant-Kirchhoff material model, impose conditions $\bar{S}^{22}=\bar{S}^{33}=0$, neglect higher-order terms w.r.t.~cross-sectional coordinates, and set off-diagonal terms in the reciprocal metric to $0$, i.e.,
\begin{equation}
	\bar{S}^{11} = E \bar{a}^{11}\bar{a}^{11} \bar{\epsilon}_{11}\approx \frac{E}{a^2} \bar{\epsilon}_{11}, \; \bar{S}^{12} = 2 \mu \bar{a}^{11}\bar{a}^{22} \bar{\epsilon}_{12} \approx 2 \frac{\mu}{a} \bar{\epsilon}_{12} , \; \bar{S}^{13} = 2 \mu \bar{a}^{11}\bar{a}^{33} \bar{\epsilon}_{13} \approx 2\frac{\mu}{a}\bar{\epsilon}_{13},
\end{equation}
where $E$ is Young's modulus. The variation of the strain energy 
\begin{equation}
	\delta 	\ii{\Pi}{}{\mathrm{B,str}} = \int_{V_\mathrm{B}}^{} \bar{S}^{1j} \delta \bar{\epsilon}_{1j}\dd{\bar V}_\mathrm{B} = \int_{L_\mathrm{B}}^{} \int_{A_\mathrm{B}}^{}  \bar{S}^{1j} \delta \bar{\epsilon}_{1j} \dd{\eta} \dd{\varsigma} \bar{J}_\mathrm{B} \dd{\xi}
\end{equation}
can now be analytically integrated by assuming $\bar{J}_\mathrm{B} = J_\mathrm{B}$, i.e.,
\begin{equation}
	\label{eq:varstrainb}
	\delta 	\ii{\Pi}{}{\mathrm{B,str}} =  \int_{\xi}^{} (N \delta \epsilon_{11} + M^i \delta \kappa_i) J_\mathrm{B} \dd{\xi},
\end{equation}
where section force and moments are
\begin{equation}
	N := \frac{EF}{a} \epsilon_{11}, \quad M^1 := \frac{\mu I_t}{\sqrt{a}} \kappa_1, \quad \quad M^2 := \frac{EI_{\varsigma\varsigma}}{a} \kappa_2, \quad \quad M^3 := \frac{EI_{\eta\eta}}{a} \kappa_3,
\end{equation}
while cross-sectional geometric properties are
\begin{equation}
		F :=\int_{A_\mathrm{B}}^{} \dd{\eta} \dd{\varsigma} , \; \ii{I}{}{\varsigma \varsigma} := \int_{A_\mathrm{B}}^{}  \varsigma ^2 \dd{\eta} \dd{\varsigma} ,  \;  \ii{I}{}{\eta \eta} := \int_{A_\mathrm{B}}^{}  \eta ^2 \dd{\eta} \dd{\varsigma} ,  \;  \ii{I}{}{t} := \int_{A_\mathrm{B}}^{}  \left(\eta ^2 + \varsigma ^2 \right) \dd{\eta} \dd{\varsigma} .
\end{equation}

To completely define the variation of the strain energy in \eqref{eq:varstrainb}, variations of the reference strains are required. The variation of the tangential basis vector is $\delta \ve{a}_1 = \delta \ve{u}_{\mathrm{B},1}$, while the variation of the other two basis vectors, $\ve{a}_\alpha = \pmb{\Lambda}  \cdot\ve{A}_\alpha$, requires variation of the rotation tensor, i.e.,
\begin{equation}
	\delta \ve{a}_\alpha = \delta \pmb{\Lambda} \cdot\ve{A}_\alpha = \delta \ve{w}_\times  \cdot\pmb{\Lambda}  \cdot\ve{A}_\alpha = \delta \ve{w}_\times  \cdot \ve{a}_\alpha = \delta \ve{w} \times \ve{a}_\alpha.
\end{equation}
%
%
where $\delta \ve{w}$ represents a multiplicative variation of the rotation vector. Components of this vector can be determined from conditions $\delta \ve{w} \times \ve{t} = \delta \ve{t}$ and $\delta\ve{w}\cdot\ve{t}=\delta\ve{a}_2^\mathrm{SR}\cdot\ve{a}_3^\mathrm{SR}+\delta \varphi$, i.e.,
\begin{equation}
	\delta \ve{w} =(\delta \varphi - \frac{\ve{T}\cdot(\ve{t}\times\delta\ve{u}_\mathrm{B,1})}{\sqrt{a}\;t_T}) \ve{t} + \frac{\ve{t}\times\delta \ve{u}_{\mathrm{B},1}}{\sqrt{a}} = \ve{t} \delta \varphi + \left(\ve{I} - \frac{\ve{t}\otimes \ve{T}}{t_T}\right) \frac{\ve{t}\times\delta \ve{u}_{\mathrm{B},1}}{\sqrt{a}},
\end{equation}
and the variations of basis vectors $\ve{a}_\alpha$ follow as
\begin{equation}
	\begin{aligned}
		\delta \ve{a}_2 &= \delta \ve{w} \times \ve{a}_2 = \ve{a}_3 \delta \varphi - \frac{1}{\sqrt{a}} \left[\ve{t}\otimes\ve{a}_2 + \frac{\ve{a}_3 \otimes (\ve{T}\times \ve{t})}{t_T}\right] \cdot \delta \ve{u}_{\mathrm{B},1}, \\
		\delta \ve{a}_3 &= \delta \ve{w} \times \ve{a}_3 = -\ve{a}_2 \delta \varphi - \frac{1}{\sqrt{a}} \left[\ve{t}\otimes\ve{a}_3 - \frac{\ve{a}_2 \otimes (\ve{T}\times \ve{t})}{t_T}\right] \cdot \delta \ve{u}_{\mathrm{B},1}.
	\end{aligned}
\end{equation}
Finally, the variations of reference strains \eqref{eq:beamstrain} can be found as
\begin{equation}
	\begin{aligned}
		\delta \epsilon_{11} &= \ve{a}_1 \cdot \delta \ve{u}_{\mathrm{B,1}}, \\ 
		\delta \kappa_1 &= \delta (\ve{a}_{2,1}\cdot\ve{a}_3) = \ve{a}_3 \cdot \delta\ve{a}_{2,1} + \ve{a}_{2,1} \cdot \delta \ve{a}_3, \\
		\delta \kappa_2 &= \delta (-\ve{a}_{1,1}\cdot\ve{a}_3) = -\ve{a}_3 \cdot \delta\ve{u}_{\mathrm{B},11} - \ve{a}_{1,1} \cdot \delta \ve{a}_3), \\
		\delta \kappa_3 &= \delta (\ve{a}_{1,1}\cdot\ve{a}_2) = \ve{a}_2 \cdot \delta\ve{u}_{\mathrm{B},11} + \ve{a}_{1,1} \cdot \delta \ve{a}_2, \\
	\end{aligned}
\end{equation}
which, with a bit of long but straightforward calculation, can be written as
\begin{equation}
	\begin{aligned}
		\delta \kappa_1 &= \delta \varphi_{,1} + \frac{1}{a \, t_T} \left[ \frac{(\ve{t}\cdot\ve{T} )_{,1}(\ve{T}\times\ve{t}) }{t_T} + \Gamma (\ve{T}\times\ve{t}) - (\ve{T}\times\ve{t})_{,1} \right] \cdot \delta \ve{u}_{\mathrm{B},1}, \\
		 &+ \frac{1}{a}\left( K_2 \ve{a}_2 + K_3 \ve{a}_3\right)\cdot \delta \ve{u}_{\mathrm{B},1} -\frac{\ve{T}\times\ve{t}}{a \, t_T} \cdot \delta \ve{u}_{\mathrm{B},11}, \\
		\delta \kappa_2 &= (\Gamma \, \ve{a}_3 - \frac{K_3 (\ve{T}\times \ve{t})}{ t_T}) \cdot \delta\ve{u}_{\mathrm{B},1} - \ve{a}_3 \cdot \delta \ve{u}_{\mathrm{B},11}  + a K_3 \, \delta \varphi, \\
		\delta \kappa_3 &= (-\Gamma \, \ve{a}_2 + \frac{K_2 (\ve{T}\times \ve{t})}{ t_T}) \cdot \delta\ve{u}_{\mathrm{B},1} + \ve{a}_2 \cdot \delta \ve{u}_{\mathrm{B},11}  - a K_2 \, \delta \varphi.  \\
	\end{aligned}
\end{equation}

Regarding the external forces, let us consider the work of the point moment, $\ve{M}$, since it is used in our numerical example. If the external moment acts at point $\xi_c$, the variation of the external work is
\begin{equation}
	\label{eq:varmom}
	\delta \Pi_\mathrm{ext} = \ve{M} \cdot \delta \ve{w} = 
	(\ve{M}\cdot\ve{t}) \delta \varphi -(\ve{M}\cdot\ve{t}) \frac{\ve{T} \times \ve{t}}{\sqrt{a}\; t_T} \cdot \delta \iv{u}{}{\text{B},1} + \frac{\ve{M} \times \ve{t} }{\sqrt{a}} \cdot \delta \iv{u}{}{\text{B},1},
\end{equation}
where all quantities are evaluated at $\xi=\xi_c$. Linearizations of variations of all required quantities are provided in Appendix \ref{appendixb}.

\section{Beam-shell interaction potential}
\label{secbsip}

Having defined the potential energy of the beam and the shell, we now turn to determining their interaction potential, see \eqqref{eq: poten}. In this section, we revisit the disk-infinite plate law and estimate the error made by approximating a disk-curved body with a disk-flat body interaction.

\subsection{Disk-infinite plate law}

\begin{figure}[h!]
	\centering
	\includegraphics[width=0.95\textwidth]{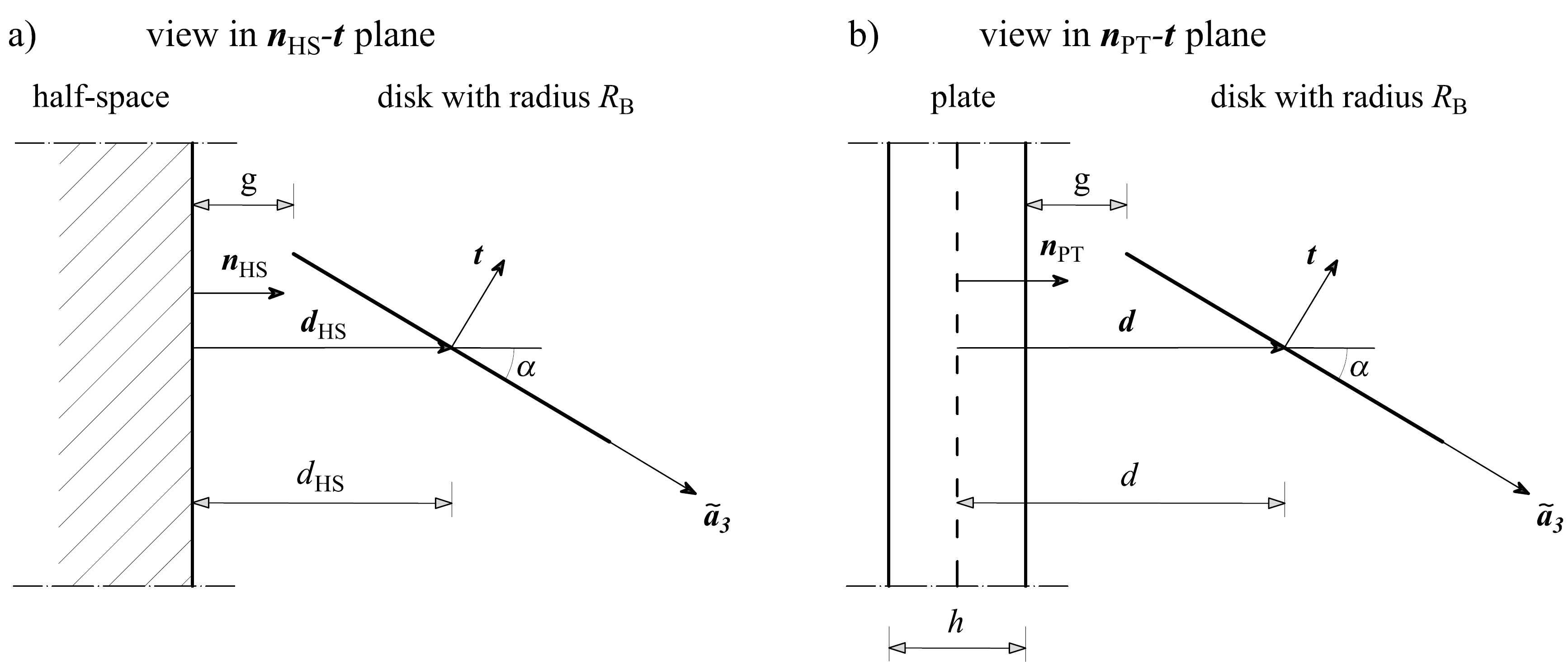} 
	\caption{a) Disk-half-space interaction. b) Disk-plate interaction.}
	\label{fig:FigDHS}
\end{figure}

Let us consider disk-half-space interaction as in Fig.~\ref{fig:FigDHS}a.
We observe the $\ivn{n}{}{HS}$-$\ve{t}$ plane, where $\ivn{n}{}{HS}$ is the normal of the half-space and $\ve{t}$ is the normal of the disk. In this plane, the disk is seen as a line. The disk-half-space law for arbitrary $m$ is derived in \cite{2025borković}. For vdW interaction, this law reduces to
\begin{equation}
	\begin{aligned}
	\Pi_{\operatorname{D-HS}}^6(d_\mathrm{HS},\alpha)&= \frac{\pi^2 R_\mathrm{B}^2}{6 \left(d_\mathrm{HS}^2-R_\mathrm{B}^2 \cos^2 \alpha\right)^{3/2}} 
	, \quad \text{with} \quad d_\mathrm{HS} = g+R_\mathrm{B} \cos \alpha.
	\end{aligned}
	\label{eq:dhsvdw}
\end{equation}
Here, $d_\mathrm{HS}$, $\alpha$, and $g$ are the distance, the angle, and the gap between the disk centroid and the half-space, respectively.

Having the D-HS law, it is straightforward to find the interaction potential between a disk and an infinite-area plate, see Fig.~\ref{fig:FigDHS}b:
\begin{equation}
	\begin{aligned}
		\Pi_{\operatorname{D-PT}}^m (d_\mathrm{PT},\alpha,h)&=  \Pi_{\operatorname{D-HS}}^m \left(d_\mathrm{PT}-h/2,\alpha\right) - \Pi_{\operatorname{D-HS}}^m \left(d_\mathrm{PT}+h/2,\alpha\right),
	\end{aligned}
	\label{eq:DPT}
\end{equation}
where $d_\mathrm{PT}= g+R_\mathrm{B} \cos \alpha + h/2 $ is the distance between the disk's centroid and the plate's middle plane; $h$ is the thickness of the plate, while $\alpha$ and $g$ are the angle and the gap between the disk and the plate, respectively. Note that the D-PT law can be expressed either w.r.t.~the distance $d_\mathrm{PT}$ or w.r.t.~the gap $g$.

An explicit expression for the D-PT interaction via steric repulsion ($m=12$), required for modeling the LJ potential, can be found in Notebook 4 of \cite{2025borkovićb}.

\subsection{Error estimate for approximating disk-curved body interaction by disk-flat body interaction}
\label{sec:error2}

Let us estimate the error we are making by approximating the interaction between a disk and curved geometries with disk-half-space or disk-plate laws. In this subsection, we are exclusively considering $m=6$.

\begin{figure}[h]
	\centering
	\includegraphics[width=\textwidth]{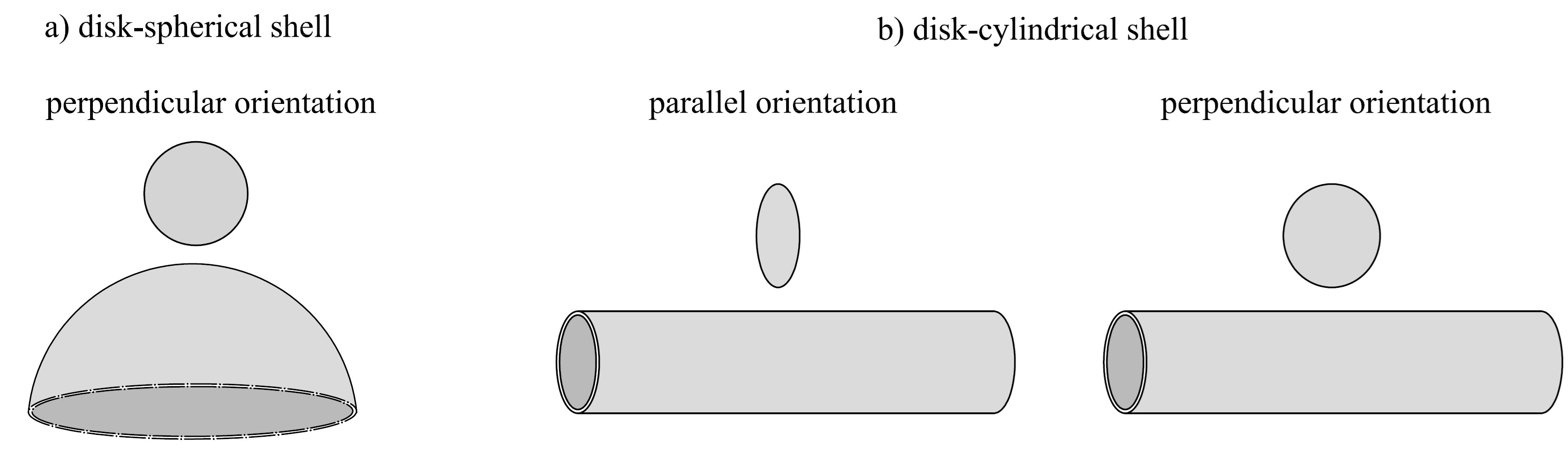}
	\caption{Some types of interaction pairs considered for error estimates.}
	\label{fig:error}
\end{figure}

Consider an interaction between a disk and a sphere (or a spherical shell) that is depicted in Fig.~\ref{fig:error}a as a disk-spherical cap pair. 
Our focus is on a special case of relative orientation between a disk and a sphere, when the disk lies in the plane passing through the sphere's center. We will refer to it as a \emph{perpendicular orientation} because the disk's plane is perpendicular to the sphere's tangent plane at the closest point. For this interaction case, an exact disk-sphere vdW law is derived in Appendix \ref{appendixf}. 

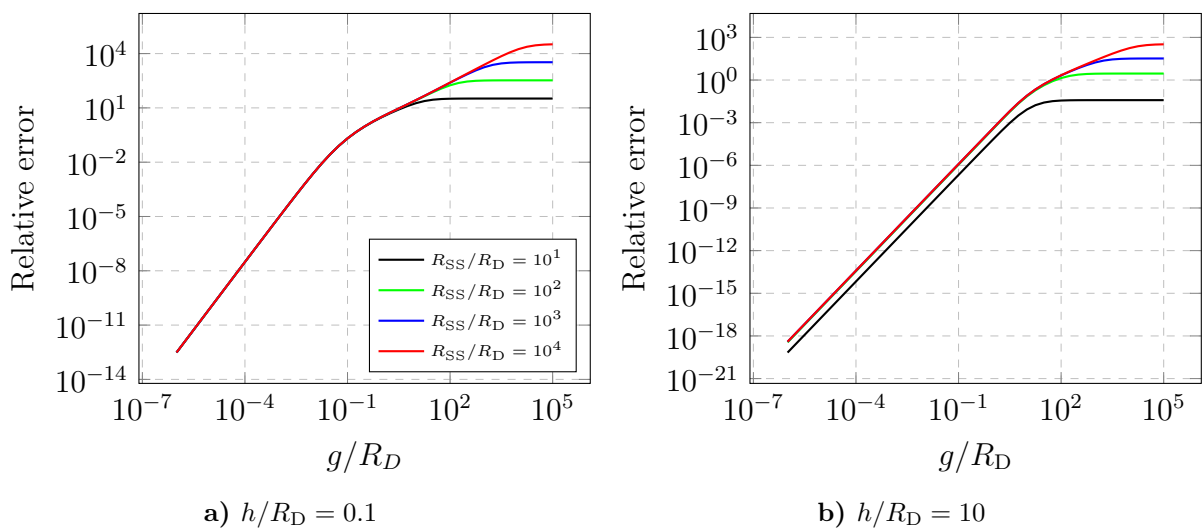
\begin{figure}[htbp]
	\centering
	\subfloat[$h/R_\mathrm{D}=0.1$]{\label{fig:a2}
		\begin{tikzpicture}
			\begin{loglogaxis}[
				xlabel = {$g/R_D$},
				ylabel = {Relative error},
				ylabel near ticks,
				legend pos=south east,
				legend cell align=left,
				legend style={font=\tiny},
				width=0.47\textwidth,
				ytick distance=1000,
				xtick distance=1000,
				clip=true,grid=both,xminorgrids=false,yminorgrids=false];
				
				\addplot[black,thick] table [col sep=comma] {data/relErrorDiskVsSphereSpherShellhRc01RsRc10.csv};
				\addplot[green,thick] table [col sep=comma] {data/relErrorDiskVsSphereSpherShellhRc01RsRc100.csv};
				\addplot[blue,thick] table [col sep=comma] {data/relErrorDiskVsSphereSpherShellhRc01RsRc1000.csv};
				\addplot[red,thick] table [col sep=comma] {data/relErrorDiskVsSphereSpherShellhRc01RsRc10000.csv};
				\legend{$R_\mathrm{SS}/R_\mathrm{D}=10^1$,$R_\mathrm{SS}/R_\mathrm{D}=10^2$,$R_\mathrm{SS}/R_\mathrm{D}=10^3$,$R_\mathrm{SS}/R_\mathrm{D}=10^4$}
			\end{loglogaxis}
		\end{tikzpicture}	
	}
	\subfloat[$h/R_\mathrm{D}=10$]{\label{fig:b1}
		\begin{tikzpicture}
			\begin{loglogaxis}[
				xlabel = {$g/R_\mathrm{D}$},
				ylabel = {Relative error},
				ylabel near ticks,
				legend cell align=left,
				legend style={font=\tiny},
				legend style={at={(0.025,0.85)},anchor=west},
				width=0.47\textwidth,
				ytick distance=1000,
				xtick distance=1000,
				clip=true,grid=both,xminorgrids=false,yminorgrids=false];
				\addplot[black,thick] table [col sep=comma] {data/relErrorDiskVsSphereSpherShellhRc10RsRc10.csv};
				\addplot[green,thick] table [col sep=comma] {data/relErrorDiskVsSphereSpherShellhRc10RsRc100.csv};
				\addplot[blue,thick] table [col sep=comma] {data/relErrorDiskVsSphereSpherShellhRc10RsRc1000.csv};
				\addplot[red,thick] table [col sep=comma] {data/relErrorDiskVsSphereSpherShellhRc10RsRc10000.csv};
			\end{loglogaxis}
	\end{tikzpicture}	}
	\caption{Relative error of the total force between disk-sphere and disk-spherical shell (see Fig.~\ref{fig:error}a).}
	\label{fig:DS-DSS}
\end{figure}

Let us scrutinize the difference between the disk-sphere (D-S) and the disk-spherical shell (D-SS) interaction. The relative difference in the total interaction force is displayed in Fig.~\ref{fig:DS-DSS} for varying values of the three ratios: $R_\mathrm{SS}/R_\mathrm{D}$, $h/R_\mathrm{D}$, and $g/R_\mathrm{D}$. Here, $R_\mathrm{SS}$ is the radius of the midsurface of the spherical shell, meaning that the comparable sphere has the radius of $R_\mathrm{SS}+h/2$.
The difference is large for large gaps, but decreases as the gap narrows. Increase in the ratio $h/R_\mathrm{D}$ shifts the curves to the right so that the differences decrease faster with the gap for large $h/R_\mathrm{D}$. Since we have analytical expressions for both D-S and D-SS interactions, we can show that the limit value of the relative difference for $g\rightarrow0$ is indeed zero in this case, i.e.,
\begin{equation}
	\lim_{g\rightarrow0} \frac{\Pi_{\mathrm{D-SS}}-\Pi_{\mathrm{D-S}}}{\Pi_{\mathrm{D-SS}}} = 0.
\end{equation}

\begin{figure}[htbp]
	\centering
	\subfloat[$h/R_\mathrm{D}=0.1$]{\label{fig:a2cs}
		\begin{tikzpicture}
			\begin{loglogaxis}[
				xlabel = {$g/R_\mathrm{D}$},
				ylabel = {Relative error},
				ylabel near ticks,
				legend cell align=left,
				legend style={font=\tiny},
				legend style={at={(0.025,0.85)},anchor=west},
				width=0.47\textwidth,
				ytick distance=100,
				xtick distance=100,
				clip=true,grid=both,xminorgrids=false,yminorgrids=false];
				
				\addplot[black,thick] table [col sep=comma] {data/relErrorDiskPTvsDiskSpherShellhRc01RsRc10.csv};
				\addplot[green,thick] table [col sep=comma] {data/relErrorDiskPTvsDiskSpherShellhRc01RsRc100.csv};
				\addplot[blue,thick] table [col sep=comma] {data/relErrorDiskPTvsDiskSpherShellhRc01RsRc1000.csv};
				\addplot[red,thick] table [col sep=comma] {data/relErrorDiskPTvsDiskSpherShellhRc01RsRc10000.csv};
				\legend{$R_\mathrm{SS}/R_\mathrm{D}=10^1$,$R_\mathrm{SS}/R_\mathrm{D}=10^2$,$R_\mathrm{SS}/R_\mathrm{D}=10^3$,$R_\mathrm{SS}/R_\mathrm{D}=10^4$}
			\end{loglogaxis}
		\end{tikzpicture}	
	}
	\subfloat[$h/R_\mathrm{D}=10$]{\label{fig:b1cs}
		\begin{tikzpicture}
			\begin{loglogaxis}[
				xlabel = {$g/R_\mathrm{D}$},
				ylabel = {Relative error},
				ylabel near ticks,
				legend cell align=left,
				legend style={font=\tiny},
				legend style={at={(0.025,0.85)},anchor=west},
				width=0.47\textwidth,
				ytick distance=100,
				xtick distance=100,
				clip=true,grid=both,xminorgrids=false,yminorgrids=false];
				\addplot[black,thick] table [col sep=comma] {data/relErrorDiskPTvsDiskSpherShellhRc10RsRc10.csv};
				\addplot[green,thick] table [col sep=comma] {data/relErrorDiskPTvsDiskSpherShellhRc10RsRc100.csv};
				\addplot[blue,thick] table [col sep=comma] {data/relErrorDiskPTvsDiskSpherShellhRc10RsRc1000.csv};
				\addplot[red,thick] table [col sep=comma] {data/relErrorDiskPTvsDiskSpherShellhRc10RsRc10000.csv};
			\end{loglogaxis}
	\end{tikzpicture}	}
	\caption{Relative error of the total force between a disk and a spherical shell (see Fig.~\ref{fig:error}a) approximated with the D-PT interaction.}
	\label{fig:DSS-DPT}
\end{figure}
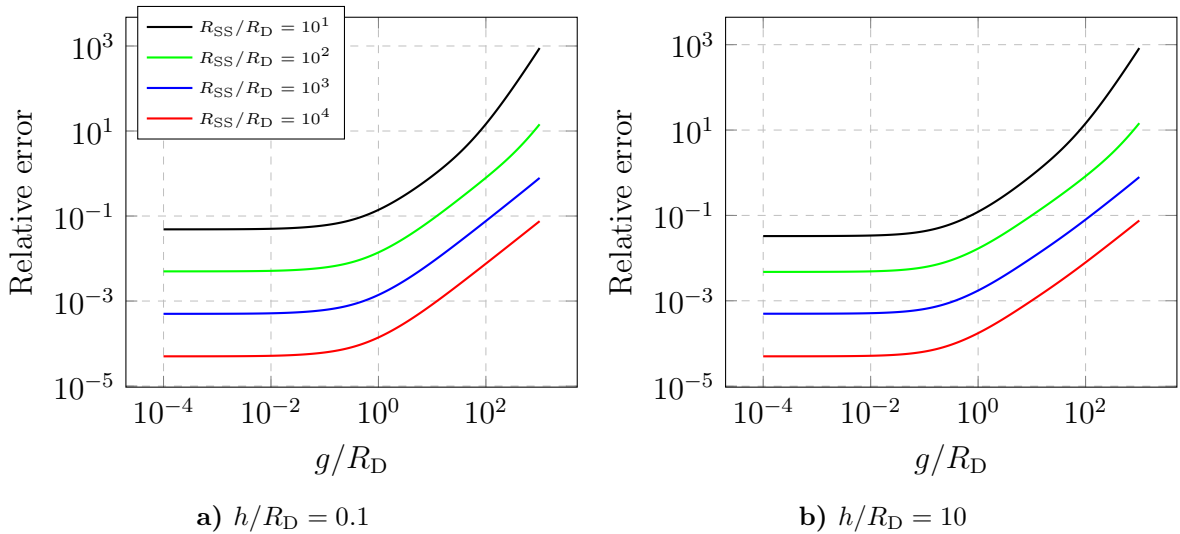

The next step is to approximate the D-SS interaction with the D-PT interaction. The relative error of this approximation is shown in Fig.~\ref{fig:DSS-DPT}.
The error again reduces with the gap, but it converges to a finite value. The analytical solution for this error is
\begin{equation}
	\label{eq:errorestim}
	\lim_{g\rightarrow0} \frac{\Pi_{\mathrm{D}-\mathrm{SS}}-\Pi_{\mathrm{D-PT}}}{\Pi_{\mathrm{D-SS}}} = 1-\sqrt{1+\frac{2 R_\mathrm{D}}{h+2R_\mathrm{SS}}}.
\end{equation}
%
%
Therefore, the accuracy improves for thick (strongly curved) shells with large $h$, while the accuracy reduces for thin shells. The error of this approximation represents an upper bound of the error for approximating the disk-shell interaction with the disk-plate interaction.

\begin{figure}[h!]
	\centering
	\subfloat[$h/R_\mathrm{D}=0.1$]{\label{fig:as2c}
		\begin{tikzpicture}
			\begin{loglogaxis}[
				xlabel = {$g/R_\mathrm{D}$},
				ylabel = {Relative error},
				ylabel near ticks,
				legend cell align=left,
				legend style={font=\tiny},
				legend style={at={(0.025,0.85)},anchor=west},
				width=0.47\textwidth,
				ytick distance=100,
				xtick distance=100,
				clip=true,grid=both,xminorgrids=false,yminorgrids=false];
				
				\addplot[black,thick] table [col sep=comma] {data/relErrorOfDPTvsDiskCylinShellParallelhRc01RsRc10.csv};
				\addplot[green,thick] table [col sep=comma] {data/relErrorOfDPTvsDiskCylinShellParallelhRc01RsRc100.csv};
				\addplot[blue,thick] table [col sep=comma] {data/relErrorOfDPTvsDiskCylinShellParallelhRc01RsRc1000.csv};
				\addplot[red,thick] table [col sep=comma] {data/relErrorOfDPTvsDiskCylinShellParallelhRc01RsRc10000.csv};
				\legend{$R_\mathrm{C}/R_\mathrm{D}=10^1$,$R_\mathrm{C}/R_\mathrm{D}=10^2$,$R_\mathrm{C}/R_\mathrm{D}=10^3$,$R_\mathrm{C}/R_\mathrm{D}=10^4$}
			\end{loglogaxis}
		\end{tikzpicture}	
	}
	\subfloat[$h/R_D=10$]{\label{fig:sb1c}
		\begin{tikzpicture}
			\begin{loglogaxis}[
				xlabel = {$g/R_D$},
				ylabel = {Relative error},
				ylabel near ticks,
				legend cell align=left,
				legend style={font=\tiny},
				legend style={at={(0.025,0.85)},anchor=west},
				width=0.47\textwidth,
				ytick distance=100,
				xtick distance=100,
				clip=true,grid=both,xminorgrids=false,yminorgrids=false];
				\addplot[black,thick] table [col sep=comma] {data/relErrorOfDPTvsDiskCylinShellParallelhRc10RsRc10.csv};
				\addplot[green,thick] table [col sep=comma] {data/relErrorOfDPTvsDiskCylinShellParallelhRc10RsRc100.csv};
				\addplot[blue,thick] table [col sep=comma] {data/relErrorOfDPTvsDiskCylinShellParallelhRc10RsRc1000.csv};
				\addplot[red,thick] table [col sep=comma] {data/relErrorOfDPTvsDiskCylinShellParallelhRc10RsRc10000.csv};
			\end{loglogaxis}
	\end{tikzpicture}	}
	\caption{Relative error of the total force between a disk and a cylindrical shell in parallel orientation (see Fig.~\ref{fig:error}b) approximated with the D-PT interaction.}
	\label{fig:DCII-DPT}
\end{figure}
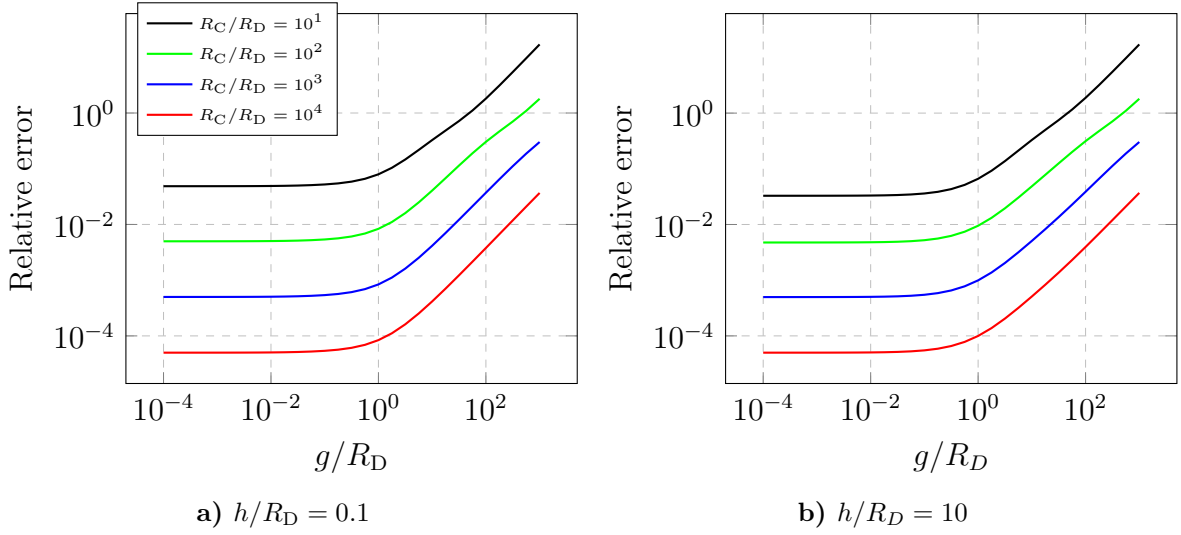

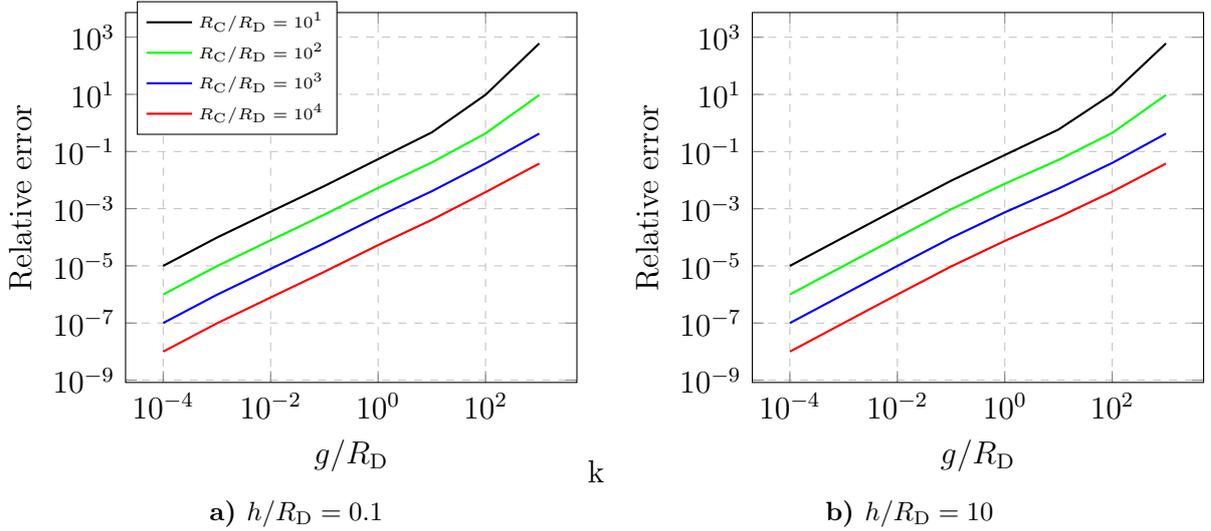
\begin{figure}[h!]
	\centering
	\subfloat[$h/R_\mathrm{D}=0.1$]{\label{fig:a2c}
		\begin{tikzpicture}
			\begin{loglogaxis}[
				xlabel = {$g/R_\mathrm{D}$},
				ylabel = {Relative error},
				ylabel near ticks,
				legend cell align=left,
				legend style={font=\tiny},
				legend style={at={(0.025,0.85)},anchor=west},
				width=0.47\textwidth,
				ytick distance=100,
				xtick distance=100,
				clip=true,grid=both,xminorgrids=false,yminorgrids=false];
				
				\addplot[black,thick] table [col sep=comma] {data/relErrorOfDPTvsDiskCylinShellOrthohRc01RsRc10.csv};
				\addplot[green,thick] table [col sep=comma] {data/relErrorOfDPTvsDiskCylinShellOrthohRc01RsRc100.csv};
				\addplot[blue,thick] table [col sep=comma] {data/relErrorOfDPTvsDiskCylinShellOrthohRc01RsRc1000.csv};
				\addplot[red,thick] table [col sep=comma] {data/relErrorOfDPTvsDiskCylinShellOrthohRc01RsRc10000.csv};
				\legend{$R_\mathrm{C}/R_\mathrm{D}=10^1$,$R_\mathrm{C}/R_\mathrm{D}=10^2$,$R_\mathrm{C}/R_\mathrm{D}=10^3$,$R_\mathrm{C}/R_\mathrm{D}=10^4$}
			\end{loglogaxis}
		\end{tikzpicture}	
k	}
	\subfloat[$h/R_\mathrm{D}=10$]{\label{fig:b1c}
		\begin{tikzpicture}
			\begin{loglogaxis}[
				xlabel = {$g/R_\mathrm{D}$},
				ylabel = {Relative error},
				ylabel near ticks,
				legend cell align=left,
				legend style={font=\tiny},
				legend style={at={(0.025,0.85)},anchor=west},
				width=0.47\textwidth,
				ytick distance=100,
				xtick distance=100,
				clip=true,grid=both,xminorgrids=false,yminorgrids=false];
				\addplot[black,thick] table [col sep=comma] {data/relErrorOfDPTvsDiskCylinShellOrthohRc10RsRc10.csv};
				\addplot[green,thick] table [col sep=comma] {data/relErrorOfDPTvsDiskCylinShellOrthohRc10RsRc100.csv};
				\addplot[blue,thick] table [col sep=comma] {data/relErrorOfDPTvsDiskCylinShellOrthohRc10RsRc1000.csv};
				\addplot[red,thick] table [col sep=comma] {data/relErrorOfDPTvsDiskCylinShellOrthohRc10RsRc10000.csv};
			\end{loglogaxis}
	\end{tikzpicture}	}
	\caption{Relative error of the total force between a disk and a cylindrical shell in perpendicular orientation (see Fig.~\ref{fig:error}b) approximated with the D-PT interaction.}
	\label{fig:DCT-DPT}
\end{figure}

Finally, let us consider a special case of disk-curved geometry interaction: a disk-cylindrical shell interaction in two orientations, see Fig.~\ref{fig:error}b. In the first orientation, the disk and the cylindrical shell are parallel, meaning that the cylinder's cross section is parallel to the disk. This is a well-known case for which no analytical solution exists for the vdW interaction. Therefore, we find the reference solution by numerically integrating the exact point-cylinder law. The relative error of approximating disk-cylindrical shell in parallel orientation with the D-PT law is displayed in Fig.~\ref{fig:DCII-DPT}.
The results are almost the same as those for the disk-spherical shell case in Fig.~\ref{fig:DSS-DPT}. Second, let us consider a disk-cylindrical shell interaction in perpendicular orientation, meaning that the disk and the cylinder's axis lie in the same plane, see Fig.~\ref{fig:error}b. The relative error of approximating such interaction with D-PT law is shown in Fig.~\ref{fig:DCT-DPT}.
For this case of relative orientation between a disk and a cylindrical shell, the error converges to zero with the gap.

To summarize, we have considered two cases of the disk-cylindrical shell interaction. They can be distinguished by examining the section that forms when the gap becomes negative. If such a section is curved, this corresponds to the upper bound of the approximation error for $g\rightarrow0$, estimated with \eqqref{eq:errorestim}. On the other hand, when their section is flat, the error goes to zero with the gap, see Fig.~\ref{fig:DCT-DPT}, which represents the lower bound of the error of our approximation.

\subsection{Error estimate for approximating the cylinder-spherical shell interaction with SBS approach}
\label{subs:css}

After estimating the error in the interaction between a disk and various curved bodies in the previous subsection, we now take the next step: integrating the D–PT law along the beam axis to obtain the SBS interaction and assess the corresponding error. Several configurations can be considered.
The most general case is the interaction between a curved beam and a curved shell. However, a numerical integration of such a potential that could serve as a reference solution is computationally prohibitive. A special case is the interaction between a curved beam and a flat shell (plate). This case is exactly described with the SBS approach, since the main assumption of our approximation (that the shell at the closest point can be replaced with an infinite plate) is then fully satisfied. Another special case is the straight beam–curved shell interaction, which is examined here.

\begin{figure}[h]
	\centering
	\includegraphics[width=0.8\textwidth]{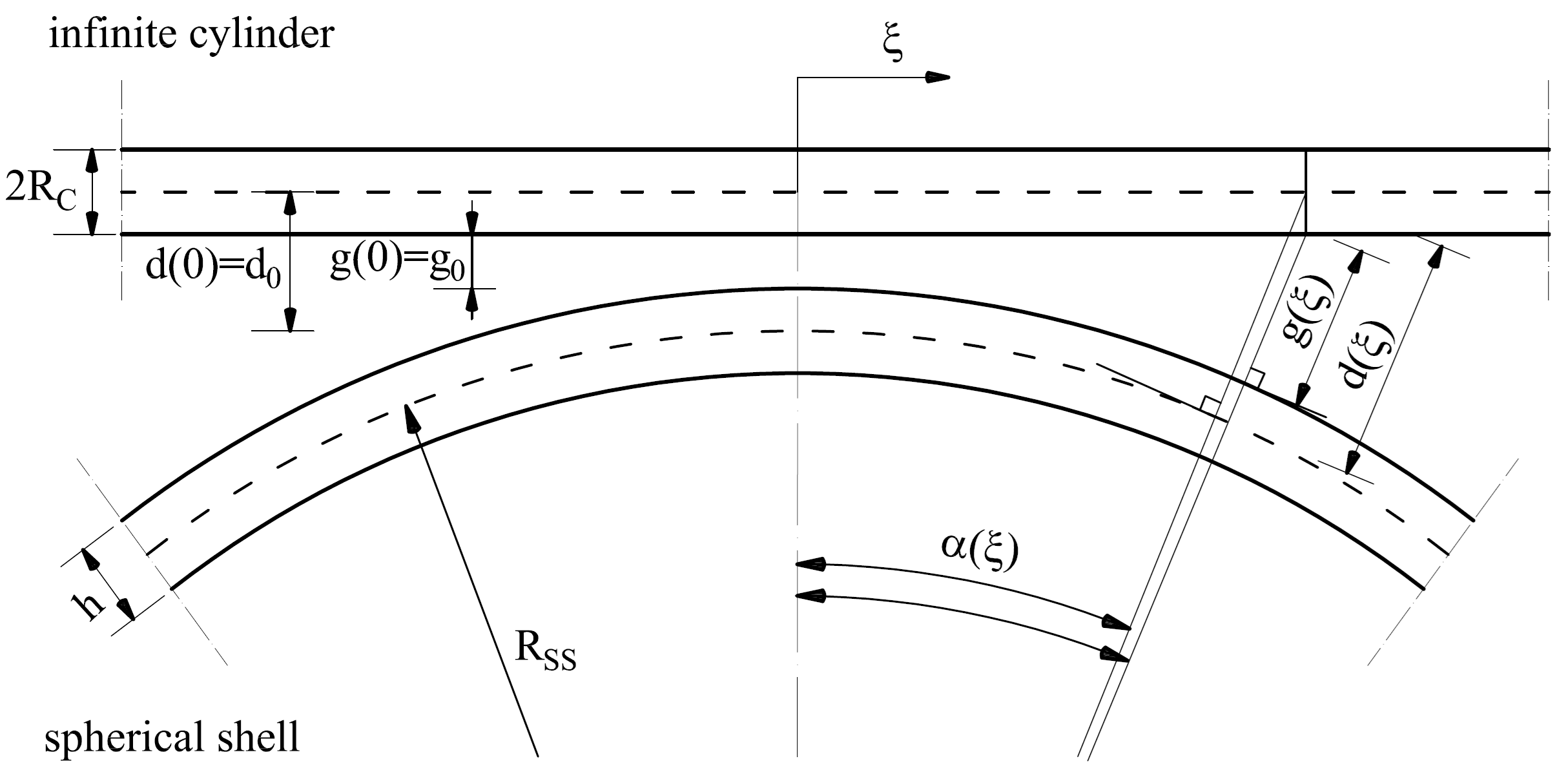} 
	\caption{Infinite cylinder-spherical shell interaction (C-SS).}
	\label{fig:C-SS}
\end{figure}

An interaction between an infinite cylinder of cross-sectional radius $R_\mathrm{C}$ and a spherical shell of thickness $h$ and midsurface radius $R_{\mathrm{SS}}$ (C-SS) is sketched in Fig.~\ref{fig:C-SS}.
The smallest gap between the cylinder's cross-section and the shell is marked with $g_0$. As a reference, we use the analytical solution for the vdW force derived by Kirsch \cite{2003kirsch}.

\begin{figure}[hbt!]
	\centering
	\begin{tikzpicture}
		\begin{axis}[
			xlabel = {Log $g_0$},
			ylabel = {Scaling factor $S(g_0)$},
			ylabel near ticks,
			legend pos=north west,
			legend cell align=left,
			legend style={font=\scriptsize},
			width=0.5\textwidth,
			xmin = -8, xmax = 8,
			clip=true,grid=both];
			\addplot[red,thick] table [col sep=comma] {data/scalingKirsch.csv};
			\addplot[blue,dashed,thick] table [col sep=comma] {data/scalingSurogat.csv};
			\legend{Exact - Kirsch, Approximate - SBS} 
		\end{axis}
	\end{tikzpicture}
	\caption{Scaling factor for the infinite cylinder-spherical shell interaction force. Comparison of the exact values with the SBS approximation. ($R_{\mathrm{SS}}=10, R_\mathrm{C}=1, h=1$)}
	\label{fig:LJpoten4xaaNEW3}
\end{figure}
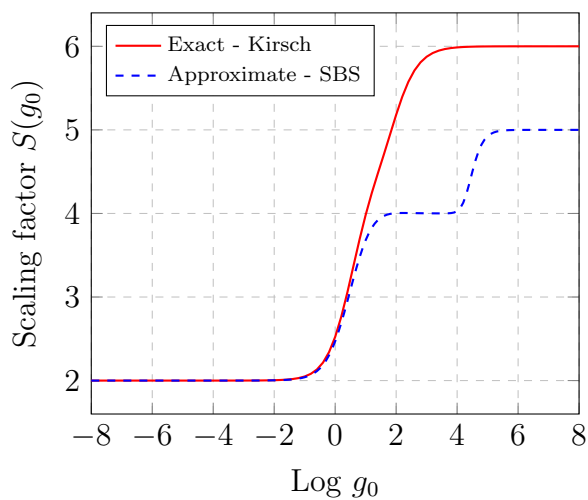

First, we compare the exact and approximate scaling factors in Fig.~\ref{fig:LJpoten4xaaNEW3}.
The most important observation is that the SBS approach provides the correct scaling of the interaction force w.r.t.~the gap for small separations, 2. Even for moderate separations, $g_0\approx R_\mathrm{C}$, the scaling is within the limits of acceptable accuracy. As the separation increases, $g_0 > R_\mathrm{C}$, the C-SS interaction reduces to the point-line interaction with the scaling factor of 6. The SBS approach fails to describe this long-range scaling; it first converges to 4 and then to 5. 
In essence, this analysis demonstrates that the SBS approach correctly captures the scaling in the decisive regime of small separations.

\begin{remark}
	\normalfont
Let us define the true scaling factor function, $\bar S(p)$, of the potential $\Pi (p)$, where $p$ is the gap \cite{2025borković,2025borkoviće}. We assume that the potential $\Pi(p)$ is a continuous function of the gap $p$ in the domain of interest, and that the new value of potential, $\Pi (\zeta a)$, due to a small change of the gap, $\omega p$, with $\omega \approx 1$, can be represented as $\omega ^{\bar S(p)} \Pi (p)$. For the limit case $\omega \rightarrow 1$, both representations should give the same value, i.e.,
	\begin{equation}
		\label{eqscalinga1}
		\begin{aligned}
			\lim_{\omega\rightarrow1}	\big[\Pi( \omega p) - \omega ^{\bar S(p)} \Pi (p)\big] = 0.
		\end{aligned}
	\end{equation}
Furthermore, these representations must have the same first derivative w.r.t.~$\omega$ for $\omega \rightarrow 1$
	\begin{equation}
		\label{eqscalingas1}
		\begin{aligned}
			\lim_{\omega \rightarrow1}	\frac{\partial \big[\Pi(\omega p) - \omega ^{\bar S(p)} \Pi (p)\big]}{\partial \omega} = 0,
		\end{aligned}
	\end{equation}	
which allows us to find an explicit expression for the true scaling factor	
	\begin{equation}
		\label{eqscalingass1}
		\lim_{\omega \rightarrow1}	\big[ \frac{\partial \Pi(\omega p) }{\partial \omega} - \bar S(p) \omega^{\bar S(p)-1} \Pi (p) \big]= 0 \implies \bar S(p) = \frac{1}{\Pi (p)} \frac{\partial \Pi (\omega p)}{\partial \omega} \vert_{\omega=1} =  \frac{p}{\Pi (p)} \frac{\partial \Pi (p)}{\partial p}.
	\end{equation}
Due to inverse power laws considered, the potential decreases with the increase in gap, meaning that the true scaling factor is negative. For convenience, we refer to the scaling factor as the negative value of the true scaling factor, $S(p)=-\bar S(p)$. 
\end{remark}

Next, let us consider the accuracy of the SBS approximation in this case. The complete D-PT expression in \eqqref{eq:DPT} takes into account the orientation of the cross section w.r.t.~to the plate, $\cos \alpha (\xi)$. Depending on the choice of representing the D-PT law either via gap $g(\xi)$ or via distance $d(\xi)$, a different closest point and angle $\alpha (\xi)$ are obtained, see Fig.~\ref{fig:C-SS}. Nevertheless, our numerical analysis shows that both choices return practically indistinguishable results. This is because the sections closest to the shell, near $\xi=0$, dominate the total interaction, and this ambiguity is negligible in these regions. In the following, we are comparing the complete D-PT law \eqref{eq:DPT} with an approximate one, where the orientation is neglected by setting $\cos \alpha (\xi)=1$.

\begin{figure}[htbp]
	\centering
	\subfloat[$h/R_\mathrm{C}=0.1$]{\label{fig:a}
		\begin{tikzpicture}
				\begin{loglogaxis}[
						xlabel = {$g_0/R_\mathrm{C}$},
						ylabel = {Relative error},
						ylabel near ticks,
						legend cell align=left,
						legend style={font=\tiny},
						legend style={at={(0.025,0.85)},anchor=west},
						width=0.47\textwidth,
						ytick distance=10,
						xtick distance=10,
						clip=true,grid=both,xminorgrids=false,yminorgrids=false];
						
							\addplot[black,thick] table [col sep=comma] {data/relErrorSurogatExacthRc01RsRc10.csv};
						\addplot[green,thick] table [col sep=comma] {data/relErrorSurogatExacthRc01RsRc100.csv};
						\addplot[blue,thick] table [col sep=comma] {data/relErrorSurogatExacthRc01RsRc1000.csv};
						\addplot[red,thick] table [col sep=comma] {data/relErrorSurogatExacthRc01RsRc10000.csv};
						
						\addplot[black,thick,dashed] table [col sep=comma] {data/relErrorSurogatApproxhRc01RsRc10.csv};
						\addplot[green,thick,dashed] table [col sep=comma] {data/relErrorSurogatApproxhRc01RsRc100.csv};
						\addplot[blue,thick,dashed] table [col sep=comma] {data/relErrorSurogatApproxhRc01RsRc1000.csv};
						\addplot[red,thick,dashed] table [col sep=comma] {data/relErrorSurogatApproxhRc01RsRc10000.csv};

						\legend{$R_{\mathrm{SS}}/R_\mathrm{C}=10^1$,$R_{\mathrm{SS}}/R_\mathrm{C}=10^2$,$R_{\mathrm{SS}}/R_\mathrm{C}=10^3$,$R_{\mathrm{SS}}/R_\mathrm{C}=10^4$}
					\end{loglogaxis}
			\end{tikzpicture}
		}
	\subfloat[$h/R_\mathrm{C}=10$.]{\label{fig:b}
		\begin{tikzpicture}
				\begin{loglogaxis}[
						xlabel = {$g_0/R_\mathrm{C}$},
						ylabel = {Relative error},
						ylabel near ticks,
						legend cell align=left,
						legend style={font=\tiny},
					    legend style={at={(0.025,0.85)},anchor=west},
						width=0.47\textwidth,
						ytick distance=10,
					xtick distance=10,
						clip=true,grid=both,xminorgrids=false,yminorgrids=false];

					\addplot[black,thick] table [col sep=comma] {data/relErrorSurogatExacthRc10RsRc10.csv};
					\addplot[green,thick] table [col sep=comma] {data/relErrorSurogatExacthRc10RsRc100.csv};
					\addplot[blue,thick] table [col sep=comma] {data/relErrorSurogatExacthRc10RsRc1000.csv};
					\addplot[red,thick] table [col sep=comma] {data/relErrorSurogatExacthRc10RsRc10000.csv};
					
						\addplot[black,thick,dashed] table [col sep=comma] {data/relErrorSurogatApproxhRc10RsRc10.csv};
					\addplot[green,thick,dashed] table [col sep=comma] {data/relErrorSurogatApproxhRc10RsRc100.csv};
					\addplot[blue,thick,dashed] table [col sep=comma] {data/relErrorSurogatApproxhRc10RsRc1000.csv};
					\addplot[red,thick,dashed] table [col sep=comma] {data/relErrorSurogatApproxhRc10RsRc10000.csv};
					\end{loglogaxis}
			\end{tikzpicture}
		}
	\caption{Relative error of the total force between an infinite cylinder and a spherical shell using SBS. The error of the exact D-PT law (considering the orientation of the cross section) is plotted with solid lines, while the error of the approximate law (neglecting the orientation of the cross section) is plotted with dashed lines.}
	\label{fig:IC-SS}
\end{figure}
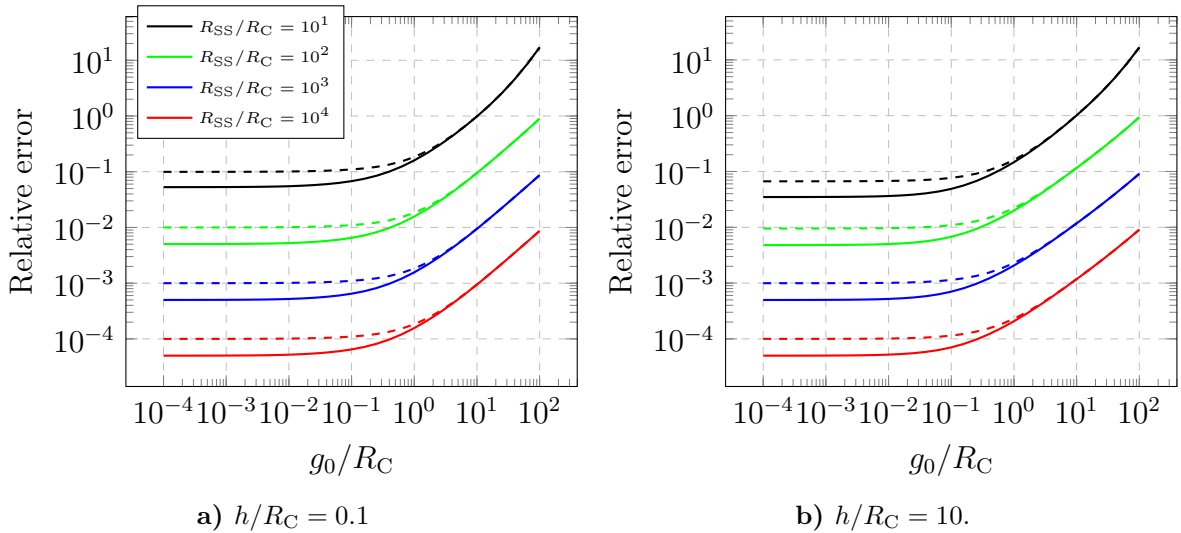

In Fig.~\ref{fig:IC-SS} we plot the relative error of SBS w.r.t.~the analytical solution for different values of the ratios $R_{\mathrm{SS}}/R_\mathrm{C}$, $h/R_\mathrm{C}$ and $g_0/R_\mathrm{C}$.
%
First, we can observe that the exact D-PT law implemented in the SBS approach provides the same error as the one estimated for the disk-spherical shell interaction. Therefore, the upper bound of the relative error for the SBS approach is well-approximated by \eqqref{eq:errorestim}. Second, the orientation of the cross section w.r.t.~the tangential plane influences the approximation error. Third, as seen in the previous subsection, the error significantly reduces for large $R_{\mathrm{SS}}/R_\mathrm{C}$ values. Fourth, the influence of the $h/R_\mathrm{C}$ ratio is noticeable only for strongly curved shells.

The relative errors for large separations, $g_0/R_C>1$ in Fig.~\ref{fig:IC-SS}, are not significant for problems dominated by small separations, such as peeling and adhesion. For other cases, such as an attraction between far-distanced bodies, the accuracy at large separations is crucial. However, this effect is not the subject of the present paper.



To summarize, the proposed SBS provides good accuracy for small separations. The upper bound of the error drops below $1\%$ for $R_{\mathrm{SS}}/R_\mathrm{C}>10$. Considering the significant gain in computational efficiency for the case of interaction between deformable beams and shells (by reducing 6D to 1D integration), we find that the upper and the lower bounds of the SBS error, as defined in this section, are very good.

\section{Variation of the SBS potential}

\label{sec:varsbs}

In this section, the D-PT potential is varied in the context of the SBS approach. While the D-HS and D-PT formulations have already benn derived in \cite{2025borković}, the SBS introduces a specific feature: the gap is determined from the unilateral closest-point condition, see Fig.~\ref{fig:FigBS}. For clarity, we will designate such a potential as the disk-surrogate plate (D-SPT) potential, $\phi$. This potential is a function of two quantities: (i) the distance between the beam's axis and the shell's midsurface, $d$, or the gap between outer surfaces, $g$, and (ii) the cosine of the angle between the disk and the surrogate plate, $\cos\alpha$, see Figs.~\ref{fig:FigDHS}b and \ref{fig:C-SS}. As already noted in Subsection \ref{subs:css}, the difference between approaches that employ $d$ and $g$ exists, but it is negligible. Moreover, employing the gap $g$ as the variable implies finding the closest-point projection between the section's outer perimeter and the shell, which complicates the formulation without a significant improvement in accuracy. Therefore, we consider the D-SPT potential as a function of the distance and the cosine of the angle, i.e., $\phi=\phi(d,\cos\alpha)$.

\subsection{Variation of the distance and the angle}

In the following, the hat symbol marks a quantity evaluated at the closest-point projection of the section's center onto a shell. The distance vector between the point on a beam axis, $\ivn{r}{}{B} (\xi)$, and the corresponding closest point on a shell, $\ivn{r}{}{S} (\hat \theta^\alpha)$, 
\begin{equation}
	\ve{d}(\ivn{r}{}{B} (\xi), \ivn{r}{}{S} (\hat \theta^\alpha)) = \veh{ d}= \ivn{r}{}{B}(\xi) - \ve{r}_\text{S} (\hat \theta^\alpha) = \ivn{r}{}{B} -  \ivnh{r}{}{S}, 
\end{equation}
is defined by the unilateral closest-point condition
\begin{equation}
	\label{eq:cpc}
	\begin{aligned}
		\veh{d} \cdot \veh{g}_\alpha &=	 0,
	\end{aligned}
\end{equation}
where $\veh{g}_\alpha$ are tangential basis vectors of the shell midsurface at the closest point, see Fig.~\ref{fig:FigBS}. The shell normal is colinear with the distance vector at the closest point, i.e.,
\begin{equation}
	\veh{g}_3 = \veh{n} = \frac{\veh{d}}{d}, \quad d =\sqrt{\veh{d}\cdot\veh{d}}.
\end{equation}
To find the variation of the potential, let us define the derivative w.r.t.~the parametric coordinate at the closest point as
\begin{equation}
	\frac{\partial  (.)}{\partial \hat\theta^\alpha} := \frac{\partial  (.)}{\partial \theta^\alpha} \bigg\rvert_{\theta^\alpha=\hat\theta^\alpha}.
\end{equation}
The variation of the distance is
\begin{equation}
	\delta d  = \veh{n} \cdot \delta \veh{d},
	\label{eq:vard}
\end{equation}
while the variation of the distance vector can be represented as
\begin{equation}
\begin{aligned}
		\delta \veh{d}&= \frac{\partial  \veh{d}}{\partial \ivn{r}{}{\text{B}}} \cdot \delta \ivn{r}{}{B}+ \frac{\partial  \veh{d}}{\partial \ivnh{r}{}{\text{S}}}  \cdot \delta \ivnh{r}{}{\text{S}} + \frac{\partial  \veh{d}}{\partial \hat{\theta}^\alpha}  \; \delta \hat{\theta}^\alpha = \delta \iv{ r}{}{\text{B}} - \delta \ivh{r}{}{\text{S}} -\ivh{g}{}{\alpha} \; \delta \hat{\theta}^\alpha, \\
		 \frac{\partial  \veh{d}}{\partial \ivn{r}{}{\text{B}}}  &:= \ve{I},\quad \frac{\partial  \veh{d}}{\partial \ivnh{r}{}{\text{S}}} := -\ve{I}, \quad \frac{\partial  \veh{d}}{\partial \hat{\theta}^\alpha} := -\iv{\hat g}{}{\alpha}.
		 \label{eq:deltad}
\end{aligned}
\end{equation}
By inserting \eqqref{eq:deltad} in \eqqref{eq:vard}, the variation of the distance becomes
\begin{equation}
	\delta d = \veh{n} \cdot (\delta \iv{ r}{}{\text{B}} -\delta \ivh{r}{}{\text{S}} - \ivh{g}{}{\alpha} \;\delta \hat{\theta}^\alpha).
\end{equation}
Finally, due to orthogonality condition \eqref{eq:cpc}, we obtain
\begin{equation}
	\delta d =  \veh{n} \cdot (\delta \iv{ r}{}{\text{B}} - \delta \ivh{r}{}{\text{S}} ) = \veh{n} \cdot (\delta \iv{u}{}{\text{B}} - \delta \ivh{u}{}{\text{S}} ),
\end{equation}
where $\delta \iv{ r}{}{\text{B}} = \delta \iv{u}{}{\text{B}}$ and $\delta \ivh{r}{}{\text{S}}=\delta \ivh{u}{}{\text{S}}$. A specific feature of the present formulation is that the distance depends not only on the positions of the interacting objects but also on the closest-point condition. In other words, it is not only the position of the shell that changes, but also the position of the closest point. Therefore, the variation w.r.t.~the shell depends on the position of the shell and the coordinates $\hat \theta^\alpha$ of the closest point. We must first vary quantities w.r.t.~the configuration for a fixed $\hat\theta^\alpha$, and then w.r.t.~$\hat \theta^\alpha$ itself.

Regarding the angle between the cross section and the surrogate plate, its cosine is represented in \cite{2025borković} as a function of the cosine of its complementary angle, i.e.,
\begin{equation}
	\cos\alpha=\sqrt{1-t_n^2},	 \quad t_n:=\ve{t} \cdot \hat{\ve{n}}=\cos\bar\alpha, \quad \text{with} \quad \alpha+\bar\alpha=\pi/2.
\end{equation}
Therefore, the variation of $\cos\alpha$ is
\begin{equation}
	\delta \cos \alpha =- \frac{t_n}{\cos\alpha}  (\hat{\ve{ n}} \cdot \delta \ve{t} + \ve{t} \cdot \delta \hat{\ve{ n}}).
\end{equation}

\subsection{Variation of the disk-surrogate plate potential}
\label{sec:vards}

The variation of the D-SPT interaction potential, $\phi=\phi(d,\cos\alpha)$, is given by
\begin{equation}
	\label{eq:var1a1a}
	\begin{aligned}
		\delta \phi &= \phi_{,d} \; \delta d + \phi_{,\cos \alpha} \; \delta \cos\alpha.
	\end{aligned}
\end{equation}
%
%
The first term gives
\begin{equation}
	\label{eq:var1aa}
	\begin{aligned}
		\phi_{,d} \; \delta d = \phi_{,d} \; \veh{n} \cdot (\delta \ivh{u}{}{\text{B}}  - \delta \ivh{ u}{}{\text{S}} ) = 
		\iv{f}{}{\text{B}} \cdot \delta \iv{u}{}{\text{B}} + \iv{f}{}{\text{S}} \cdot \delta \ivh{ u}{}{\text{S}} = \ve{f} \cdot \delta \iv{u}{}{\text{B}} - \ve{f} \cdot \delta \ivh{ u}{}{\text{S}},
	\end{aligned}
\end{equation}
where $\ve{f}:=\phi_{,d} \;\veh{n}$ is the D-SPT interaction force. Actually, this quantity has dimensions of [force/length], but we will refer to it as the \emph{interaction force} in the following.
The second term gives
\begin{equation}
	\label{eq:var1a1a11}
	\begin{aligned}
		\phi_{,\cos \alpha} \; \delta \cos\alpha = -\phi_{,\cos \alpha} \; \frac{t_n}{\cos\alpha}  (\veh{n} \cdot \delta \ve{t} + \ve{t} \cdot \delta \veh{n}) = w_c (\veh{n} \cdot \delta \ve{t} + \ve{t} \cdot \delta \veh{n}),
	\end{aligned}
\end{equation}
where 
\begin{equation}
	\label{eq:var1a1a11x}
	\begin{aligned}
		w_c:=  -\phi_{,\cos \alpha} \; \frac{t_n}{\cos\alpha} =  -\phi^R_{,\cos \alpha} t_n \quad \text{and} \quad \phi^R_{,\cos \alpha}:=\frac{\phi_{,\cos \alpha}}{\cos \alpha}. 
	\end{aligned}
\end{equation}
We introduce a reduced potential $\phi^R$ to avoid a singularity issue when $\cos\alpha\rightarrow0$ \cite{2025borković}. The variation of the tangent is given in \eqqref{eq:inctn}, while the variation (linear increment) of the normal is given in Appendix \ref{appendixd}. With these expressions at hand, we get
\begin{equation}
	\label{eq:var1a1a110}
	\begin{aligned}
		\phi_{,\cos \alpha} \; \delta \cos\alpha &= \frac{w_c}{\sqrt{j_\mathrm{B}}} (\veh{n} - t_n \ve{t} )\cdot \delta \ve{u}_{\mathrm{B},1}\\
		& +\frac{w_c}{d} (\ve{t} - (\ve{t} \cdot \ivh{g}{}{\alpha}) \bar{\hat{\ve{g}}}^\alpha - t_n \veh{n}) (\delta \iv{u}{}{\text{B}} -\delta \ivh{u}{}{\text{S}}) - w_c (\ve{t} \cdot \bar{\hat{\ve{g}}}^\beta ) \veh{n} \cdot \delta \ivh{u}{}{\text{S},\beta}  \\
		&= \ve{f}_{\mathrm{CB}}\cdot \delta \ve{u}_{\mathrm{B},1}
	+\ve{f}_\mathrm{C} \cdot (\delta \iv{u}{}{\text{B}} -\delta \iv{\hat u}{}{\text{S}}) -\ve{f}_{\mathrm{CS}}^\beta \cdot \delta \iv{\hat u}{}{\text{S},\beta},  \\
		\ve{f}_{\mathrm{CB}} &:=\frac{w_c}{\sqrt{j}} (\veh{n} - t_n \ve{t} ) = w_c \tilde{\ve{f}}_{\mathrm{CB}},\\
		\ve{f}^\beta_{\mathrm{CS}} &:=w_c (\ve{t} \cdot \bar{\hat{\ve{g}}}^\beta ) \veh{n} = w_c \tilde{\ve{f}}^\beta_{\mathrm{CS}},\\
		\ve{f}_\mathrm{C} &:=\frac{w_c}{d} (\ve{t} - (\ve{t} \cdot \ivh{g}{}{\alpha}) \bar{\hat{\ve{g}}}^\alpha - t_n \veh{ n}) = w_c \tilde{\ve{f}}_\mathrm{C}.
	\end{aligned}
\end{equation}
Due to the variation of potential w.r.t.~the angle, we obtain a contribution to the interaction force that is work-conjugate to the virtual displacements, $\ve{f}_\mathrm{C}$, and interaction moments that are work-conjugate to the variations of the tangent vectors of the beam and the shell, $\ve{f}_{\mathrm{CB}}$ and $\ve{f}_{\mathrm{CS}}$. For the implementation, it is convenient to express
\begin{equation}
	\label{eq:var1aff1a110}
	\begin{aligned}
		\delta t_n & = \delta (\ve{t}\cdot\veh{n}) = \tilde{\ve{f}}_{\mathrm{CB}}\cdot \delta \ve{u}_{\mathrm{B},1}
		+\tilde{\ve{f}}_\mathrm{C} \cdot (\delta \iv{u}{}{\text{B}} -\delta \ivh{u}{}{\text{S}}) -\tilde{\ve{f}}_\mathrm{CS} \cdot \delta \ivh{ u}{}{\text{S},\beta} .
	\end{aligned}
\end{equation}

The variation of the D-SPT potential should be integrated along the beam length to find the complete variation of the beam-shell potential, i.e.,
\begin{equation}
	\label{eq:varBS}
	\begin{aligned}
		\delta \Pi_{\mathrm{B-S}} &= \int_{L_\mathrm{B}}\delta \phi \dd{s_\mathrm{B}}= \int_{L_\mathrm{B}}(\ve{f}+\ve{f}_\mathrm{C} ) \cdot \delta \iv{u}{}{\text{B}} \dd{s_\mathrm{B}}-\int_{L_\mathrm{B}}(\ve{f}+\ve{f}_\mathrm{C} ) \dd{s_\mathrm{B}} \cdot \delta \ivh{u}{}{\text{S}} \\
		&+\int_{L_\mathrm{B}}\ve{f}_{\mathrm{CB}}\cdot \delta \ve{u}_{\mathrm{B},1} \dd{s_\mathrm{B}}-\int_{L_\mathrm{B}}\ve{f}_{\mathrm{CS}} \dd{s_\mathrm{B}} \cdot  \delta \ivh{ u}{}{\text{S},\beta} 
	\end{aligned}
\end{equation}
%
%
In the following, we will consider three formulations: (i) full formulation (FF) that employs complete \eqqref{eq:varBS}, (ii) reduced formulation (RF1) given by \eqqref{eq:var1aa} that neglects contributions from the variation of angle, and (iii) reduced formulation (RF2) that neglects both the variation w.r.t.~angle and the orientation of the cross section by adopting $\phi=\phi(d,1)$.

Details on the linearization of SBS contribution to the total potential energy are presented in Appendices \ref{appendixd} and \ref{appendixe}, while the variation of orthogonality condition is given in Appendix \ref{appendixc}.

\section{Numerical examples}
\label{secnum}

To assess the applicability of the derived formulation, we consider two examples that model fundamental cases of interaction between fibers (modeled as beams) and membranes (modeled as shells): (i) peeling of an adhering beam from a shell, and (ii) bending of a shell by an adhering beam. In both examples, the beam and the shell interact via the LJ potential.

The beam and the shell are discretized using the isogeometric approach, employing quartic B-splines with $C^1$ interelement continuity. All contributions to the total potential are integrated with $p+1=5$ Gauss integration points. The level of spatial discretization is marked as $n_\mathrm{el}=n_\mathrm{el,B}\times n_\mathrm{el,S1}\times n_\mathrm{el,S2}$, where $n_\mathrm{B}$ is the number of elements along the beam, while $n_\mathrm{el,S\alpha}$ $(\alpha=1,2)$ are the numbers of elements along the shell's $\theta^\alpha$ directions. No cutoff is applied, and the contributions of all points are fully taken into account.

%
%

\subsection{Peeling of an adhering beam from a shell}
\label{sec:ex1}

\begin{figure}[h!]
	\centering
	\includegraphics[width=0.9\textwidth]{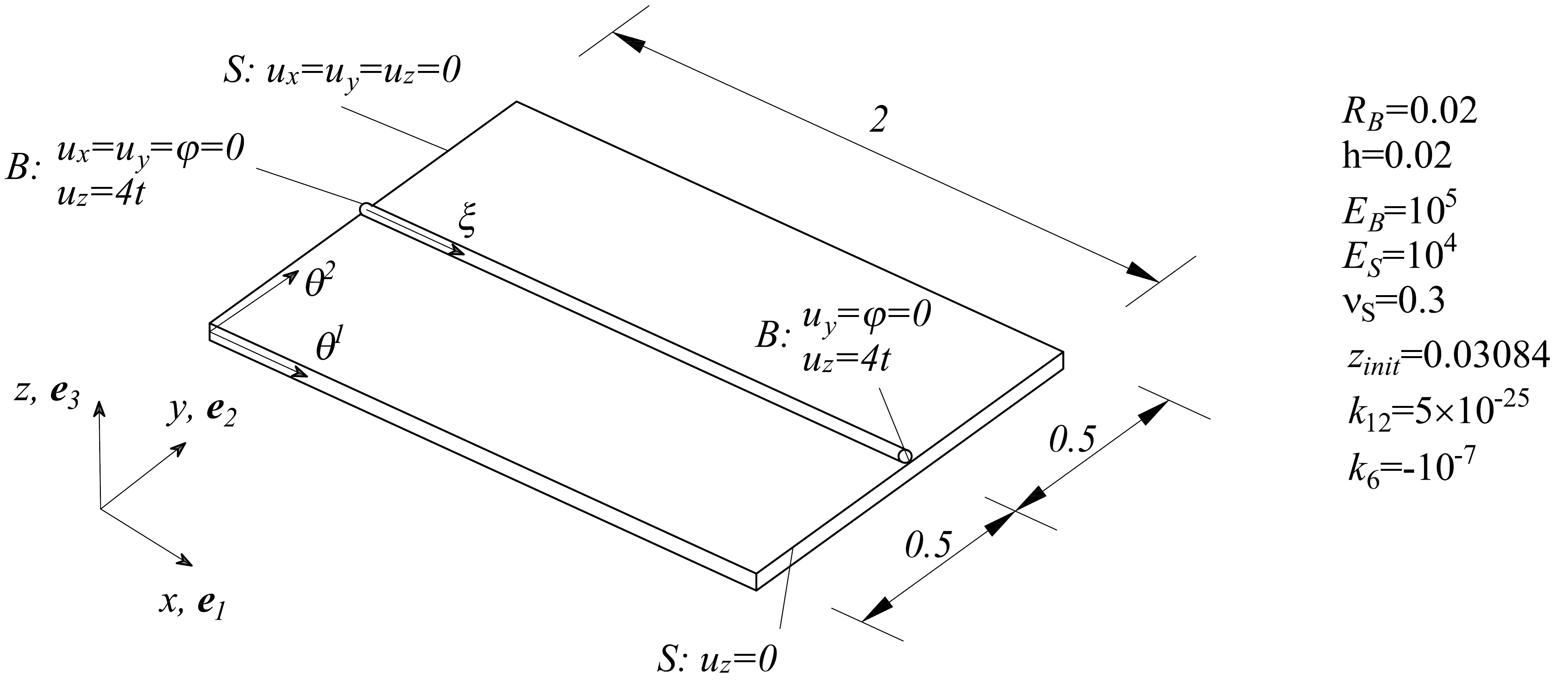}
	\caption{Peeling of an adhering beam from a shell. Problem setup.}
	\label{fig:num1setup}
\end{figure}

\begin{figure}[h!]
	\centering
	\includegraphics[width=0.95\textwidth]{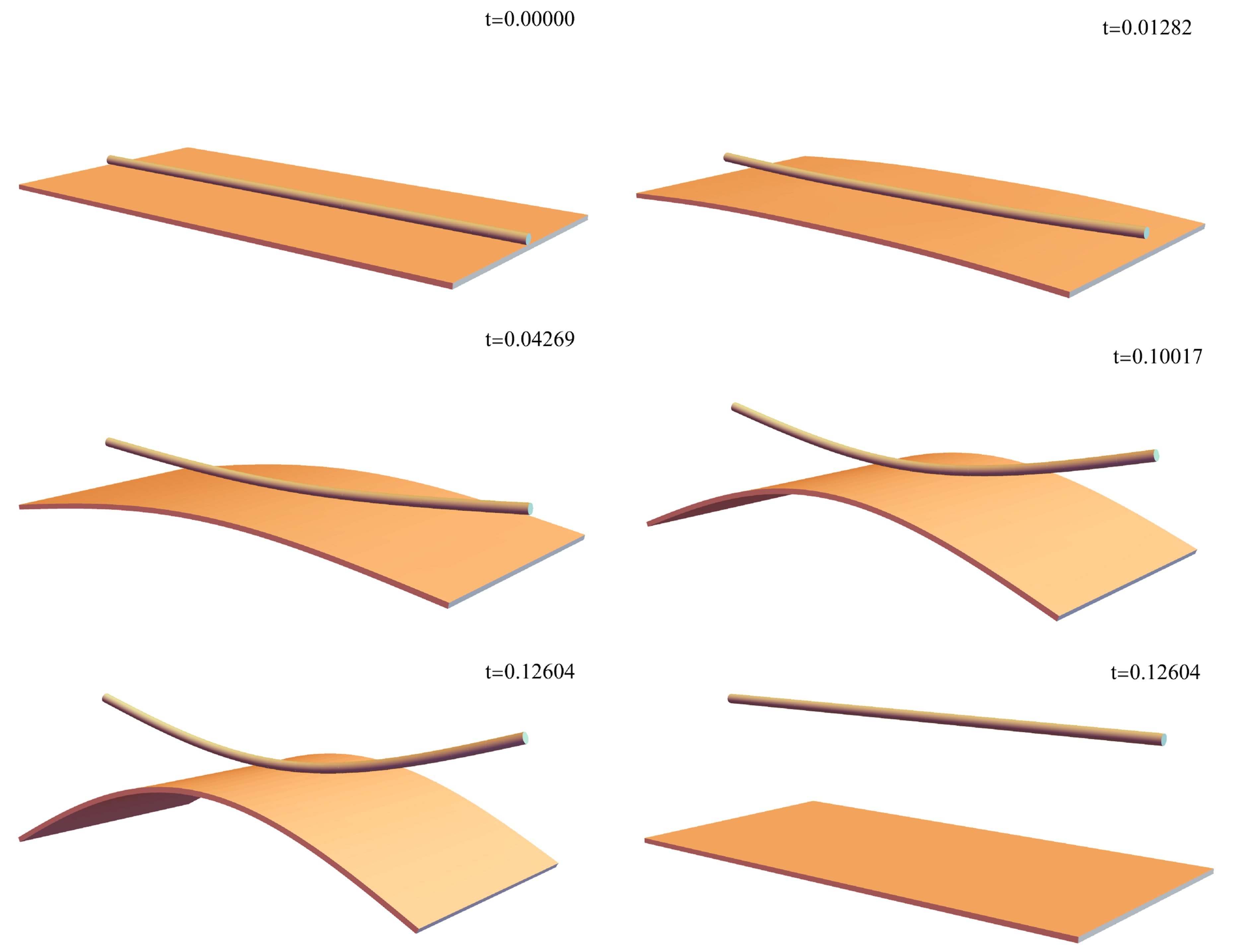}
	\caption{Peeling of an adhering beam from a shell. Six characteristic configurations.}
	\label{fig:num1conf}
\end{figure}

We consider the peeling of a beam from a shell. The setup and corresponding parameters are displayed in Fig.~\ref{fig:num1setup}.
At $t=0$, the beam and the shell are close to their equilibrium separation $z_\mathrm{init}$. We prescribe displacements at both beam ends along the $z-$axis, $u_z=4t$. The obtained configurations are shown in Fig.~\ref{fig:num1conf}.
The behavior is similar to that observed for beams \cite{2024borkovićb}. The process begins with an initial separation of the beam ends, leading to an almost complete peeling of the beam. Only a small region near the middle of the beam remains adhered to the shell, which induces pulling; consequently, both the beam and the shell undergo visible deformation, and the process ends with a pull-off.

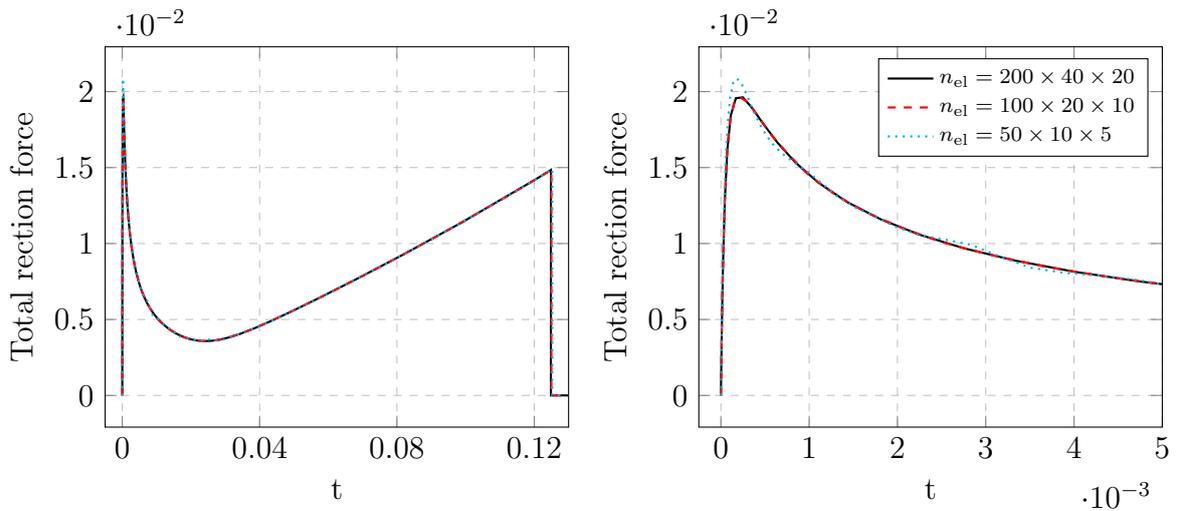
\begin{figure}[h!]
	\centering
	\subfloat[The complete response.]{\label{figj:ha}
		\begin{tikzpicture}
			\begin{axis}[
				xlabel = {t},
				ylabel = {Total rection force},
				legend pos= north east,
				legend cell align=left,
				width=0.48\textwidth,
				xmin=-0.005,xmax=0.13,
				xtick distance=0.04,	
				clip=false,
				grid=both,
				xticklabel style={/pgf/number format/fixed, /pgf/number format/precision=3}
				]			
				\addplot[black,thick,restrict x to domain=0:0.13] table [col sep=comma] {data/sPeelReactSUMM200x40x20IPT881.csv};	
				\addplot[red,thick,dashed,restrict x to domain=0:0.13] table [col sep=comma] {data/sPeelReactSUMM100x20x10IPT881.csv};
				\addplot[cyan,thick,dotted,restrict x to domain=0:0.13]table [col sep=comma]{data/sPeelReactSUMM50x10x5IPT881.csv};
			\end{axis}
			
		\end{tikzpicture}
	}
	\subfloat[The response for $t \in (0,0.005)$.]{\label{fijgh:b}
		\begin{tikzpicture}
			\begin{axis}[
				xlabel = {t},
				ylabel = {Total rection force
				},
				legend pos= north east,
				legend style={font=\scriptsize},
				legend cell align=left,
				width=0.48\textwidth,
				xmin=-0.00025,xmax=0.005,
				clip=true,
				grid=both]
				
				\addplot[black,thick,restrict x to domain=0:0.01] table [col sep=comma] {data/sPeelReactSUMM200x40x20IPT881.csv};
				
				\addplot[red,thick,dashed,restrict x to domain=0:0.01] table [col sep=comma] {data/sPeelReactSUMM100x20x10IPT881.csv};
				\addplot[cyan,thick,dotted,restrict x to domain=0:0.01]table [col sep=comma]{data/sPeelReactSUMM50x10x5IPT881.csv};
				
				\legend{$n_\mathrm{el}=200\times40\times20$,$n_\mathrm{el}=100\times20\times10$,$n_\mathrm{el}=50\times10\times5$}   
			\end{axis}
			
		\end{tikzpicture}
	}
	\caption{Peeling of an adhering beam from a shell. Sum of the reaction forces on a beam for different FE meshes.}
	\label{fig:ex1conv}
\end{figure}

As a primary variable, we observe the total reaction force on the beam -- the sum of reaction forces on beam ends. In contrast to the beam-beam simulations in \cite{2024borkovićb}, the reaction forces on beam ends in this case are not equal. The total reaction force for three FE meshes is shown in Fig.~\ref{fig:ex1conv}.
From these results, we adopt the mesh $n_\mathrm{el}=100\times20\times10$.

\begin{figure}[h!]
	\centering
	\includegraphics[width=1\textwidth]{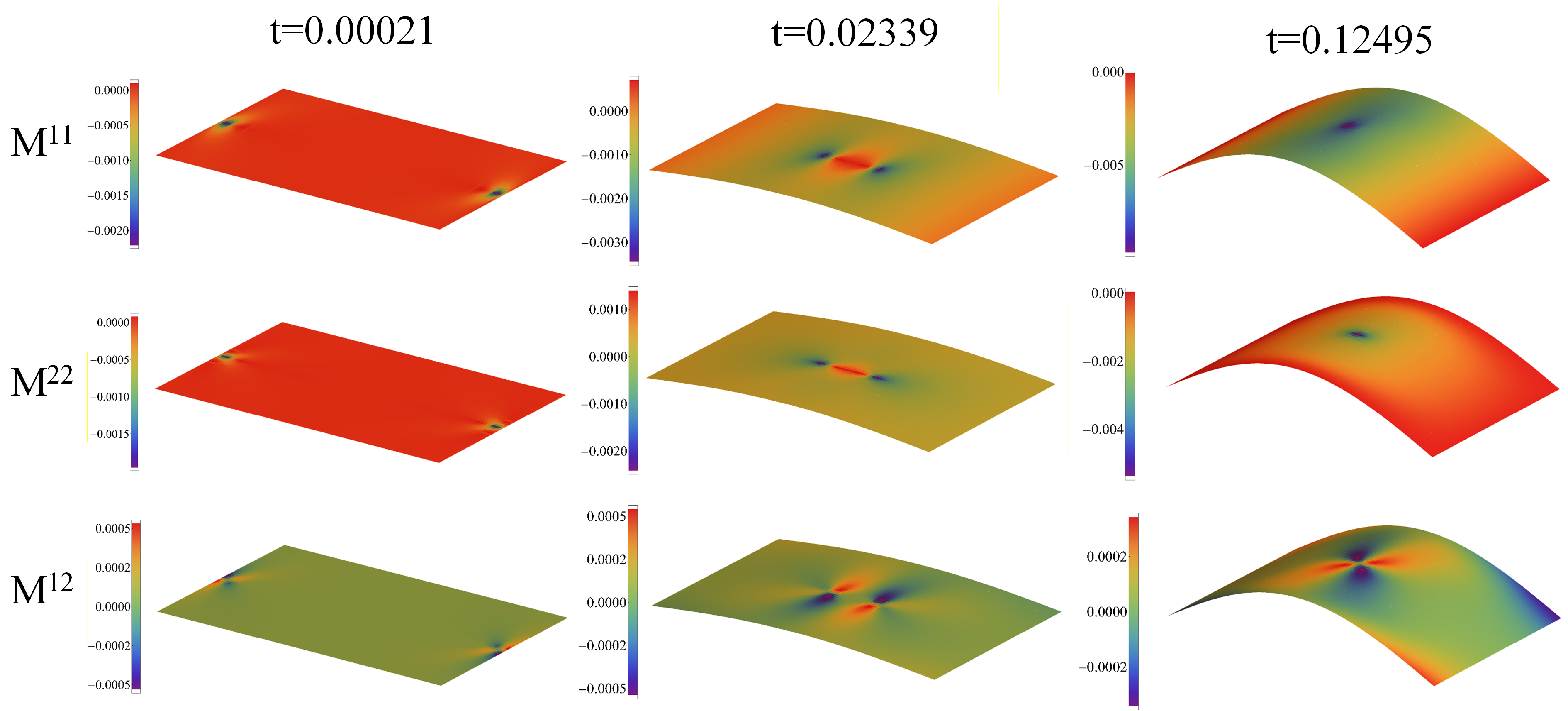}
	\caption{Peeling of an adhering beam from a shell. Shell's section moments.}
	\label{fig:ex1SF}
\end{figure}

By observing the total reaction force on the beam, three characteristic time instances w.r.t.~the reaction force can be defined: (i) maximum reaction force at the initiation of peeling when $t\approx 0.00021$, (ii) minimum reaction force at $t\approx 0.02339$, (iii) reaction force at pull-off when $t\approx 0.12495$. The shell's section moments at these characteristic quasi-time instances are plotted in Fig.~\ref{fig:ex1SF}.
%
Stress concentrations occur at locations of the peeling front.

\begin{figure}[h!]
	\centering
	\begin{tikzpicture}
			\begin{axis}[
					xlabel = {$\xi$},
					ylabel = {Interaction force magnitude
						},
					legend style={at={(0.8,0.15)},anchor=west},
					legend style={font=\scriptsize},
					legend cell align=left,
					height=0.5\textwidth,
					width=0.95\textwidth,
					xmin=0.0,xmax=1,
					clip=false,
					grid=both]
					\addplot[black,thick] table [col sep=comma] {data/ipFPeelEB105ES104IP1IP881M100x20P1.csv};
					\addplot[red,thick] table [col sep=comma] {data/ipFPeelEB105ES104IP1IP881M100x20P2.csv};
					\addplot[green,thick] table [col sep=comma] {data/ipFPeelEB105ES104IP1IP881M100x20P23.csv};
					\addplot[blue,thick] table [col sep=comma] {data/ipFPeelEB105ES104IP1IP881M100x20P3.csv};

					\legend{$t=0.00021$,$t=0.00304$,$t=0.02339$,$t=0.12495$}   
				\end{axis}
		\end{tikzpicture}
	\caption{Peeling of an adhering beam from a shell. Distribution of the magnitude of the interaction force along the beam at four quasi-time instances.}
	\label{fig:ex1IPdistr}
\end{figure}
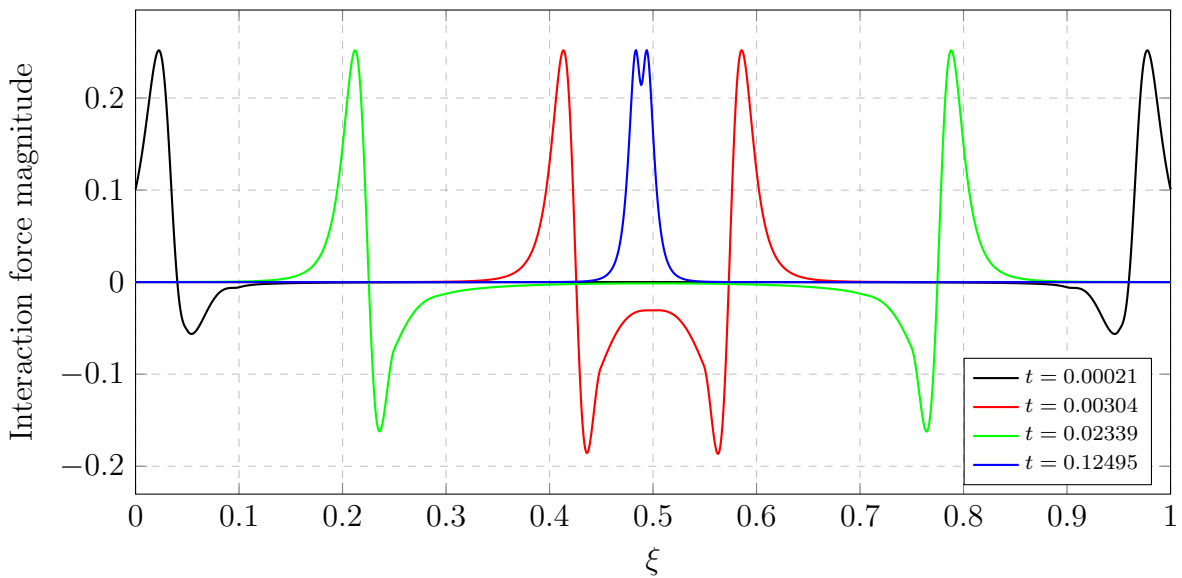

To gain a better insight into this problem, let us examine the interaction force. Distributions of magnitudes of the interaction force along the beam are plotted in Fig.~\ref{fig:ex1IPdistr} for four quasi-time instances.
Distributions of the interaction force are similar to those observed for the beam-beam peeling example \cite{2024borkovićb}. The most obvious difference is a small loss of symmetry, since the distribution of the interaction force is not symmetric w.r.t.~$\xi=0.5$. This loss of symmetry is consistent with the observed difference in the reaction forces at the beam ends.

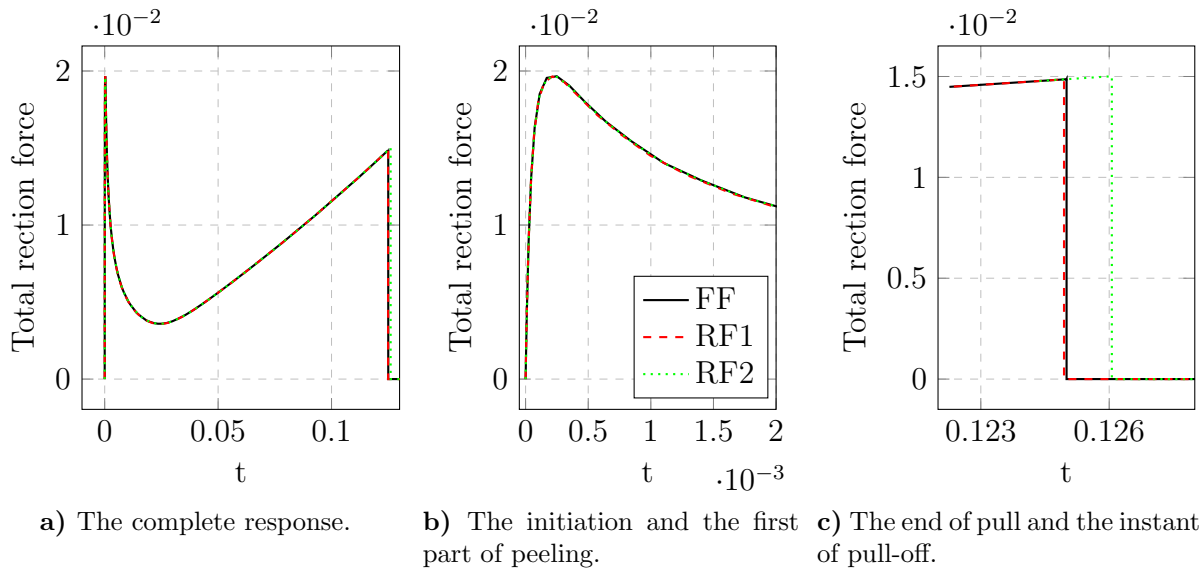
\begin{figure}[h!]
	\centering
	\subfloat[The complete response.]{\label{fig:ha}
		\begin{tikzpicture}
			\begin{axis}[
				xlabel = {t},
				ylabel = {Total rection force
				},
				legend cell align=left,
				legend pos = south east,
				width=0.36\textwidth,
				height=0.4\textwidth,
				xmin=-0.01,xmax=0.13,
				clip=true,
				grid=both,
				xticklabel style={/pgf/number format/fixed, /pgf/number format/precision=3}]
				\addplot[black,thick] table [col sep=comma] {data/sPeelReactSUMM100x20x10IPT880.csv};
				\addplot[red,dashed,thick] table [col sep=comma] {data/sPeelReactSUMM100x20x10IPT881.csv};
				\addplot[green,thick,dotted] table [col sep=comma] {data/sPeelReactSUMM100x20x10IPT882.csv};

			\end{axis}
		
		\end{tikzpicture}
	} \:
	\subfloat[The initiation and the first part of peeling.]{\label{figh:b}
		\begin{tikzpicture}
			\begin{axis}[
				xlabel = {t},
				ylabel = {Total rection force
				},
				legend pos= south east,
				legend cell align=left,
				width=0.31\textwidth,
					height=0.4\textwidth,
				xmin=-0.00005,xmax=0.002,
				clip=true,
				grid=both,
				xticklabel style={/pgf/number format/fixed, /pgf/number format/precision=3}]
				\addplot[black,thick,restrict x to domain=0:0.01] table [col sep=comma] {data/sPeelReactSUMM100x20x10IPT880.csv};
				\addplot[red,dashed,thick,restrict x to domain=0:0.01] table [col sep=comma] {data/sPeelReactSUMM100x20x10IPT881.csv};
				\addplot[green,thick,dotted,restrict x to domain=0:0.01] table [col sep=comma] {data/sPeelReactSUMM100x20x10IPT882.csv};
				\legend{FF,RF1,RF2} 
			
			\end{axis}

		\end{tikzpicture}
	}\:
	\subfloat[The end of pull and the instant of pull-off.]{\label{figh:fb}
		\begin{tikzpicture}
			\begin{axis}[
				xlabel = {t},
				ylabel = {Total rection force
				},
				legend pos= north east,
				legend cell align=left,
				width=0.31\textwidth,
				height=0.4\textwidth,
				xmin=0.122,xmax=0.128,
				xtick distance=0.003,	
				clip=true,
				grid=both,
				xticklabel style={/pgf/number format/fixed, /pgf/number format/precision=3}
				]
				\addplot[black,thick,restrict x to domain=0.12:0.13] table [col sep=comma] {data/sPeelReactSUMM100x20x10IPT880.csv};
				\addplot[red,dashed,thick,restrict x to domain=0.12:0.13] table [col sep=comma] {data/sPeelReactSUMM100x20x10IPT881.csv};
				\addplot[green,thick,dotted,restrict x to domain=0.12:0.13] table [col sep=comma] {data/sPeelReactSUMM100x20x10IPT882.csv};
				  
			\end{axis}
			
		\end{tikzpicture}
	}
	\caption{Peeling of an adhering beam from a shell. Reaction force on the beam vs.~quasi-time. Comparison of three formulations.}
	\label{fig:ex1form}
\end{figure}

We then compare the total reaction force obtained using the three formulations defined in Subsection \ref{sec:vards}, see Fig.~\ref{fig:ex1form}.
The results indicate that all three formulations yield very similar responses, although small differences can be observed. First, at the initiation of peeling, when the whole beam is still adhering to the shell, all formulations give practically the same result. This implies that, for cases in which the beam is expected to remain adhered to the shell throughout the simulation, the simplest formulation, RF2, can be employed. Second, RF2 overestimates the instant and the force at pull-off due to neglecting the orientation of the beam cross sections w.r.t.~the surrogate plate. This effect becomes noticeable at the instant of pull-off since both the beam and the shell are then maximally curved. Third, although there is a small discrepancy between RF1 and FF, it is negligible in the considered setting. 

\begin{figure}[h!]
	\centering
	\subfloat[Complete responses.]{\label{figj:ha}
		\begin{tikzpicture}
					\begin{axis}[
				xlabel = {t},
				ylabel = {Total rection force
				},
				legend pos= north east,
				legend style={font=\scriptsize},
				legend cell align=left,
				width=0.6\textwidth,
				height=0.4\textwidth,
				xmin=-0.005,xmax=0.25,
				xticklabel style={/pgf/number format/fixed, /pgf/number format/precision=3},
				clip=true,
				grid=both]		
				
				\addplot[blue,thick] table [col sep=comma] {data/sPeelReactSUMMEB103ES104100x20x10IPT881.csv};		
				\addplot[green,thick] table [col sep=comma] {data/sPeelReactSUMMEB104ES104100x20x10IPT881.csv};	
				\addplot[black,thick] table [col sep=comma] {data/sPeelReactSUMM100x20x10IPT881.csv};
				\addplot[red,thick] table [col sep=comma] {data/sPeelReactSUMMEB106ES104100x20x10IPT881.csv};
				\legend{ $E_B=10^3$ $E_S=10^4$, $E_B=10^4$ $E_S=10^4$,$E_B=10^5$ $E_S=10^4$ , $E_B=10^6$ $E_S=10^4$
				}   
			\end{axis}
			
		\end{tikzpicture}
	}
	\subfloat[Responses for $t \in (0,0.01)$.]{\label{fijgh:b}
		\begin{tikzpicture}
		\begin{axis}[
	xlabel = {t},
	ylabel = {Total rection force
	},
	legend pos= north east,
	legend style={font=\scriptsize},
	legend cell align=left,
	width=0.35\textwidth,
	height=0.4\textwidth,
	xmin=-0.0005,xmax=0.01,
	clip=true,
	grid=both]		
	
		\addplot[blue,thick] table [col sep=comma] {data/sPeelReactSUMMEB103ES104100x20x10IPT881.csv};		
\addplot[green,thick] table [col sep=comma] {data/sPeelReactSUMMEB104ES104100x20x10IPT881.csv};	
\addplot[black,thick] table [col sep=comma] {data/sPeelReactSUMM100x20x10IPT881.csv};
\addplot[red,thick] table [col sep=comma] {data/sPeelReactSUMMEB106ES104100x20x10IPT881.csv};
 
\end{axis}
			
		\end{tikzpicture}
	}
	\caption{Peeling of an adhering beam from a shell. Reaction force on beam vs.~quasi-time for different Young's moduli.}
	\label{fig:ex1modul}
\end{figure}
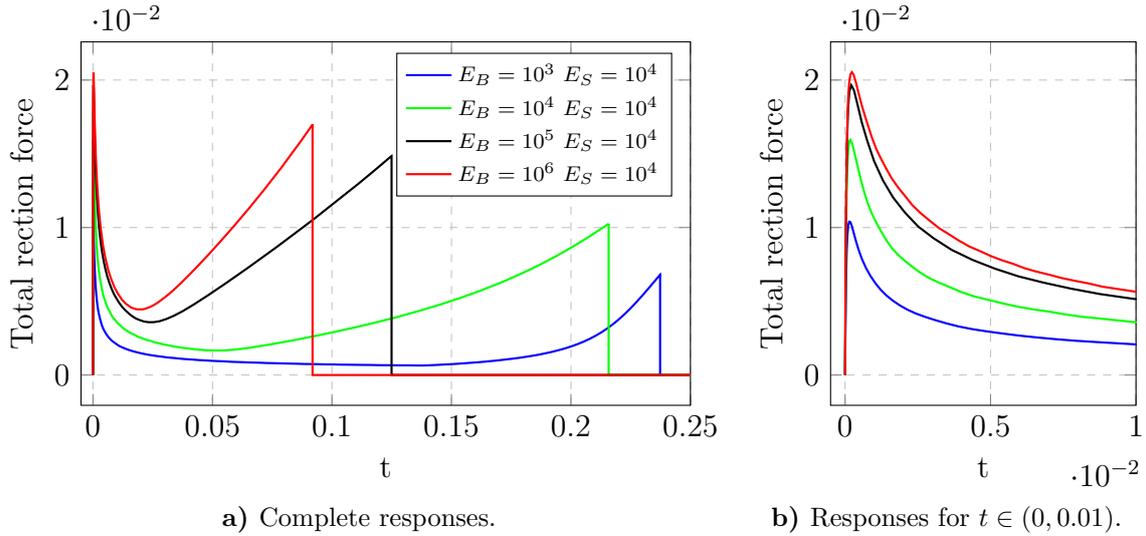

Next, we investigate the influence of Young's modulus on the reaction force using the RF1 formulation, see Fig. \ref{fig:ex1modul}.
The ratio $E_\mathrm{B}/E_\mathrm{S}$ significantly affects the structural response. Both the peak force at the initiation of peeling and the pull-off force decrease with the increase in $E_\mathrm{B}/E_\mathrm{S}$ ratio. On the other hand, the pull-off instant and the deformation of both the beam and the shell increase with $E_\mathrm{B}/E_\mathrm{S}$.

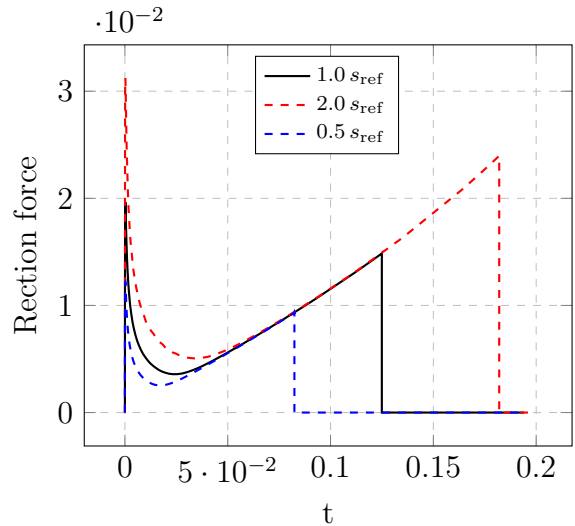
\begin{figure}[h!]
	\centering
	\begin{tikzpicture}
		\begin{axis}[
			xlabel = {t},
			ylabel = {Rection force
			},
			legend style={at={(0.35,0.85)},anchor=west},
			legend style={font=\scriptsize},
			legend cell align=left,
			width=0.5\textwidth,
			xtick distance=0.05,	
			clip=false,
			grid=both]
			\addplot[black,thick,restrict x to domain=0:0.2] table [col sep=comma] {data/sPeelReactSUMM100x20x10IPT881.csv};
			\addplot[red,thick,dashed,restrict x to domain=0:0.2] table [col sep=comma] {data/sPeelReactSUMMEBIP2x100x20x10IPT881.csv};
			\addplot[blue,dashed,thick,restrict x to domain=0:0.2] table [col sep=comma] {data/sPeelReactSUMMEBIP05x100x20x10IPT881.csv};

			\legend{$1.0 \, s_\mathrm{ref}$,$2.0 \,  s_\mathrm{ref}$,$0.5 \, s_\mathrm{ref}$}   
		\end{axis}
	\end{tikzpicture}
	\caption{Peeling of an adhering beam from a shell. Reaction force vs.~quasi-time for three different interaction strengths.}
	\label{fig:num1INT}
\end{figure}

Finally, let us vary the interaction strength, $s_i$, by multiplying the LJ potential defined in Fig.~\ref{fig:num1setup} by factors 0.5 and 2. The results are shown in Fig.~\ref{fig:num1INT}.
All characteristic quantities are proportional to the interaction strength: peak reaction forces and the instant of pull-off. The proportionality of the peak reaction forces, $f_p$, to the interaction strength is a power function with an approximate exponent $0.678$, i.e., $f_p \propto s_i^{0.678}$. 

\subsection{Shell bending by an adhering fiber}

\begin{figure}[h!]
	\centering
	\includegraphics[width=0.95\textwidth]{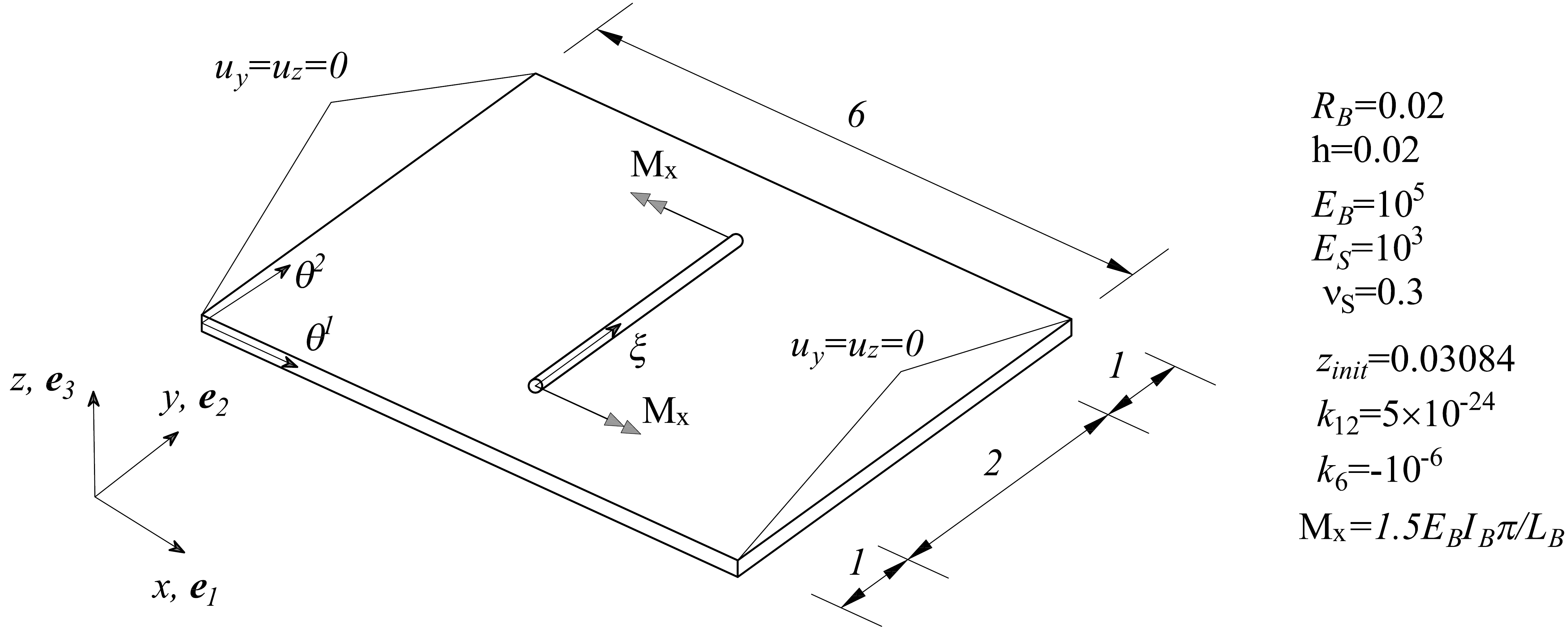} 
	\caption{Shell bending by an adhering fiber. Problem setup.}
	\label{fig:ex2setup}
\end{figure}

We observe a beam/fiber that adheres to a shell. There are no kinematic boundary conditions on the beam, while the shell's $Y$ and $Z$ displacement components at the corners are fixed. The beam is loaded with two point moments, $M_x$, applied in opposite directions at its ends, Fig.~\ref{fig:ex2setup}.
As the moments increase, the beam bends and, due to adhesion, the shell deforms correspondingly. By exploiting symmetry of the setup, we discretize only a quarter of the shell and one half of the beam.

\begin{figure}[h!]
	\centering
	\includegraphics[width=0.95\textwidth]{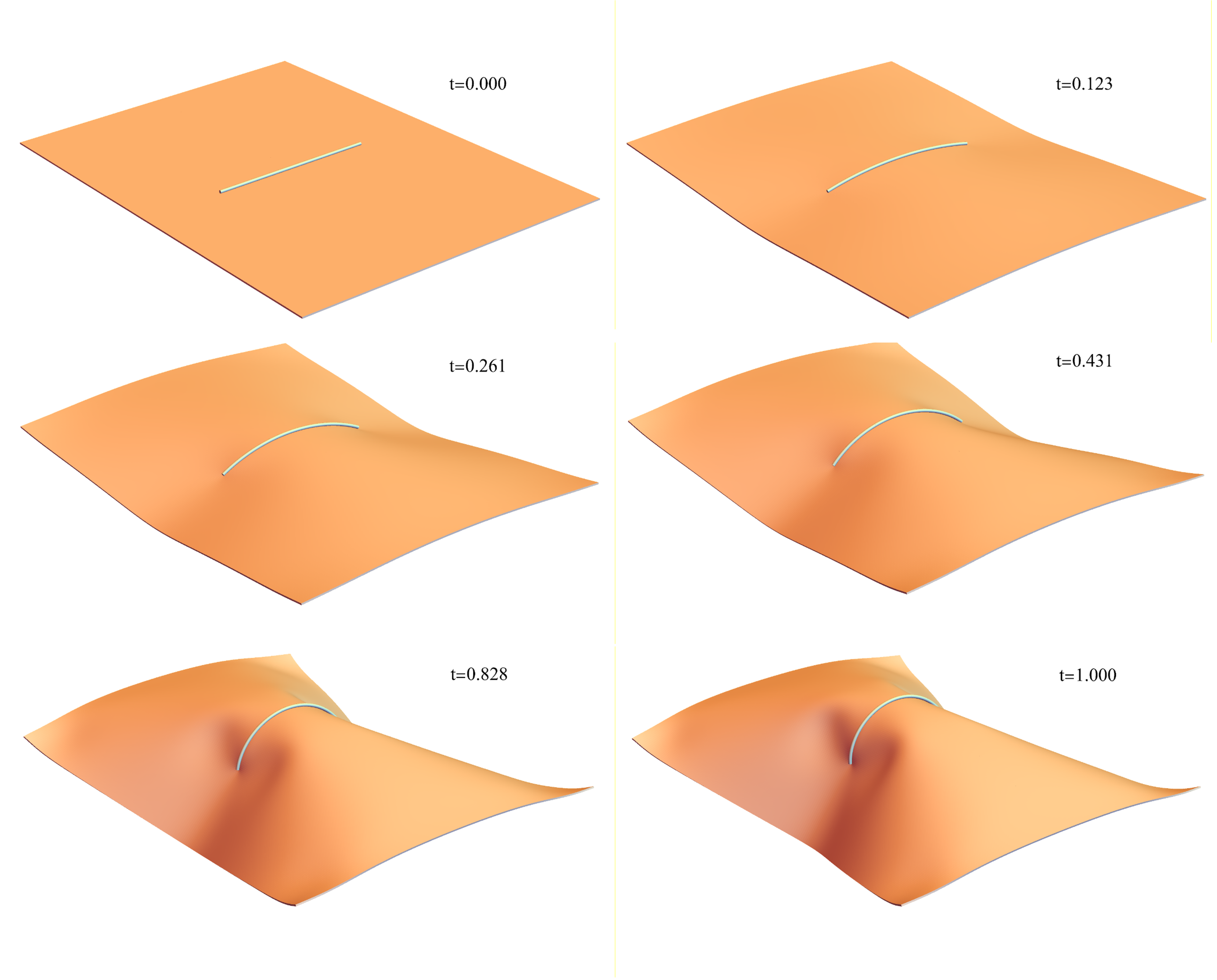} 
	\caption{Shell bending by an adhering fiber. Six characteristic configurations.}
	\label{fig:ex2conf}
\end{figure}

Six characteristic configurations are presented in Fig.~\ref{fig:ex2conf} where the observed behavior is consistent with expectations: The beam adheres to the shell and bends it.
Although the applied point moments would bend a free beam into a three-quarter circle, the shell resists this bending, resulting in a final beam configuration closer to a half-circle. Based on the results obtained in Subsection \ref{sec:ex1} and the fact that the beam remains adhered to the shell throughout the simulation, we employ the RF2 formulation exclusively.

\begin{figure}[h!]
	\centering
	\subfloat[Displacement component $u_z$ at the beam endpoint.]{\label{fig:a}
		\begin{tikzpicture}
		\begin{axis}[
			xlabel = {t},
			ylabel = {Displacement},
			legend pos= north west,
			legend cell align=left,
			legend style={font=\scriptsize},
			width=0.5\textwidth,
			clip=false,
			xmin=-0.05,xmax=1.05,
			grid=both,
			xticklabel style={/pgf/number format/fixed, /pgf/number format/precision=2},
			yticklabel style={/pgf/number format/fixed, /pgf/number format/precision=2},]

				\addplot[red,only marks,line width=0.25pt,solid,mark=star,mark size=1.25pt] table [col sep=comma] {data/dispZBksi0SqueezeNewEB105ES103IP10IP882M100x50.csv};
			
			\addplot[blue,only marks,line width=0.25pt,solid,mark=square,mark size=1.25pt] table [col sep=comma] {data/dispZBksi0SqueezeNewEB105ES103IP10IP882M200x70.csv};
			
			\addplot[green,only marks,solid,mark=triangle,mark size=1.25pt] table [col sep=comma] {data/dispZBksi0SqueezeNewEB105ES103IP10IP882M400x100.csv};
			
			\addplot[black,dashed,line width=0.25pt,solid,mark=o,mark size=1.25pt] table [col sep=comma] {data/dispZBksi0SqueezeNewEB105ES103IP10IP882M800x140.csv};
	\legend{
	$n_{el}=100\times50\times50$,
	$n_{el}=200\times70\times70$,
	$n_{el}=400\times100\times100$,
	$n_{el}=800\times140\times140$}
			
		\end{axis}
		\end{tikzpicture}
	}
	\subfloat[Interaction force at the beam endpoint.]{\label{fig:b}
		\begin{tikzpicture}
				\begin{axis}[
			xlabel = {t},
			ylabel = {Interaction force},
			legend pos= north east,
			legend cell align=left,
			width=0.5\textwidth,
			clip=false,
			xmin=-0.05,xmax=1.05,
			grid=both]

%
%
%

			\addplot[red,only marks,line width=0.25pt,solid,mark=star,mark size=1.25pt] table [col sep=comma] {data/ipFksi0SqueezeNewEB105ES103IP10IP882M100x50.csv};
			
			\addplot[blue,only marks,line width=0.25pt,solid,mark=square,mark size=1.25pt] table [col sep=comma] {data/ipFksi0SqueezeNewEB105ES103IP10IP882M200x70.csv};
			
			\addplot[green,only marks,solid,mark=triangle,mark size=1.25pt] table [col sep=comma] {data/ipFksi0SqueezeNewEB105ES103IP10IP882M400x100.csv};
			
			\addplot[black,dashed,line width=0.25pt,solid,mark=o,mark size=1.25pt] table [col sep=comma] {data/ipFksi0SqueezeNewEB105ES103IP10IP882M800x140.csv};
	
		\end{axis}
		\end{tikzpicture}
	}
	\caption{Shell bending by an adhering fiber. Comparison of different FE meshes.}
	\label{fig:convEx2}
\end{figure}
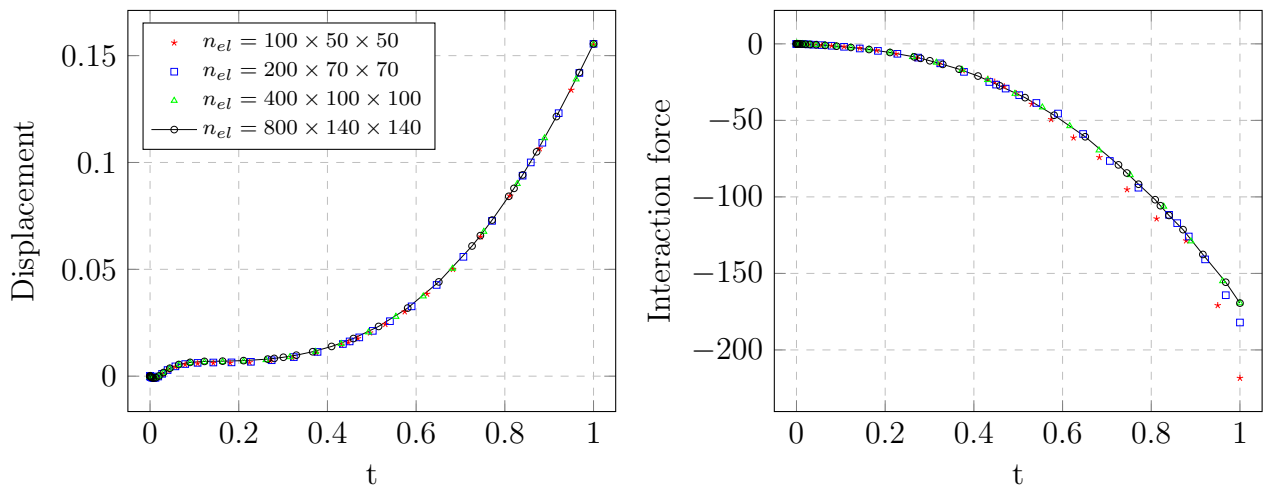

Let us consider some results regarding $h$-refinement. For this, we observe the beam endpoint $\xi=0$, where the displacement component along the $Z$-axis and the interaction force are followed. The results for four selected FE meshes are displayed in Fig.~\ref{fig:convEx2}. For sparse meshes, differences w.r.t.~the densest/reference mesh in displacement primarily occur during the first half of the response, whereas for the interaction force, differences are observed in the latter half. However, these differences are of different orders: they are clearly pronounced for the interaction force but barely noticeable for the displacement component. Based on this analysis, we adopt a mesh with 400 elements along the beam and 100 elements along each edge of the shell, giving  $n_{el}=400\times100\times100$.

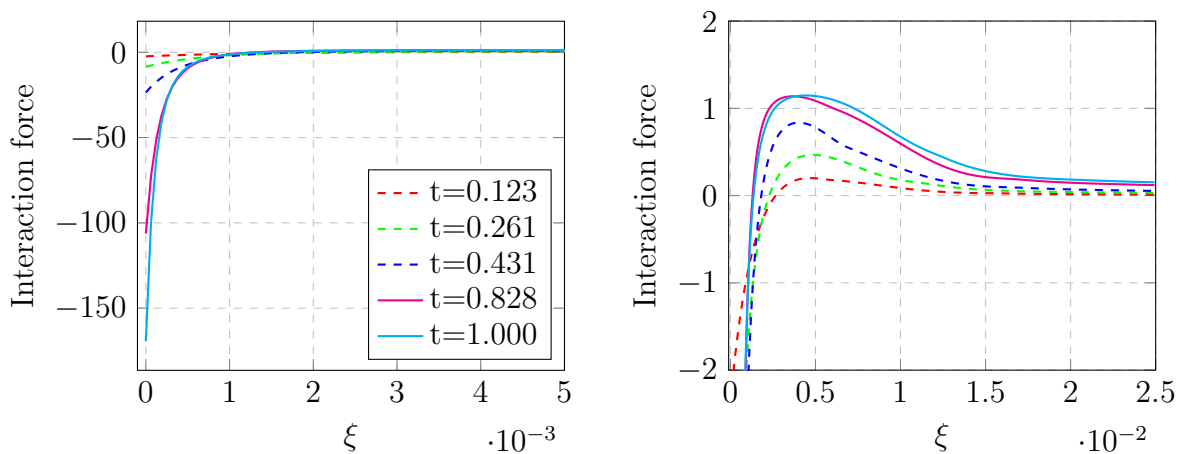
\begin{figure}[h!]
	\centering
	\subfloat[Distribution of the interaction force for $\xi \in ( 0,0.005 ) $.]{\label{fig:a}
		\begin{tikzpicture}
				\begin{axis}[
			xlabel = {$\xi$},
			ylabel = {Interaction force
			},
			legend pos= south east,
			legend cell align=left,
			width=0.45\textwidth,
			xmin=-0.0001,xmax=0.005,
			clip=true,
			grid=both]
			\addplot[red,dashed,thick] table [col sep=comma] {data/ipFDistrT0123.csv};
			\addplot[green,dashed,thick] table [col sep=comma] {data/ipFDistrT0261.csv};
			\addplot[blue,dashed,thick] table [col sep=comma] {data/ipFDistrT0431.csv};
			\addplot[magenta,thick] table [col sep=comma] {data/ipFDistrT0828.csv};
			\addplot[cyan,thick] table [col sep=comma] {data/ipFDistrT1.csv};
			\legend{t=0.123,t=0.261,t=0.431,t=0.828,t=1.000}   
		\end{axis}
		\end{tikzpicture}
	} \;
	\subfloat[Distribution of the interaction force for $\xi \in ( 0,0.025 ) $ with clipped peak values.]{\label{fig:b}
		\begin{tikzpicture}
				\begin{axis}[
			xlabel = {$\xi$},
			ylabel = {Interaction force
			},
			legend pos= south east,
			legend cell align=left,
			width=0.45\textwidth,
			xmin=-0.0001,xmax=0.025,
			ymin = -2, ymax = 2,
			clip=true,
			grid=both]
			\addplot[red,dashed,thick] table [col sep=comma] {data/ipFDistrT0123.csv};
			\addplot[green,dashed,thick] table [col sep=comma] {data/ipFDistrT0261.csv};
			\addplot[blue,dashed,thick] table [col sep=comma] {data/ipFDistrT0431.csv};
			\addplot[magenta,thick] table [col sep=comma] {data/ipFDistrT0828.csv};
			\addplot[cyan,thick] table [col sep=comma] {data/ipFDistrT1.csv};
		\end{axis}
		\end{tikzpicture}
	}
	\caption{Shell bending by an adhering fiber. Distributions of the interaction force on the beam for five quasi-time instances.}
	\label{fig:distIF}
\end{figure}

With the adopted mesh, we consider the distribution of the interaction force at five quasi-time instances. To improve visibility of this highly nonlinear quantity with a steep gradient, its distribution is shown in the close vicinity to the beam endpoint ($\xi=0$) in Fig.~\ref{fig:distIF}a and with clipped peak values in Fig.~\ref{fig:distIF}b.
First, the steep gradient of the interaction force near the beam endpoint is clearly visible, and it increases with quasi-time. Second, the maximum adhesion force is reached already for $t=0.828$ and the resultant adhesion force for $t>0.828$ increases by extending the area with the maximum attraction force value.

\begin{figure}[]
	\centering
	\includegraphics[width=1\textwidth]{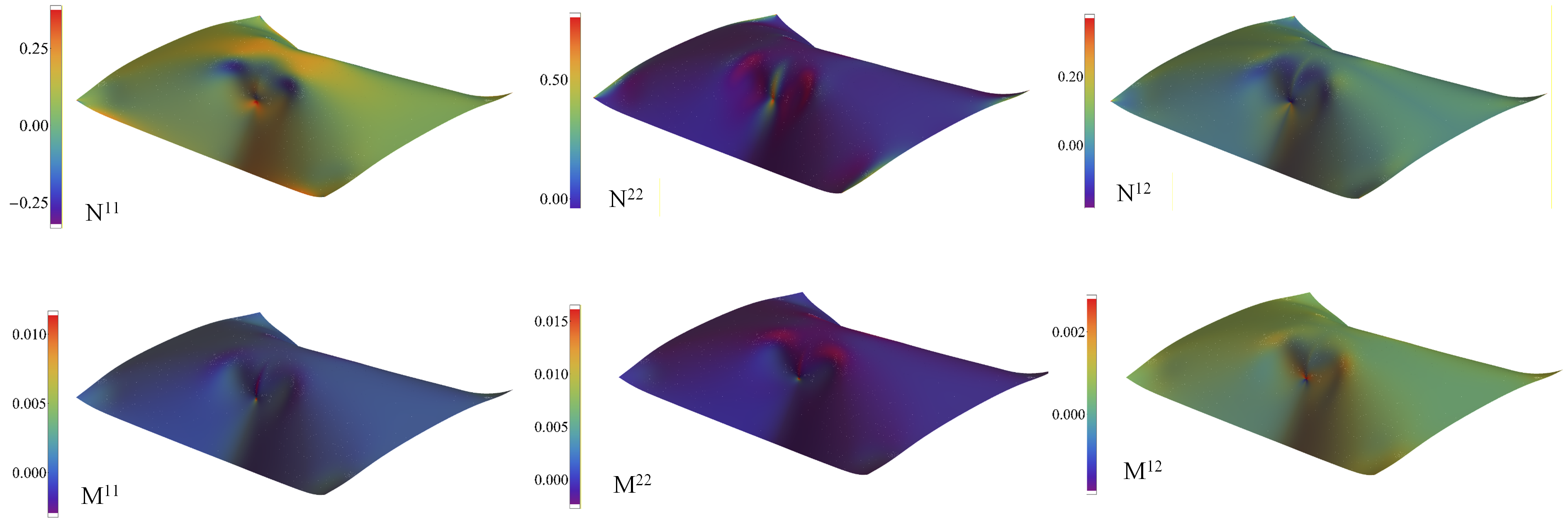} 
	\caption{Shell bending by an adhering fiber. Shell section forces for $t=1$.}
	\label{fig:shellSF}
\end{figure}

Next, the shell's section forces in the final configuration are shown in Fig.~\ref{fig:shellSF}.
Due to the strong contact forces exerted by the beam's endpoints, extreme values for the shell's section forces are located in the vicinity of those points. Furthermore, beam leaves a visible trace on these plots.

%
%
%
%
%
%

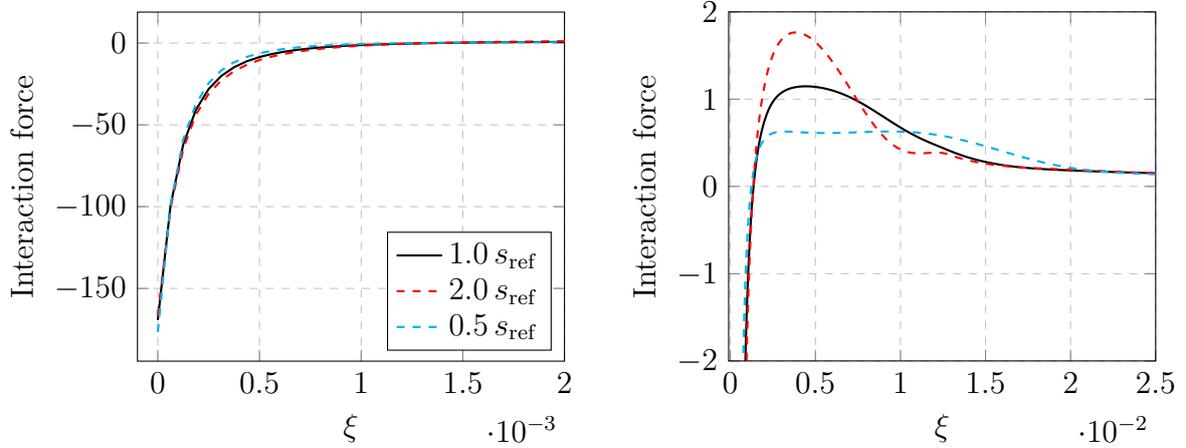
\begin{figure}[h!]
			\centering
			\subfloat[Distribution of the interaction force for $\xi \in (0,0.002)$.]{\label{fig:a}
				\begin{tikzpicture}
					\begin{axis}[
						xlabel = {$\xi$},
						ylabel = {Interaction force
						},
						legend pos= south east,
						legend cell align=left,
						width=0.45\textwidth,
						xmin=-0.0001,xmax=0.002,
						clip=true,
						grid=both]
	
						\addplot[black,thick] table [col sep=comma] {data/ipFDistrT1.csv};
						\addplot[red,dashed,thick] table [col sep=comma] {data/ipFDistrT1ip20.csv};
						\addplot[cyan,dashed,thick] table [col sep=comma] {data/ipFDistrT1ip5.csv};
						\legend{$1.0  \, s_\mathrm{ref}$,$2.0 \, s_\mathrm{ref}$,$0.5\, s_\mathrm{ref}$} 
					\end{axis}
				\end{tikzpicture}
			} \;
			\subfloat[Distribution of the interaction force for $\xi \in (0,0.025)$ with clipped peak values.]{\label{fig:b}
				\begin{tikzpicture}
					\begin{axis}[
						xlabel = {$\xi$},
						ylabel = {Interaction force
						},
						legend pos= south east,
						legend cell align=left,
						width=0.45\textwidth,
						xmin=-0.0001,xmax=0.025,
						ymin = -2, ymax = 2,
						clip=true,
						grid=both]
	
					\addplot[black,thick] table [col sep=comma] {data/ipFDistrT1.csv};
				\addplot[red,dashed,thick] table [col sep=comma] {data/ipFDistrT1ip20.csv};
				\addplot[cyan,dashed,thick] table [col sep=comma] {data/ipFDistrT1ip5.csv};
					\end{axis}
				\end{tikzpicture}
			}
			\caption{Shell bending by an adhering fiber. Distribution of the interaction force on the beam for three interaction strengths.}
			\label{fig:distIFcomp}
\end{figure}

Finally, we examine the effect of interaction strength, $s_i$, on the structural response. As in the first example, in addition to the reference interaction value  $s_\mathrm{ref}$ defined in Fig.~\ref{fig:ex2setup}, we consider simulations with $2 \, s_\mathrm{ref}$ and $0.5 \, s_\mathrm{ref}$. The final configurations for the considered interaction strengths are indistinguishable, and we focus on the interaction force. Distributions of the interaction force at the final configuration for three levels of interaction strength are shown in Fig.~\ref{fig:distIFcomp}.
Regarding the peak repulsive values near the beam's endpoints, the lowest interaction strength, $0.5 \, s_\mathrm{ref}$, produces the highest repulsive interaction force, and vice versa. This is balanced by a different distribution of the repulsive force in this area, see Fig.~\ref{fig:distIFcomp}a. On the other hand, the peak of the attractive force is proportional to the interaction strength, see Fig.~\ref{fig:distIFcomp}b. As the interaction strength reduces, the attractive force gets smeared over the larger area.

Let us emphasize that there is small frictionless sliding between the beam and the shell in this example. Namely, due to bending, the beam's surface at the contact interface contracts, while the shell's surface at the contact interface extends. To provide a resistance to this relative tangential motion between interacting bodies, we plan to introduce a tangential contact potential in our further studies.

%
%
%
%
%
%
%
%

\section{Conclusions}
\label{seccon}

This paper presents the first computational formulation for potential-based interactions between deformable beams and shells. The approach is founded on coarse-graining the interaction potential and on the analytical pre-integration of the disk-surrogate plate interaction, which approximates the disk-shell interaction. In this way, the original 6D integral of the interaction potential is reduced to a 1D integral.

To enable rigorous error estimation, we derived an analytical disk-sphere vdW law. The obtained upper bound of the approximation error for small separations is within the reasonable limits, considering the problem at hand. Thorough numerical analysis shows that the developed beam-shell formulation is capable of modeling fundamental cases of interaction between deformable fibers and membranes. 

The main advantages of the proposed approach are efficiency and good approximation of the interaction force acting on the beam. The main disadvantage is a rather poor approximation of the interaction force on the shell, since it is approximated here with a set of point forces. This issue could be addressed, for example, by employing the two half-pass approach \cite{2013sauera}, in which the forces on the shell could be calculated from a fiber-surrogate cylinder interaction.

The present work has potential applications in several fields of engineering, for example, in the study of cellular membranes and their interactions with various objects, such as fibers and proteins. Future research will focus on the tangential component of the interaction force and accompanying energy dissipation.

\section*{Acknowledgments}


This research was funded by the Austrian Science Fund (FWF) \href{https://doi.org/10.55776/P36019}{10.55776/P36019}. For the purpose of open access, the authors have applied a CC BY public copyright licence to any Author Accepted Manuscript version arising from this submission.

\section*{Appendix}

\appendix

\section{Linearization of the variation of shell strain energy}

\label{appendixa}

\setcounter{equation}{0}
\renewcommand\theequation{A\arabic{equation}}

To simplify the writing, in this Appendix, we adopt $\iv{u}{}{}=\ivn{u}{}{S}$. The linear increment of the variation of the shell's strain energy, \eqqref{shellvp}, is
\begin{equation}
	\Delta \delta \Pi_\mathrm{S,str} =  \int_{A}^{} \big(h C^{\alpha\beta\nu\gamma} \Delta \ii{E}{}{\nu\gamma} \delta E_{\alpha\beta} + \frac{h^3}{12} C^{\alpha\beta\nu\gamma} \Delta \ii{K}{}{\nu\gamma}
	\delta K_{\alpha\beta}\big) J_\text{S} \dd{\theta^1} \dd{\theta^2},
\end{equation}
where the linear increment of the variation of the membrane strain is
\begin{equation}
	\Delta\delta E_{\alpha\beta} = \frac{1}{2} \left(\Delta \iv{u}{}{,\alpha} \cdot \delta\iv{u}{}{,\beta} + \Delta \iv{u}{}{,\beta} \cdot \delta\iv{u}{}{,\alpha}\right).
\end{equation}
The increment of the bending strain is more involved
\begin{equation}
	\label{eq:dinck}
	\Delta \delta K_{\alpha\beta} = \Delta [\iv{n}{}{} \cdot \left(\delta \iv{u}{}{,\alpha,\beta} - \ii{\Gamma}{\gamma}{\alpha\beta} \cdot \delta \iv{u}{}{,\gamma}\right)] = \Delta\iv{n}{}{} \cdot \left(\delta \iv{u}{}{,\alpha,\beta} - \ii{\Gamma}{\gamma}{\alpha\beta} \cdot \delta \iv{u}{}{,\gamma}\right)+ \Delta  \ii{\Gamma}{\gamma}{\alpha\beta} \iv{n}{}{} \cdot \delta \iv{u}{}{,\gamma},
\end{equation}
since we need to find the increment of the Christoffel symbol. Increments of the contravariant base vectors follow from two conditions,
\begin{equation}
	\begin{aligned}
		\iv{g}{\alpha}{}\cdot\iv{n}{}{} = 0	 \quad \text{and} \quad  \iv{g}{\alpha}{} \cdot \iv{g}{}{\beta} = \delta^\alpha_\beta.
	\end{aligned}
\end{equation}
The first condition gives us the component along the normal, $(\Delta \iv{g}{\alpha}{})_3$, while the second condition gives us the components along tangent vectors, $(\Delta \iv{g}{\alpha}{})_\beta$, i.e.,
\begin{equation}
	\begin{aligned}
		\Delta (\iv{g}{\alpha}{}\cdot\iv{n}{}{}) &= \Delta \iv{g}{\alpha}{}\cdot\iv{n}{}{}+\iv{g}{\alpha}{}\cdot\Delta\iv{n}{}{} =0 \rightarrow (\Delta \iv{g}{\alpha}{})_3 = g^{\alpha\beta} (\iv{n}{}{} \otimes \ve{n}) \cdot \Delta \iv{u}{}{,\beta}, \\
		\Delta	(\iv{g}{\alpha}{} \cdot \iv{g}{}{\beta} )&= \Delta \iv{g}{\alpha}{} \cdot \iv{g}{}{\beta} + \iv{g}{\alpha}{} \cdot \Delta \iv{g}{}{\beta} = 0 \rightarrow 	(\Delta \iv{g}{\alpha}{})_\beta = - (\iv{g}{\beta}{} \otimes \iv{g}{\alpha}{}) \cdot \Delta \iv{u}{}{,\beta}.
	\end{aligned}
\end{equation}
Therefore, the total linear increment of contravariant basis vectors is
\begin{equation}
	\label{eq:incgrec}
	\Delta \iv{g}{\alpha}{} = [g^{\alpha\beta} (\iv{n}{}{} \otimes \ve{n})  - (\iv{g}{\beta}{} \otimes \iv{g}{\alpha}{})] \cdot \Delta \iv{u}{}{,\beta},
\end{equation}
and the increment of Christoffel symbol follows
\begin{equation}
	\label{eq:dchris}
	\Delta \ii{\Gamma}{\gamma}{\alpha\beta} = \Delta (\iv{g}{}{\alpha,\beta} \cdot \iv{g}{\gamma}{}) = \Delta\iv{g}{}{\alpha,\beta} \cdot \iv{g}{\gamma}{}+ \iv{g}{}{\alpha,\beta} \cdot \Delta\iv{g}{\gamma}{} = \iv{g}{\gamma}{} \cdot \Delta \iv{u}{}{,\alpha\beta} + [\ii{b}{}{\alpha\beta} \,g^{\gamma\delta} \iv{n}{}{} - \ii{\Gamma}{\delta}{\alpha\beta} \iv{g}{\gamma}{}] \cdot \Delta \iv{u}{}{\delta}.
\end{equation}
By inserting Eqs.~\eqref{eq:incgrec} and \eqref{eq:dchris} into \eqref{eq:dinck}, the increment of the variation of the bending strain is
\begin{equation}
	\begin{aligned}
		\Delta \delta K_{\alpha\beta} &= -\delta\iv{u}{}{,\alpha\beta} \cdot (\iv{g}{\delta}{}\otimes\ve{n})\cdot \Delta\iv{u}{}{,\delta} - \delta \iv{u}{}{,\delta} (\ve{n}\otimes\iv{g}{\delta}{}) \cdot\Delta\iv{u}{}{,\alpha\beta} - b_{\alpha\beta} \, g^{\gamma\delta} \delta \iv{u}{}{,\gamma} \cdot (\ve{n}\otimes\ve{n}) \cdot \Delta \iv{u}{}{,\delta}  \\
		&+\Gamma_{\alpha\beta}^\gamma \, \delta\iv{u}{}{,\gamma} \cdot (\iv{g}{\delta}{}\otimes\ve{n}) \cdot \Delta \iv{u}{}{,\delta}  +\Gamma_{\alpha\beta}^\delta \, \delta\iv{u}{}{,\gamma} \cdot (\iv{n}{}{}\otimes\iv{g}{\gamma}{}) \cdot \Delta \iv{u}{}{,\delta}.
	\end{aligned}
\end{equation}
These expressions agree with those obtained in \cite{2017sauerf,2023gfrerer}.

\section{Linearization of the variation of beam potential energy}

\label{appendixb}

\setcounter{equation}{0}
\renewcommand\theequation{B\arabic{equation}}

Linear increments of variations of the BE beam strains and external work due to the point moment are derived in this Appendix. Regarding the axial strain of the beam axis, its linear increment is 
\begin{equation}
	\Delta\delta \epsilon_{11} = \delta \iv{a}{}{1} \cdot \Delta \iv{a}{}{1} = \delta \iv{u}{}{\text{B},1} \cdot \Delta \iv{u}{}{\text{B},1}.
\end{equation}
For curvatures, the expressions are more involved. We omit lengthy derivations and present only final expressions. Let us introduce a vector of unknowns $\ve{q} := [ \ve{u}_{\text{B},1} , \ve{u}_{\text{B},11} ,  \varphi]$, so that the linear increments are
\begin{equation}
	\begin{aligned}
		\Delta \delta \kappa_i &= \delta \ve{q} \cdot \ve{C}^i \cdot  \Delta \ve{q}, \quad i=1,2,3.
	\end{aligned}
\end{equation}
Non-zero blocks of matrices $\ve{C}^i=[\ve{c}^i_{jk}]$ are given by the following expressions:
\begin{equation}
	\begin{aligned}
		\ve{c}^1_{11} &:= -\frac{1}{g}\left[\Gamma \, \ve{T}_\mathrm{\times} + K_2(\iv{a}{}{1}\otimes\iv{a}{}{2} + 2 \iv{a}{}{2}\otimes\iv{a}{}{1}) + K_3 (\iv{a}{}{1}\otimes\iv{a}{}{3} + 2\iv{a}{}{3}\otimes\iv{a}{}{1})\right] \\
		&+\frac{1}{a \,t_T} \left[2\Gamma \,\iv{T}{}{\times}-\iv{T}{}{\times,1} -4\Gamma (\ve{T}\times\ve{t})\otimes\ve{t}\right] \\
		&+ \frac{1}{a \, t_T^2} \left[(\ve{t} \cdot \ve{T})_{,1} \iv{T}{}{\times} + (\iv{T}{}{}\times\ve{t})\otimes\left(\iv{T}{}{,1} +(2+\ve{T} \cdot \ve{t})\iv{t}{}{,1} - 3(\ve{t}\cdot\ve{T})_{,1} \ve{t} - 2 \Gamma \, (\ve{T}\cdot\iv{t}{}{I})\right) \right. \\
		&\left. + (\ve{T}\times\ve{t})_{,1} \otimes (2\ve{t} + (\ve{T}\cdot\ve{t})\ve{t}+\ve{T})\right] \\
		&-\frac{2}{a \, t_T^3 } (\ve{T}\times\ve{t}) \otimes (\ve{T} \cdot \ve{t}_I), \\
		\ve{c}^1_{21} &:= \frac{(\ve{T}\times\ve{t})\otimes(\ve{T} \cdot \ve{t}_I)}{a \, t_T^2} + \frac{2(\ve{T}\times\ve{t})\otimes \ve{t} - \ve{T}_\times}{a \, t_T}, \\
		\ve{c}^1_{12} &:= \frac{(\ve{T}\times\ve{t})\otimes(\ve{T} \cdot \ve{t}_I)}{a \, t_T^2} + \frac{2(\ve{T}\times\ve{t})\otimes \ve{t} - \ve{T}_\times}{a \, t_T} + \frac{1}{a} \ve{T}_\times ,
	\end{aligned}
\end{equation}

\begin{equation}
	\begin{aligned}
		\ve{c}^2_{11} & := \frac{\Gamma}{\sqrt{a}} \left(\frac{\ve{a}_2\otimes(\ve{T}\times\ve{t})+(\ve{T}\times\ve{t})\otimes\ve{a}_2}{t_T}-\ve{a}_3\otimes\ve{t} -\ve{t}\otimes\ve{a}_3 \right) \\
		& - K_2 \,\ve{a}_3\otimes\ve{a}_3 -\frac{K_2 \, (\ve{T}\times\ve{t})\otimes(\ve{T}\times\ve{t})}{t_T^2} \\
		&+ K_3\, \ve{a}_3\otimes\ve{a}_2  + \frac{K_3}{t_T} \left[2 (\ve{T}\times\ve{t})\otimes\ve{t}-\ve{T}_\times+\frac{(\ve{T}\times\ve{t})\otimes(\ve{T}\cdot\ve{t}_I)}{t_T}\right], \\
		\ve{c}^2_{12} &:= \frac{1}{\sqrt{a}} \left(\ve{a}_3\otimes\ve{t} - \frac{(\ve{T}\times\ve{t})\otimes\ve{a}_2}{t_T}\right), \\
		\ve{c}^2_{21} &:= \frac{1}{\sqrt{a}} \left(\ve{t}\otimes\ve{a}_3 - \frac{\ve{a}_2\otimes(\ve{T}\times\ve{t})}{t_T}\right), \\
		\ve{c}^2_{13} &: = \ve{c}^3_{31}:=- \Gamma \, \ve{a}_2 + \frac{K_2 (\ve{T}\times\ve{t})}{\sqrt{a} \, t_T}, \\
		\ve{c}^2_{23} &:=\ve{c}^3_{32}:= \iv{a}{}{2}, \\
		\ve{c}^2_{33} &:= -K_2,
	\end{aligned}
\end{equation}
\begin{equation}
	\label{eq:link3}
	\begin{aligned}
		\ve{c}^3_{11} &:= \frac{\Gamma}{\sqrt{a}} \left(\ve{a}_2\otimes\ve{t} +\ve{t}\otimes\ve{a}_2 + \frac{\ve{a}_3\otimes(\ve{T}\times\ve{t})+(\ve{T}\times\ve{t})\otimes\ve{a}_3}{t_T}\right) \\
		&+ K_2\, \ve{a}_2\otimes\ve{a}_3  + \frac{K_2}{t_T} \left[\ve{T}_\times - 2 (\ve{T}\times\ve{t})\otimes\ve{t}-\frac{(\ve{T}\times\ve{t})\otimes(\ve{T}\cdot\ve{t}_I)}{t_T}\right] \\
		& - K_3 \,\ve{a}_2\otimes\ve{a}_2 -\frac{K_3 \, (\ve{T}\times\ve{t})\otimes(\ve{T}\times\ve{t})}{t_T^2}, \\
		\ve{c}^3_{12} &:= -\frac{1}{\sqrt{a}} \left(\ve{a}_2\otimes\ve{t} + \frac{(\ve{T}\times\ve{t})\otimes\ve{a}_3}{t_T}\right), \\
		\ve{c}^3_{21} &:= -\frac{1}{\sqrt{a}} \left(\ve{t}\otimes\ve{a}_2 + \frac{\ve{a}_3\otimes(\ve{T}\times\ve{t})}{t_T}\right), \\
		\ve{c}^3_{13} &:= \ve{c}^3_{31}:=- \Gamma \, \ve{a}_3 + \frac{K_3 (\ve{T}\times\ve{t})}{\sqrt{a}\, t_T}, \\
		\ve{c}^3_{23} &:=\ve{c}^3_{32}:= \iv{a}{}{3}, \\
		\ve{c}^3_{33} &:= -K_3,
	\end{aligned}
\end{equation}
where we have introduced designation $\ve{t}_I := \ve{I} - \ve{t} \otimes\ve{t}$. 

It is interesting to note that the symmetry of all terms of the obtained tangent operator is not evident in these expressions. However, our numerical investigations show that the resulting tangent operator is indeed symmetric. Let us investigate this analytically by focusing on the term $\ve{c}_{11}^3$ in \eqqref{eq:link3}, and its part that is not evidently symmetric:
\begin{equation}
			\tilde{\ve{c}}^3_{11} = \ve{a}_2\otimes\ve{a}_3  + \frac{1}{t_T} \left[\ve{T}_\times - 2 (\ve{T}\times\ve{t})\otimes\ve{t}-\frac{(\ve{T}\times\ve{t})\otimes(\ve{T}\cdot\ve{t}_I)}{t_T}\right].
\end{equation}
By defining arbitrary vectors, $\ve{p} = p_1 \ve{t} + p_2 \ve{a}_2 + p_3 \ve{a}_3$ and $\ve{m} = m_1 \ve{t} + m_2 \ve{a}_2 + m_3 \ve{a}_3$, the symmetry requirement is
\begin{equation}
	\label{eq:sym}
	\ve{p} \cdot \tilde{\ve{c}}^3_{11} \cdot \ve{m} = \ve{m} \cdot \tilde{\ve{c}}^3_{11} \cdot \ve{p}.
\end{equation}
With a bit of calculation, it proves that the expression \eqref{eq:sym} is satisfied and that the term $c_{11}^3$ is symmetric. The same conclusion holds for the complete tangent operator.

Linear increment of the variation of external work, defined by \eqqref{eq:varmom}, is
\begin{equation}
	\begin{aligned}
	\Delta \delta \Pi_\mathrm{ext} &=  \Delta \left[(\ve{M}\cdot\ve{t}) \delta \varphi -(\ve{M}\cdot\ve{t}) \frac{\ve{T} \times \ve{t}}{\sqrt{a}\; t_T} \cdot \delta \iv{u}{}{\text{B},1} + \frac{\ve{M} \times \ve{t} }{\sqrt{a}} \cdot \delta \iv{u}{}{\text{B},1}\right] \\
		&= \delta \varphi \; \ve{l}_1 \cdot \Delta \iv{u}{}{\text{B},1}+\delta \iv{u}{}{\text{B},1} \cdot \ve{l}_2 \cdot \Delta \iv{u}{}{\text{B},1} +\delta \iv{u}{}{\text{B},1} \cdot \ve{l}_3 \cdot \Delta \iv{u}{}{\text{B},1},
	\end{aligned}
\end{equation}
where
\begin{equation}
	\begin{aligned}
		\ve{l}_1&:=\frac{\ve{M}\cdot \ve{t}_I}{\sqrt{a}},\\
		\ve{l}_2&:=\frac{1}{a \; t_T } \left[(\ve{M}\cdot\ve{t}) \left(2 (\ve{T}\times\ve{t})\otimes\ve{t} -\ve{T}_\times + \frac{(\ve{T}\times\ve{t})\otimes (\ve{t}_I \cdot\ve{T})}{t_T}\right) -(\ve{T}\times\ve{t})\otimes (\ve{M} \cdot \ve{t}_I) \right],\\
		\ve{l}_3&:=\frac{1}{a} (\ve{M}_\times - 2 (\ve{M}\times \ve{t}) \otimes \ve{t}.\\
	\end{aligned}
\end{equation}

\section{Variation of the orthogonality condition}

\label{appendixc}

\setcounter{equation}{0}
\renewcommand\theequation{C\arabic{equation}}

Let us find the variation (that is analogous to the linear increment) of $\hat{\theta}^\alpha$ by finding the variation of the orthogonality condition 
\begin{equation}
	\label{eq:der3r3}
	\delta (\veh{ d}\cdot \veh{g}_\alpha)=\delta h_\alpha =\frac{\partial h_\alpha}{\partial \iv{r}{}{\text{B}}} \cdot \delta\iv{u}{}{\text{B}} +  \frac{\partial h_\alpha}{\partial \ivh{r}{}{\text{S}}} \cdot \delta\ivh{u}{}{\text{S}}  +\frac{\partial h_\alpha}{\partial \ivh{g}{}{\alpha}} \cdot \delta\ivh{u}{}{\mathrm{S},\alpha}  + \frac{\partial h_\alpha}{\partial \hat{\theta}^\beta} \; \delta \hat{\theta}^\beta = 0,
\end{equation}
which gives
\begin{equation}
	\label{eq:der32v3}
	\iv{\hat g}{}{\alpha} \cdot \delta\iv{u}{}{\text{B}} -\ivh{ g}{}{\alpha} \cdot   \delta\ivh{u}{}{\text{S}} +\veh{d} \cdot \delta \ivh{ g}{}{\mathrm{S},\alpha} -\hat C_{\alpha \beta} \; \delta \hat{\theta}^\beta = 0,
\end{equation}
where we have introduced a symmetric shift tensor $\hat C_{\alpha \beta}$ defined by
\begin{equation}
	\label{eq:derortvho}
	\frac{\partial h_\alpha}{\partial \hat{\theta}^\beta}= - \hat{\ve{g}}_\beta \cdot \hat{\ve{g}}_\alpha + \veh{ d} \cdot \hat{\ve{g}}_{\alpha,\beta} = -\hat g_{\alpha \beta} +d \; \hat \Gamma^3_{\alpha \beta} = -(\hat{g}_{\alpha \beta } -d \; \hat b_{\alpha \beta}) =: -\hat{C}_{\alpha \beta} = -(\hat{C}^{\alpha \beta})^{-1}.
\end{equation}
Now, the variation of the parametric coordinates is
\begin{equation}
	\label{eq:derv3s23}
	\delta \hat{\theta}^\beta = \hat C^{\alpha \beta} [(\delta\iv{u}{}{\text{B}} -   \delta\iv{\hat u}{}{\text{S}} ) \cdot \iv{\hat g}{}{\alpha} +\ve{\hat d} \cdot \delta \iv{\hat u}{}{\mathrm{S},\alpha} ] =(\delta\iv{u}{}{\text{B}} -   \delta\ivh{u}{}{\text{S}} ) \cdot \bar{\hat{\ve{g}}}^\beta +\hat C^{\alpha \beta} \veh{ d} \cdot \delta \ivh{u}{}{\mathrm{S},\alpha},
\end{equation}
where quantity $\bar{\hat{\ve{g}}}^\beta = \hat C^{ \alpha \beta} \iv{\hat g}{}{\alpha}$
represents contravariant basis vectors of a shell shifted to the beam axis. Importantly, these vectors are orthogonal to the normal vector of the shell and the distance vector.

\section{Linearization of variation for the RF2 formulation}

\label{appendixd}

\setcounter{equation}{0}
\renewcommand\theequation{D\arabic{equation}}

Let us first linearize the RF2 expressions. Within this reduced formulation, we completely disregard the influence of the angle. Such an approach corresponds to the linearization in standard contact mechanics, where the interaction only depends on the normal gap between two contact surfaces \cite{2006wriggers}. The linear increment of the virtual work at one beam integration point consists of two terms defined in \eqqref{eq:var1aa}:
%
%
\begin{equation}
	\begin{aligned}
		\Delta (\iv{f}{}{} \cdot \delta \iv{u}{}{\text{B}})&=\Delta \ve{f} \cdot \delta \iv{u}{}{\text{B}},\\
		\Delta (\iv{f}{}{} \cdot \delta \ivh{u}{}{\text{S}})&=\Delta \ve{f} \cdot \delta \ivh{u}{}{\text{S}} +\ve{f} \cdot \Delta \delta \ivh{u}{}{\text{S}} =\Delta \ve{f} \cdot \delta \ivh{u}{}{\text{S}} +\ve{f} \cdot \delta \ivh{u}{}{\mathrm{S},\alpha} \; \Delta \theta^\alpha.
	\end{aligned}
\end{equation}
The increment of variation of the displacement of the shell exists, since this vector implicitly depends on the closest-point projection and therefore, on the parametric coordinate
\begin{equation}
	\begin{aligned}
		\Delta \delta \veh{u}_\mathrm{S} = \frac{\partial (\delta \veh{u}_\mathrm{S})}{\partial \hat{\theta}^\alpha} \Delta \hat{\theta}^\alpha = \delta \frac{\partial \veh{u}_\mathrm{S}}{\partial \hat{\theta}^\alpha} \Delta \hat{\theta}^\alpha = \delta \ivh{u}{}{\mathrm{S},\alpha} \; \Delta \hat{\theta}^\alpha.
	\end{aligned}
\end{equation}
The linear increment of the interaction force is
\begin{equation}
	\begin{aligned}
		\Delta \ve{f}&= \Delta (\phi_{,d} \ve{\hat n})=  \phi_{,dd} \ve{\hat n} \Delta d +  \phi_{,d} \Delta \ve{\hat n}.
	\end{aligned}
\end{equation}
Let us now derive and sort all required increments:
\begin{equation}
	\label{eq:incnorm}
	\begin{aligned}
		\Delta \veh{d} &=  \Delta \iv{u}{}{\text{B}} - \Delta \ivh{u}{}{\text{S}} -\ivh{g}{}{\alpha} \; \Delta \hat{\theta}^\alpha  = \Delta \iv{u}{}{\text{B}} -\Delta \ivh{u}{}{\text{S}} - \ivh{g}{}{\alpha}[ (\Delta\iv{u}{}{\text{B}} -   \Delta\ivh{u}{}{\text{S}} ) \cdot \bar{\hat{\ve{g}}}^\alpha +\hat C^{\alpha\beta} \veh{d}\cdot \Delta \ivh{u}{}{\mathrm{S},\beta}] \\
		&= (\ve{I} - \ivh{g}{}{\alpha}\otimes \bar{\hat{\ve{g}}}^\alpha) (\Delta \iv{u}{}{\text{B}} -\Delta \ivh{u}{}{\text{S}})  - (\bar{\hat{\ve{g}}}^\beta \otimes \veh{d}) \cdot \Delta \ivh{u}{}{\text{S},\beta},\\
		\Delta d &=  \veh{n} \cdot \Delta \veh{d} = \veh{n} \cdot (\Delta \iv{u}{}{\text{B}} - \Delta \ivh{u}{}{\text{S}} ),\\
		\Delta \veh{n} &=  \Delta  \frac{\veh{d}}{d} = \frac{d \Delta \veh{d} -  \Delta d \veh{d}}{d^2} = \frac{1}{d} (\Delta\veh{d} -  \Delta  d \veh{n})	\\
		&= \frac{1}{d} [ (\ve{I} - \ivh{g}{}{\alpha}\otimes \bar{\hat{\ve{g}}}^\alpha) (\Delta \iv{u}{}{\text{B}} -\Delta \ivh{u}{}{\text{S}})  - (\bar{\hat{\ve{g}}}^\beta \otimes \veh{d}) \cdot \Delta \ivh{u}{}{\text{S},\beta}  -  (\veh{n} \otimes \veh{n}) (\Delta \iv{u}{}{\text{B}} -\Delta \ivh{u}{}{\text{S}})  ] \\
		&= \frac{1}{d} (\ve{I} - \ivh{g}{}{\alpha}\otimes \bar{\hat{\ve{g}}}^\alpha - \veh{n} \otimes \veh{n}) (\Delta \iv{u}{}{\text{B}} -\Delta \ivh{u}{}{\text{S}}) - (\bar{\hat{\ve{g}}}^\beta \otimes \veh{n}) \cdot \Delta \ivh{u}{}{\text{S},\beta}.  \\
	\end{aligned}
\end{equation}
Contrary to the standard approach in contact mechanics, where the increment of normal is obtained from the orthogonality condition \cite{2006wriggers}, here we derive it from the distance vector.

With these expressions at hand, the increment of variation w.r.t.~the beam is
\begin{equation}
	\begin{aligned}
		\Delta \iv{f}{}{} \cdot \delta \iv{u}{}{\text{B}}&=  \delta \iv{u}{}{\text{B}} \cdot \big[ \phi_{,dd} \; \veh{n} \Delta d +  \phi_{,d} \Delta \veh{n}   \big]  \\
		&=  \delta \iv{u}{}{\text{B}} \cdot \{ \phi_{,dd} \; \veh{n} (\veh{n} \cdot \Delta \iv{u}{}{\text{B}} -\veh{n} \cdot \Delta \ivh{u}{}{\text{S}} )\\
		&+ \phi_{,d} [\frac{1}{d} (\ve{I} - \ivh{g}{}{\alpha}\otimes \bar{\hat{\ve{g}}}^\alpha - \veh{n} \otimes \veh{n}) (\Delta \iv{u}{}{\text{B}} -\Delta \ivh{u}{}{\text{S}}) - (\bar{\hat{\ve{g}}}^\beta \otimes \veh{n}) \cdot \Delta \ivh{u}{}{\text{S},\beta}]  \} \\
		&= \delta \iv{u}{}{\text{B}} \cdot ( \ve{k}_{11} \cdot \Delta \iv{u}{}{\text{B}} +\ve{k}_{12}\cdot \Delta \ivh{u}{}{\text{S}} +\ve{k}_{13} \cdot \Delta \ivh{u}{}{\text{S},\beta}  ),
	\end{aligned}
\end{equation}
and the increment of variation w.r.t.~the shell is
\begin{equation}
	\begin{aligned}
		&\Delta \ve{f} \cdot \delta \iv{\hat u}{}{\text{S}} +\ve{f} \cdot \delta \iv{\hat u}{}{\text{S},\alpha} \; \Delta \theta^\alpha  \\
		&=\delta \iv{\hat u}{}{\text{S}} \cdot \big[ \phi_{,dd} \; \ve{\hat n} \Delta d +  \phi_{,d} \Delta \ve{\hat n}   \big]  +\phi_{,d} \ve{\hat n} \cdot \delta \iv{\hat u}{}{\text{S},\alpha} \; [(\Delta\iv{u}{}{\text{B}} -   \Delta\iv{\hat u}{}{\text{S}} ) \cdot \bar{\hat{\ve{g}}}^\alpha +\hat C^{\alpha \beta} \ve{\hat d} \cdot \Delta \iv{\hat u}{}{\text{S},\beta} ]  \\
		&= \delta \iv{\hat u}{}{\text{S}} \cdot \{ \phi_{,dd} \; \ve{\hat n} (\ve{\hat n} \cdot \Delta \iv{u}{}{\text{B}} -\ve{\hat n} \cdot \Delta \iv{\hat u}{}{\text{S}} )\\
		&+  \phi_{,d} [\frac{1}{d} (\ve{I} - \iv{\hat g}{}{\alpha}\otimes \bar{\hat{\ve{g}}}^\alpha - \ve{\hat n} \otimes \ve{\hat n}) (\Delta \iv{u}{}{\text{B}} -\Delta \iv{\hat u}{}{\text{S}}) - (\bar{\hat{\ve{g}}}^\beta \otimes \ve{\hat n}) \cdot \Delta \iv{\hat u}{}{\text{S},\beta}]  \} \\
		&+  \delta \iv{\hat u}{}{\text{S},\alpha} \cdot \{ \phi_{,d} [(\ve{\hat n} \otimes \bar{\hat{\ve{g}}}^\alpha)\cdot (\Delta\iv{u}{}{\text{B}} -   \Delta\iv{\hat u}{}{\text{S}} ) +  (\ve{\hat n} \otimes\ve{\hat d})\hat C^{\alpha \beta} ]\cdot \Delta \iv{\hat u}{}{\mathrm{S},\beta} \} \\
		&= \delta \iv{\hat u}{}{\text{S}} \cdot (\ve{k}_{21} \cdot \Delta \iv{u}{}{\text{B}}+\ve{k}_{22}\cdot \Delta \iv{\hat u}{}{\text{S}} + \ve{k}_{23} \cdot \Delta \iv{\hat u}{}{\text{S},\beta}  ) \\
		&+\delta \iv{\hat u}{}{\text{S},\alpha} \cdot (  \ve{k}_{31} \cdot \Delta\iv{u}{}{\text{B}} + \ve{k}_{32} \cdot \Delta\iv{\hat u}{}{\text{S}} + \ve{k}_{33} \cdot \Delta \iv{\hat u}{}{\mathrm{S},\beta} ),
	\end{aligned}
\end{equation}
where we have introduced the following notation:
\begin{equation}
	\label{eq:linvar3}
	\begin{aligned}
		\ve{k}_{11} &:= \phi_{,dd} \; (\veh{n} \otimes \veh{n}) + \phi_{,d} \frac{1}{d} (\ve{I} - \ivh{g}{}{\alpha}\otimes \bar{\hat{\ve{g}}}^\alpha - \veh{n} \otimes \veh{n}), \\
		\ve{k}_{12} &:= -\ve{k}_{11}, \\
		\ve{k}_{13} &:= -\phi_{,d}  (\bar{\hat{\ve{g}}}^\beta \otimes \veh{n}), \\
		\ve{k}_{21} &:= -\ve{k}_{11}, \\
		\ve{k}_{22} &:= \ve{k}_{11}, \\
		\ve{k}_{23} &:= -\ve{k}_{13}  \\
		\ve{k}_{31} &:= -\phi_{,d} (\veh{n} \otimes \bar{\hat{\ve{g}}}^\alpha), \\
		\ve{k}_{32} &:= -\ve{k}_{31},\\
		\ve{k}_{33} &:= -d(\veh{n} \otimes\veh{n}) \hat C^{\alpha \beta}.
	\end{aligned}
\end{equation}
If we define the vector of generalized unknowns $\ve{u}_\alpha := [\ve{u}_\text{B} \; \ivh{u}{}{\text{S}} \; \ivh{ u}{}{\mathrm{S},\alpha}]$, then the linear increment of variation of potential can be represented as
\begin{equation}
	\label{eq:linvar2}
	\begin{aligned}
		\Delta \delta \phi &= 	\Delta (\iv{f}{}{} \cdot \delta \iv{u}{}{\text{B}}) -\Delta (\iv{f}{}{} \cdot \delta \ivh{u}{}{\text{S}})= \delta \ve{u}_\alpha \cdot \ve{K}^\text{RF2}_\text{T} \cdot \Delta \ve{u}_\beta,
	\end{aligned}
\end{equation}
where $\ve{K}^\text{RF2}_\text{T}$ is the tangent operator for the RF2 formulation with blocks defined in \eqqref{eq:linvar3}. The obtained result is symmetric w.r.t.~variation and linearization, as expected. 

%



\section{Linearization of variation for the RF1 formulation}

\label{appendixe}

\setcounter{equation}{0}
\renewcommand\theequation{E\arabic{equation}}

In the RF1 formulation, we disregard the variation w.r.t~the angle as in RF2, but consider the potential as a function of both the distance and the angle, see \eqqref{eq:var1aa}. Variations of beam and shell positions do not depend on the angle between the beam's cross section and the surrogate plate. Therefore, in addition to the terms derived for the linearization of the RF2 in Appendix \ref{appendixd}, we need the linear increment of the interaction force w.r.t.~the cosine, i.e.,
\begin{equation}
	\begin{aligned}
		\Delta_{\cos \alpha} \ve{f} &= \Delta_{\cos \alpha} (\phi_{,d} \hat{\ve n}) =\phi_{,d \cos \alpha} \Delta \cos \alpha = - \phi_{,d \cos \alpha} \frac{t_n}{\cos \alpha}  (\hat{\ve{n} }\cdot \Delta \ve{t} + \ve{t} \cdot \Delta \hat{\ve{n}}) \hat{\ve{n}} \\
		&= w_d [(\hat{\ve{n} } \otimes\hat{\ve{n} })  \Delta \ve{t} + (\hat{\ve{n} } \otimes \ve{t} ) \Delta \hat{\ve{n}}] , \quad \text{with} \quad w_d:= - \phi_{,d \cos \alpha} \frac{t_n}{\cos\alpha}.
	\end{aligned}
\end{equation}
To address singularity when $\cos\alpha=1$, we cancel terms in $w_d$, i.e.,
\begin{equation}
	\begin{aligned}
		w_d= - \phi^R_{,d \cos \alpha}t_n, \quad  \phi^R_{,d \cos \alpha}:= \phi_{,d \cos \alpha}/\cos\alpha.
	\end{aligned}
\end{equation}
The linear increment of the normal is given in \eqref{eq:incnorm}, while the increment of the tangent is
\begin{equation}
	\label{eq:inctn}
	\begin{aligned}
		\Delta \ve{t} &= \frac{1}{\sqrt{a}} \left(\ve{I} - \ve{t}\otimes\ve{t}\right) \Delta \ve{u}_{\mathrm{B},1}.
	\end{aligned}
\end{equation}
The linear increment of the interaction force w.r.t.~the cosine can be written as
\begin{equation}
	\begin{aligned}
		\Delta_{\cos \alpha} \ve{f} &=  \frac{w_d}{\sqrt{j_\mathrm{B}}} \left[\hat{\ve{n} } \otimes\hat{\ve{n} }  - t_n (\hat{\ve{n}} \otimes\ve{t})\right] \Delta \ve{u}_{B,1} \\
		&+ \frac{w_d}{d} [\hat{\ve{n} } \otimes \ve{t} - (\ve{t} \cdot \iv{\hat g}{}{\alpha})(\hat{\ve{n} }\otimes \bar{\hat{\ve{g}}}^\alpha) - t_n( \ve{\hat n} \otimes \ve{\hat n})] (\Delta \iv{u}{}{\text{B}} -\Delta \ivh{u}{}{\text{S}}) \\
		&- w_d (\ve{t} \cdot \bar{\hat{\ve{g}}}^\beta) ( \veh{n}\otimes \veh{n}) \cdot \Delta \ivh{u}{}{\text{S},\beta}.
	\end{aligned}
\end{equation}
If we introduce the vector of generalized unknowns $\ve{u}_\alpha^\mathrm{RF1} := [\ve{u}_\text{B} \; \ve{u}_{\text{B},1} \; \iv{\hat u}{}{\text{S}} \; \iv{\hat u}{}{S,\alpha}]$, then the linear increment of variation of potential w.r.t.~the cosine is
\begin{equation}
	\label{eq:linvar}
	\begin{aligned}
		\Delta_{\cos \alpha} \delta \phi &= 	\Delta_{\cos \alpha} \iv{f}{}{} \cdot (\delta \iv{u}{}{\text{B}} -\delta \ivh{u}{}{\text{S}})= \delta \ve{u}_\alpha^\mathrm{RF1} \cdot \ve{K}_\text{T}^\mathrm{RF1} \cdot \Delta \ve{u}_\beta^\mathrm{RF1},
	\end{aligned}
\end{equation}
where $\ve{K}_\text{T}^\mathrm{RF1}$ is the tangent operator with the following blocks:
\begin{equation}
	\label{eq:linvars}
	\begin{aligned}
		\veq{k}_{11} &:=  \frac{w_d}{d} [\hat{\ve{n} } \otimes \ve{t} - (\ve{t} \cdot \iv{\hat g}{}{\alpha})(\hat{\ve{n} }\otimes \bar{\hat{\ve{g}}}^\alpha) - t_n( \ve{\hat n} \otimes \ve{\hat n})], \\
		\veq{k}_{12} &:= \frac{w_d}{\sqrt{j}} \left[\hat{\ve{n} } \otimes\hat{\ve{n} }  - t_n (\hat{\ve{n}} \otimes\ve{t})\right], \\
		\veq{k}_{13} &:= -\veq{k}_{11},\\
		\veq{k}_{14} &:=- w_d (\ve{t} \cdot \bar{\hat{\ve{g}}}^\beta) ( \ve{\hat n}\otimes \ve{\hat n}),  \\
		\veq{k}_{2i} &:= \veq{k}_{4i}:=0, \quad i=1,2,3,4, \\
		\veq{k}_{31} &:= -\veq{k}_{11},\\
		\veq{k}_{32} &:= -\veq{k}_{12},\\
		\veq{k}_{33} &:= \veq{k}_{11}, \\
		\veq{k}_{34} &:= -\veq{k}_{14}.
	\end{aligned}
\end{equation}
The complete tangent operator of the RF1 formulation is $\ve{K}_\text{T}^\mathrm{RF1}+\ve{K}_\text{T}^\mathrm{RF2}$, c.f.~Eqs.~\eqref{eq:linvar2} and \eqref{eq:linvar3}. The tangent operator for the RF1 formulation is not symmetric because the variation of the potential w.r.t.~the angle has been disregarded.

Regarding the linearization of the FF formulation, numerical simulations show that the RF1 tangent provides good convergence for the FF formulation. Since the linearization of FF is very involved, and the results are indistinguishable from those of the RF1 (see Section \ref{sec:ex1}), we omit its linearization.

\section{Derivation of the disk-sphere interaction potential}

\label{appendixf}

\setcounter{equation}{0}
\renewcommand\theequation{F\arabic{equation}}

\begin{figure}[hpt]
	\centering
	\includegraphics[width=0.6\textwidth]{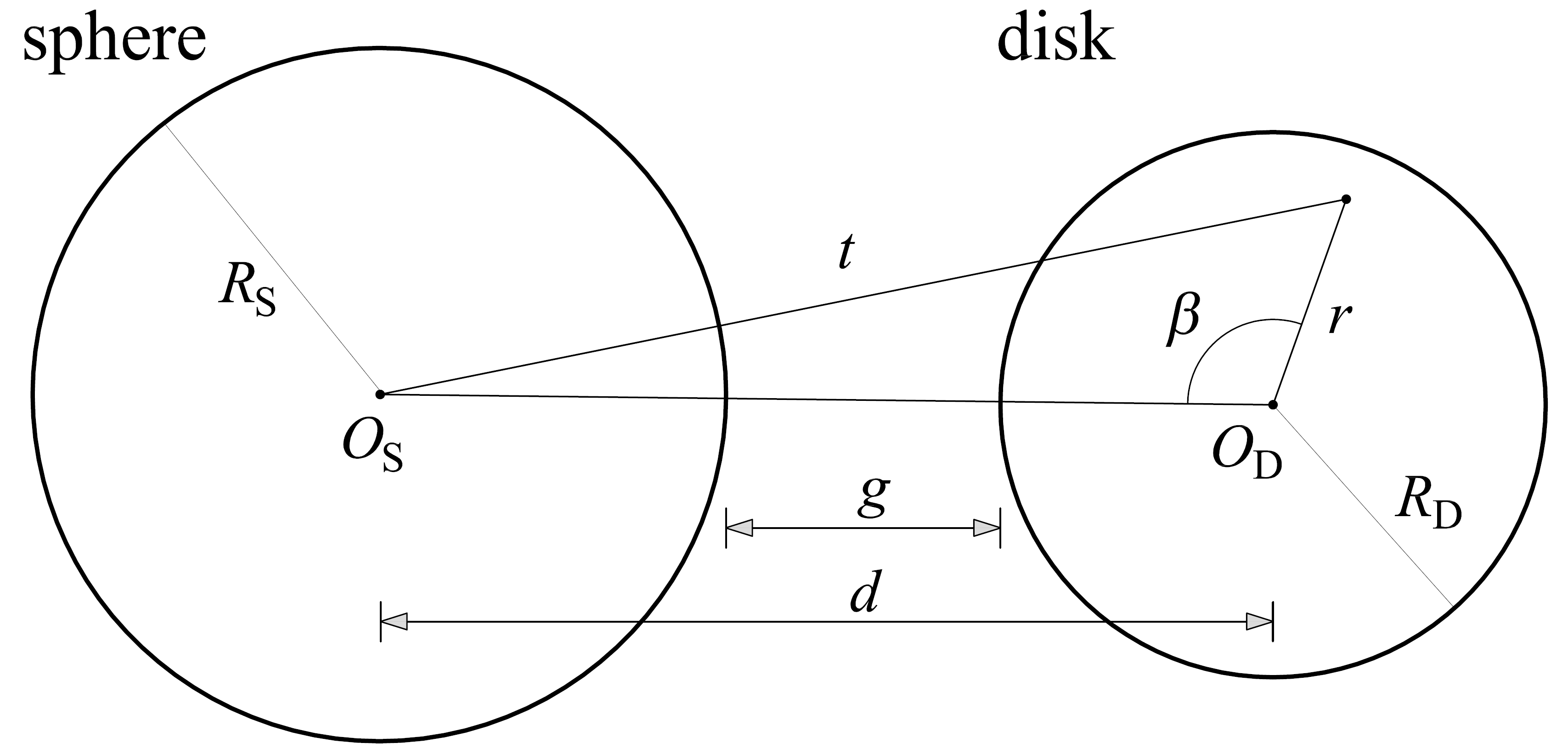} 
	\caption{Interaction between a disk and a sphere in a perpendicular orientation.}
	\label{fig:DS}
\end{figure}

In this Appendix, we derive an analytical vdW law for the disk-sphere interaction, which was used in Section \ref{sec:error2} for error estimation. We start from the point-sphere vdW law \cite{2003kirsch},
\begin{equation}
	\Pi_\mathrm{P-S}(t) = \frac{4 \pi R_\mathrm{S}}{3 t^3 (t+2R_\mathrm{S})},
\end{equation}
where $t$ is the distance between the sphere's center and the point. A schematic of the disk-sphere interaction is shown in Fig.~\ref{fig:DS}, illustrating the special case in which the disk is perpendicular to the tangential plane of the sphere at the closest point. 
In other words, the disk lies in the plane passing through the sphere's center. Using the designations defined in Fig.~\ref{fig:DS}, we can express the distance between the sphere's center and an arbitrary point of the disk, $t$, as
\begin{equation}
	t = \sqrt{d^2+r^2-2dr\cos\beta}, \quad \text{with} \quad g=d-R_\mathrm{S} - R_\mathrm{D},
\end{equation}
where $g$ is the gap between the sphere and the disk. By integrating the point-sphere potential over the area of the disk, we obtain the disk-sphere vdW interaction potential
\begin{equation}
	\Pi_\mathrm{D-S} =2 \int_{0}^{\pi}  \int_{0}^{R_\mathrm{D}} 	\Pi_\mathrm{P-S} (t) \, r \dd{r} \dd{\beta},
\end{equation}
which has the analytical solution
\begin{equation}
\Pi_\mathrm{D-S} = \frac{\pi^2}{3 R_\mathrm{S}} \left(1+\frac{R_\mathrm{S}^6+(R_\mathrm{D}^2 - 3 d^2)R_\mathrm{S}^4 + 3 (d^4 - R_\mathrm{D}^4 ) R_\mathrm{S}^2 - (d^2 - R_\mathrm{D}^2 )^3}{\left\{\left[d^2-(R_\mathrm{D} + R_\mathrm{S})\right]^2 \left[d^2+(R_\mathrm{D} + R_\mathrm{S})\right]^2\right\}^{3/2}}\right).
\end{equation}
To the best of our knowledge, this is the first derivation of an analytical disk-sphere vdW law. Arguably, this law could be used as a second-order approximation for short-range vdW beam-shell interactions, where the section-shell interaction could be approximated with a section-spherical shell interaction.

\printbibliography

@article{1994perdomo,
  title = {Interaction between cells and elastin fibers: {{An}} ultrastructural and immunocytochemical study},
  shorttitle = {Interaction between cells and elastin fibers},
  author = {Perdomo, J.J. and Gounon, P. and Schaeverbeke, M. and Schaeverbeke, J. and Groult, V. and Jacob, M.P. and Robert, L.},
  date = {1994},
  journaltitle = {J Cell Physiol},
  volume = {158},
  number = {3},
  pages = {451--458},
  doi = {10.1002/jcp.1041580309},
  langid = {english}
}

@book{1995lipowsky,
  title = {Structure and dynamics of membranes: {{I}}. {{From}} cells to vesicles / {{II}}. {{Generic}} and specific interactions},
  shorttitle = {Structure and {{Dynamics}} of {{Membranes}}},
  author = {Lipowsky, R. and Sackmann, E.},
  date = {1995-06-15},
  eprint = {V1H9oVpNC7wC},
  eprinttype = {googlebooks},
  publisher = {Elsevier},
  langid = {english},
  pagetotal = {537}
}

@article{1997argento,
  title = {Surface formulation for molecular interactions of macroscopic bodies},
  author = {Argento, C. and Jagota, A. and Carter, W. C.},
  date = {1997-07-01},
  journaltitle = {J Mech Phys Solids},
  volume = {45},
  number = {7},
  pages = {1161--1183},
  doi = {10.1016/S0022-5096(96)00121-4},
  langid = {english}
}

@article{2003kirsch,
  title = {Calculation of the van der {{Waals}} force between a spherical particle and an infinite cylinder},
  author = {Kirsch, V. A},
  date = {2003-07-01},
  journaltitle = {Adv. Colloid Interface Sci.},
  volume = {104},
  number = {1},
  pages = {311--324},
  doi = {10.1016/S0001-8686(03)00053-8}
}

@article{2005delrio,
  title = {The role of van der {{Waals}} forces in adhesion of micromachined surfaces},
  author = {DelRio, Frank W. and family=Boer, given=Maarten P., prefix=de, useprefix=true and Knapp, James A. and David Reedy, E. and Clews, Peggy J. and Dunn, Martin L.},
  date = {2005-08},
  journaltitle = {Nature Mater},
  volume = {4},
  number = {8},
  pages = {629--634},
  publisher = {Nature Publishing Group},
  doi = {10.1038/nmat1431},
  langid = {english}
}

@article{2005hughes,
  title = {Isogeometric analysis: {{CAD}}, finite elements, {{NURBS}}, exact geometry and mesh refinement},
  shorttitle = {Isogeometric analysis},
  author = {Hughes, T. J. R. and Cottrell, J. A. and Bazilevs, Y.},
  date = {2005-10-01},
  journaltitle = {Comput Methods Appl Mech Eng},
  volume = {194},
  number = {39},
  pages = {4135--4195},
  doi = {10.1016/j.cma.2004.10.008},
  langid = {english}
}

@book{2005parsegian,
  title = {Van der {{Waals Forces}}: {{A Handbook}} for {{Biologists}}, {{Chemists}}, {{Engineers}}, and {{Physicists}}},
  shorttitle = {Van der {{Waals Forces}}},
  author = {Parsegian, V. Adrian},
  date = {2005},
  publisher = {Cambridge University Press},
  location = {Cambridge},
  doi = {10.1017/CBO9780511614606}
}

@book{2006wriggers,
  title = {Computational {{Contact Mechanics}}},
  author = {Wriggers, P.},
  date = {2006},
  publisher = {Springer-Verlag},
  location = {Berlin Heidelberg},
  doi = {10.1007/978-3-540-32609-0},
  langid = {english}
}

@article{2007sauer,
  title = {A contact mechanics model for quasi-continua},
  author = {Sauer, Roger A. and Li, S.},
  date = {2007},
  journaltitle = {Int. J. Numer. Methods Eng.},
  volume = {71},
  number = {8},
  pages = {931--962},
  doi = {10.1002/nme.1970},
  langid = {english}
}

@article{2008alavinasab,
  title = {Computational modeling of nano-structured glass fibers},
  author = {Alavinasab, A. and Jha, R. and Ahmadi, G. and Cetinkaya, C. and Sokolov, I.},
  date = {2008-12-01},
  journaltitle = {Comput Mater Sci},
  volume = {44},
  number = {2},
  pages = {622--627},
  doi = {10.1016/j.commatsci.2008.05.004},
  langid = {english}
}

@article{2008sauer,
  title = {An atomistically enriched continuum model for nanoscale contact mechanics and its application to contact scaling},
  author = {Sauer, Roger A. and Li, S.},
  date = {2008-07},
  journaltitle = {J. Nanosci. Nanotechnol.},
  volume = {8},
  number = {7},
  eprint = {19051933},
  eprinttype = {pubmed},
  pages = {3757--3773},
  doi = {PMID: 19051933},
  langid = {english}
}

@article{2010saarikangas,
  title = {Regulation of the actin cytoskeleton-plasma membrane interplay by phosphoinositides},
  author = {Saarikangas, J. and Zhao, H. and Lappalainen, P.},
  date = {2010-01-01},
  journaltitle = {Physiol. Rev.},
  volume = {90},
  number = {1},
  pages = {259--289},
  publisher = {American Physiological Society},
  doi = {10.1152/physrev.00036.2009}
}

@incollection{2010wriggers,
  title = {Contact between {{Beams}} and {{Shells}}},
  booktitle = {New {{Trends}} in {{Thin Structures}}: {{Formulation}}, {{Optimization}} and {{Coupled Problems}}},
  author = {Wriggers, Peter},
  editor = {De Mattos Pimenta, Paulo and Wriggers, Peter},
  date = {2010},
  pages = {155--174},
  publisher = {Springer},
  location = {Vienna},
  doi = {10.1007/978-3-7091-0231-2_6},
  langid = {english}
}

@book{2011israelachvili,
  title = {Intermolecular and surface forces - 3rd {{Edition}}},
  author = {Israelachvili, J. N.},
  date = {2011},
  publisher = {Academic Press},
  url = {https://www.elsevier.com/books/intermolecular-and-surface-forces/israelachvili/978-0-12-391927-4},
  urldate = {2021-08-03}
}

@article{2012murrell,
  title = {F-actin buckling coordinates contractility and severing in a biomimetic actomyosin cortex},
  author = {Murrell, M. P. and Gardel, M. L.},
  date = {2012-12-18},
  journaltitle = {PNAS},
  volume = {109},
  number = {51},
  eprint = {23213249},
  eprinttype = {pubmed},
  pages = {20820--20825},
  publisher = {National Academy of Sciences},
  doi = {10.1073/pnas.1214753109},
  langid = {english}
}

@article{2012yoo,
  title = {Scalable fabrication of silicon nanotubes and their application to energy storage},
  author = {Yoo, J.-K. and Kim, J. and Jung, Y.S. and Kang, K.},
  date = {2012},
  journaltitle = {Adv Mater},
  volume = {24},
  number = {40},
  pages = {5452--5456},
  doi = {10.1002/adma.201201601}
}

@article{2013sauera,
  title = {A computational contact formulation based on surface potentials},
  author = {Sauer, Roger A. and De Lorenzis, L.},
  date = {2013-01-01},
  journaltitle = {Comput Methods Appl Mech Eng},
  volume = {253},
  pages = {369--395},
  doi = {10.1016/j.cma.2012.09.002},
  langid = {english}
}

@article{2014sauerb,
  title = {A geometrically exact finite beam element formulation for thin film adhesion and debonding},
  author = {Sauer, Roger A. and Mergel, Janine C.},
  date = {2014-09-01},
  journaltitle = {Finite Elem. Anal. Des.},
  volume = {86},
  pages = {120--135},
  doi = {10.1016/j.finel.2014.03.009},
  langid = {english}
}

@article{2015meier,
  title = {A locking-free finite element formulation and reduced models for geometrically exact {{Kirchhoff}} rods},
  author = {Meier, Christoph and Popp, Alexander and Wall, Wolfgang A.},
  date = {2015-06-15},
  journaltitle = {Computer Methods in Applied Mechanics and Engineering},
  volume = {290},
  pages = {314--341},
  doi = {10.1016/j.cma.2015.02.029},
  langid = {english}
}

@article{2015walani,
  title = {Endocytic proteins drive vesicle growth via instability in high membrane tension environment},
  author = {Walani, N. and Torres, J. and Agrawal, A.},
  date = {2015-03-24},
  journaltitle = {PNAS},
  volume = {112},
  number = {12},
  eprint = {25775509},
  eprinttype = {pubmed},
  pages = {E1423-E1432},
  publisher = {National Academy of Sciences},
  doi = {10.1073/pnas.1418491112},
  langid = {english}
}

@article{2016fan,
  title = {A three-dimensional surface stress tensor formulation for simulation of adhesive contact in finite deformation},
  author = {Fan, H. and Li, S.},
  date = {2016},
  journaltitle = {Int. J. Numer. Methods Eng.},
  volume = {107},
  number = {3},
  pages = {252--270},
  doi = {10.1002/nme.5169},
  langid = {english}
}

@article{2016wang,
  title = {Electrospun nanofiber membranes},
  author = {Wang, Xuefen and Hsiao, Benjamin S},
  date = {2016-05-01},
  journaltitle = {Current Opinion in Chemical Engineering},
  series = {Nanotechnology / {{Separation Engineering}}},
  volume = {12},
  pages = {62--81},
  doi = {10.1016/j.coche.2016.03.001}
}

@article{2017ambrosetti,
  title = {Physical adsorption at the nanoscale: {{Towards}} controllable scaling of the substrate-adsorbate van der {{Waals}} interaction},
  shorttitle = {Physical adsorption at the nanoscale},
  author = {Ambrosetti, Alberto and Silvestrelli, Pier Luigi and Tkatchenko, Alexandre},
  date = {2017-06-12},
  journaltitle = {Phys. Rev. B},
  volume = {95},
  number = {23},
  pages = {235417},
  publisher = {American Physical Society},
  doi = {10.1103/PhysRevB.95.235417}
}

@article{2017islama,
  title = {Morphology and mechanics of fungal mycelium},
  author = {Islam, M. R. and Tudryn, G. and Bucinell, R. and Schadler, L. and Picu, R. C.},
  date = {2017-10-12},
  journaltitle = {Sci Rep},
  volume = {7},
  number = {1},
  pages = {13070},
  publisher = {Nature Publishing Group},
  doi = {10.1038/s41598-017-13295-2},
  issue = {1},
  langid = {english}
}

@article{2017lu,
  title = {A new design for an artificial cell: polymer microcapsules with addressable inner compartments that can harbor biomolecules, colloids or microbial species},
  shorttitle = {A new design for an artificial cell},
  author = {Lu, A.X. and Oh, H. and Terrell, J.L. and Bentley, W.E. and Raghavan, S.R.},
  date = {2017-09-25},
  journaltitle = {Chem Sci},
  volume = {8},
  number = {10},
  pages = {6893--6903},
  publisher = {The Royal Society of Chemistry},
  doi = {10.1039/C7SC01335C},
  langid = {english}
}

@article{2017sauer,
  title = {A stabilized finite element formulation for liquid shells and its application to lipid bilayers},
  author = {Sauer, Roger A. and Duong, Thang X. and Mandadapu, Kranthi K. and Steigmann, David J.},
  date = {2017-02-01},
  journaltitle = {Journal of Computational Physics},
  volume = {330},
  pages = {436--466},
  doi = {10.1016/j.jcp.2016.11.004},
  langid = {english}
}

@article{2017sauerf,
  title = {On the theoretical foundations of thin solid and liquid shells},
  author = {Sauer, Roger A and Duong, Thang X},
  date = {2017-03-01},
  journaltitle = {Math. Mech. Solids},
  volume = {22},
  number = {3},
  pages = {343--371},
  publisher = {SAGE Publications Ltd STM},
  doi = {10.1177/1081286515594656},
  langid = {english}
}

@article{2017venkataram,
  title = {Unifying {{Microscopic}} and {{Continuum Treatments}} of van der {{Waals}} and {{Casimir Interactions}}},
  author = {Venkataram, Prashanth S. and Hermann, Jan and Tkatchenko, Alexandre and Rodriguez, Alejandro W.},
  date = {2017-06-29},
  journaltitle = {Phys. Rev. Lett.},
  volume = {118},
  number = {26},
  pages = {266802},
  publisher = {American Physical Society},
  doi = {10.1103/PhysRevLett.118.266802}
}

@article{2018borkovic,
  title = {Rotation-free isogeometric analysis of an arbitrarily curved plane {{Bernoulli}}–{{Euler}} beam},
  author = {Borković, A. and Kovačević, S. and Radenković, G. and Milovanović, S. and Guzijan-Dilber, M.},
  date = {2018-06-01},
  journaltitle = {Computer Methods in Applied Mechanics and Engineering},
  volume = {334},
  pages = {238--267},
  doi = {10.1016/j.cma.2018.02.002},
  langid = {english}
}

@article{2018brely,
  title = {The influence of substrate roughness, patterning, curvature, and compliance in peeling problems},
  author = {Brely, Lucas and Bosia, Federico and Pugno, Nicola M},
  date = {2018-01},
  journaltitle = {Bioinspir. Biomim.},
  volume = {13},
  number = {2},
  pages = {026004},
  publisher = {IOP Publishing},
  doi = {10.1088/1748-3190/aaa0e5},
  langid = {english}
}

@article{2018charles-orszag,
  title = {Adhesion to nanofibers drives cell membrane remodeling through one-dimensional wetting},
  author = {Charles-Orszag, A. and Tsai, F.-C. and Bonazzi, D. and Manriquez, V. and Sachse, M. and Mallet, A. and Salles, A. and Melican, k. and Staneva, R. and Bertin, A. and Millien, C. and Goussard, S. and Lafaye, P. and Shorte, S. and Piel, M. and Krijnse-Locker, J. and Brochard-Wyart, F. and Bassereau, P. and Duménil, G.},
  date = {2018-10-25},
  journaltitle = {Nat Commun},
  volume = {9},
  number = {1},
  pages = {4450},
  publisher = {Nature Publishing Group},
  doi = {10.1038/s41467-018-06948-x},
  issue = {1},
  langid = {english}
}

@article{2018franquelim,
  title = {Membrane sculpting by curved {{DNA}} origami scaffolds},
  author = {Franquelim, H. G. and Khmelinskaia, A. and Sobczak, J.-P. and Dietz, H. and Schwille, P.},
  date = {2018-02-23},
  journaltitle = {Nat Commun},
  volume = {9},
  number = {1},
  eprint = {29476101},
  eprinttype = {pubmed},
  pages = {811},
  doi = {10.1038/s41467-018-03198-9},
  langid = {english},
  pmcid = {PMC5824810}
}

@article{2018gov,
  title = {Guided by curvature: shaping cells by coupling curved membrane proteins and cytoskeletal forces},
  shorttitle = {Guided by curvature},
  author = {Gov, N. S.},
  date = {2018-05-26},
  journaltitle = {Philos Trans R Soc Lond B Biol Sci},
  volume = {373},
  number = {1747},
  pages = {20170115},
  publisher = {Royal Society},
  doi = {10.1098/rstb.2017.0115}
}

@article{2018nishiyama,
  title = {Molecular interactions in nanocellulose assembly},
  author = {Nishiyama, Y.},
  date = {2018-02-13},
  journaltitle = {Philos Trans Royal Soc A},
  volume = {376},
  number = {2112},
  pages = {20170047},
  publisher = {Royal Society},
  doi = {10.1098/rsta.2017.0047}
}

@article{2018yuan,
  title = {Adhesion of carbon nanotubes on elastic substrates with finite thickness},
  author = {Yuan, Xuebo and Wang, Youshan},
  date = {2018-10-19},
  journaltitle = {J. Appl. Phys.},
  volume = {124},
  number = {15},
  pages = {155306},
  doi = {10.1063/1.5048240}
}

@article{2020brunelli,
  title = {Nanofiber membranes as biomimetic and mechanically stable surface coatings},
  author = {Brunelli, M. and Alther, S. and Rossi, R. M. and Ferguson, S. J. and Rottmar, M. and Fortunato, G.},
  date = {2020-03-01},
  journaltitle = {Materials Science and Engineering: C},
  volume = {108},
  pages = {110417},
  doi = {10.1016/j.msec.2019.110417}
}

@article{2020gaynetob,
  title = {Master-master frictional contact and applications for beam-shell interaction},
  author = {Gay Neto, Alfredo and Wriggers, Peter},
  date = {2020-12-01},
  journaltitle = {Comput Mech},
  volume = {66},
  number = {6},
  pages = {1213--1235},
  doi = {10.1007/s00466-020-01890-6},
  langid = {english}
}

@article{2020grill,
  title = {A computational model for molecular interactions between curved slender fibers undergoing large {{3D}} deformations with a focus on electrostatic, van der {{Waals}}, and repulsive steric forces},
  author = {Grill, Maximilian J. and Wall, Wolfgang A. and Meier, C.},
  date = {2020},
  journaltitle = {Int. J. Numer. Methods Eng.},
  volume = {121},
  number = {10},
  pages = {2285--2330},
  doi = {10.1002/nme.6309},
  langid = {english}
}

@article{2020maa,
  title = {A finite element-based coarse-grained model for cell–nanomaterial interactions by combining absolute nodal coordinate formula and {{Brownian}} dynamics},
  author = {Ma, T. and Liu, Y. and Lin, G. and Wang, C. and Tan, H.},
  date = {2020-12-07},
  journaltitle = {J Appl Mech},
  volume = {88},
  number = {4},
  doi = {10.1115/1.4049143}
}

@article{2020vo,
  title = {Geometrically nonlinear multi-patch isogeometric analysis of planar curved {{Euler}}–{{Bernoulli}} beams},
  author = {Vo, Duy and Nanakorn, Pruettha},
  date = {2020-07},
  journaltitle = {Comput. Methods Appl. Mech. Eng.},
  volume = {366},
  pages = {113078},
  doi = {10.1016/j.cma.2020.113078},
  langid = {english}
}

@article{2020zhangb,
  title = {Measurement of {{Adhesion}} of {{In Situ Electrospun Nanofibers}} on {{Different Substrates}} by a {{Direct Pulling Method}}},
  author = {Zhang, Bin and Yan, Xu and Xu, Yuan and Zhao, Huai-Song and Yu, Miao and Long, Yun-Ze},
  date = {2020},
  journaltitle = {Adv. Mater. Sci. Eng.},
  volume = {2020},
  number = {1},
  pages = {7517109},
  doi = {10.1155/2020/7517109},
  langid = {english}
}

@article{2021čanadija,
  title = {Deep learning framework for carbon nanotubes: mechanical properties and modeling strategies},
  shorttitle = {Deep learning framework for carbon nanotubes},
  author = {Čanadija, M.},
  date = {2021-09-06},
  journaltitle = {Carbon},
  doi = {10.1016/j.carbon.2021.08.091},
  langid = {english}
}

@article{2021grilla,
  title = {Investigation of the peeling and pull-off behavior of adhesive elastic fibers via a novel computational beam interaction model},
  author = {Grill, Maximilian J. and Meier, C. and Wall, Wolfgang A.},
  date = {2021-06-11},
  journaltitle = {J. Adhes.},
  volume = {97},
  number = {8},
  pages = {730--759},
  publisher = {Taylor \& Francis},
  doi = {10.1080/00218464.2019.1699795}
}

@article{2021radenković,
  title = {Nonlinear static isogeometric analysis of arbitrarily curved {{Kirchhoff-Love}} shells},
  author = {Radenković, G. and Borković, A. and Marussig, B.},
  date = {2021-02-15},
  journaltitle = {Int. J. Mech. Sci.},
  volume = {192},
  pages = {106143},
  doi = {10.1016/j.ijmecsci.2020.106143},
  langid = {english}
}

@article{2021slepukhin,
  title = {Topological defects produce kinks in biopolymer filament bundles},
  author = {Slepukhin, V. M. and Grill, Maximilian J. and Hu, Q. and Botvinick, E.L. and Wall, W.A. and Levine, A.J.},
  date = {2021-04-13},
  journaltitle = {PNAS},
  volume = {118},
  number = {15},
  eprint = {33876768},
  eprinttype = {pubmed},
  publisher = {National Academy of Sciences},
  doi = {10.1073/pnas.2024362118},
  langid = {english}
}

@article{2021tozzia,
  title = {A theory of ordering of elongated and curved proteins on membranes driven by density and curvature},
  author = {Tozzi, C. and Walani, N. and Roux, A.-L.L. and Roca-Cusachs, P. and Arroyo, M.},
  date = {2021-04-01},
  journaltitle = {Soft Matter},
  volume = {17},
  number = {12},
  pages = {3367--3379},
  publisher = {The Royal Society of Chemistry},
  doi = {10.1039/D0SM01733G},
  langid = {english}
}

@article{2022borkovićc,
  title = {Geometrically exact static isogeometric analysis of an arbitrarily curved spatial {{Bernoulli}}–{{Euler}} beam},
  author = {Borković, A. and Marussig, B. and Radenković, G.},
  date = {2022-02-15},
  journaltitle = {Comput. Methods Appl. Mech. Eng.},
  volume = {390},
  pages = {114447},
  doi = {10.1016/j.cma.2021.114447},
  langid = {english}
}

@article{2023borković,
  title = {Geometrically exact isogeometric {{Bernoulli}}–{{Euler}} beam based on the {{Frenet}}–{{Serret}} frame},
  author = {Borković, A. and Gfrerer, M. H. and Marussig, B.},
  date = {2023-02-15},
  journaltitle = {Comput. Methods Appl. Mech. Eng.},
  volume = {405},
  pages = {115848},
  doi = {10.1016/j.cma.2022.115848},
  langid = {english}
}

@article{2023gfrerer,
  title = {Rigorous code verification for non-linear {{Kirchhoff}}–{{Love}} shells based on tangential differential calculus with application to {{Isogeometric Analysis}}},
  author = {Gfrerer, M. H.},
  date = {2023-12-01},
  journaltitle = {Finite Elements in Analysis and Design},
  volume = {227},
  pages = {104041},
  doi = {10.1016/j.finel.2023.104041}
}

@article{2023grill,
  title = {Analytical disk–cylinder interaction potential laws for the computational modeling of adhesive, deformable (nano)fibers},
  author = {Grill, Maximilian J. and Wall, Wolfgang A. and Meier, Christoph},
  date = {2023-05-01},
  journaltitle = {Int. J. Solids Struct.},
  volume = {269},
  pages = {112175},
  doi = {10.1016/j.ijsolstr.2023.112175},
  langid = {english}
}

@article{2023meier,
  title = {Generalized section-section interaction potentials in the geometrically exact beam theory: {{Modeling}} of intermolecular forces, asymptotic limit as strain-energy function, and formulation of rotational constraints},
  shorttitle = {Generalized section-section interaction potentials in the geometrically exact beam theory},
  author = {Meier, Christoph and Grill, Maximilian J. and Wall, Wolfgang A.},
  date = {2023-05-08},
  journaltitle = {Int. J. Solids Struct.},
  pages = {112255},
  doi = {10.1016/j.ijsolstr.2023.112255},
  langid = {english}
}

@article{2024borkovićb,
  title = {A novel section–section potential for short-range interactions between plane beams},
  author = {Borković, A. and Gfrerer, M. H. and Sauer, R. A. and Marussig, B. and Bui, T. Q.},
  date = {2024-09-01},
  journaltitle = {Comput. Methods Appl. Mech. Eng.},
  volume = {429},
  pages = {117143},
  doi = {10.1016/j.cma.2024.117143}
}

@inproceedings{2024borkovićg,
  title = {A note on beam-to-beam contact dynamics},
  booktitle = {Int. {{Conf}}. {{Contemp}}. {{Theory Pract}}. {{Constr}}.},
  author = {Borković, Aleksandar and Jočković, Miloš and Tatar, Dijana and Milovanović, Snježana},
  date = {2024-06-12},
  volume = {16},
  number = {1},
  pages = {337--350},
  publisher = {{University of Banja Luka, Faculty of Architecture, Civil Engineering and Geodesy}},
  doi = {10.61892/stp202401080B},
  langid = {english}
}

@article{2024grill,
  title = {Asymptotically consistent and computationally efficient modeling of short-ranged molecular interactions between curved slender fibers undergoing large {{3D}} deformations},
  author = {Grill, Maximilian J. and Wall, Wolfgang A. and Meier, Christoph},
  date = {2024-04-15},
  journaltitle = {Adv. Model. Simul. Eng. Sci.},
  volume = {11},
  number = {1},
  pages = {7},
  doi = {10.1186/s40323-023-00257-9}
}

@article{2024mokhalingam,
  title = {Continuum contact model for friction between graphene sheets that accounts for surface anisotropy and curvature},
  author = {Mokhalingam, Aningi and Gupta, Shakti S. and Sauer, Roger A.},
  date = {2024},
  journaltitle = {Phys. Rev. B},
  volume = {109},
  pages = {035435},
  doi = {10.1103/PhysRevB.109.035435}
}

@article{2024sky,
  title = {Intrinsic mixed-dimensional beam-shell-solid couplings in linear {{Cosserat}} continua via tangential differential calculus},
  author = {Sky, Adam and Hale, Jack S. and Zilian, Andreas and Bordas, Stéphane P. A. and Neff, Patrizio},
  date = {2024-12-01},
  journaltitle = {Computer Methods in Applied Mechanics and Engineering},
  volume = {432},
  pages = {117384},
  doi = {10.1016/j.cma.2024.117384}
}

@article{2025borković,
  title = {New analytical laws and applications of interaction potentials with a focus on van der {{Waals}} attraction},
  author = {Borković, A. and Gfrerer, M. H. and Sauer, R. A.},
  date = {2025-03-24},
  journaltitle = {Applied Mathematical Modelling},
  pages = {116100},
  doi = {10.1016/j.apm.2025.116100}
}

@software{2025borkovićb,
  title = {Supplementary notebooks for the manuscript regarding the new analytical laws and applications of interaction potentials with a focus on van der {{Waals}} attraction},
  author = {Borković, Aleksandar and Gfrerer, M. H.},
  date = {2025-03-06},
  doi = {10.3217/81zxn-1np78},
  organization = {Graz University of Technology}
}

@book{2025borkoviće,
  title = {Computational modeling of potential-based interactions between fibers},
  author = {Borković, Aleksandar},
  date = {2025},
  publisher = {{Faculty of Architecture, Civil Engineering and Geodesy, University of Banja Luka}},
  location = {Banja Luka},
  url = {https://pub.unibl.org/s/eng/item/66563}
}

@article{2025song,
  title = {Multi-material {{4D}} printing and {{3D}} patterned metallization enables smart architectures},
  author = {Song, Kewei and Xiong, Chunfeng and Zhang, Ze and Wu, Kunlin and Wan, Weiyang and Wang, Yifan and Umezu, Shinjiro and Sato, Hirotaka},
  date = {2025-04-15},
  journaltitle = {Composites Part B: Engineering},
  volume = {295},
  pages = {112218},
  doi = {10.1016/j.compositesb.2025.112218}
}

@article{2026borkovića,
  title = {Efficient snap-to-contact computations for van der {{Waals}} interacting fibers},
  author = {Borković, A. and Gfrerer, M. H. and Sauer, R. A. and Marussig, B.},
  date = {2026-03-01},
  journaltitle = {European Journal of Mechanics - A/Solids},
  volume = {116},
  pages = {105919},
  doi = {10.1016/j.euromechsol.2025.105919}
}

\end{document}